\documentclass[a4paper,reqno]{amsart}

\usepackage{graphics,graphicx,color}
\usepackage{srcltx, bm}
\usepackage{tikz}
\usepackage{sgame}
\usepackage{amssymb,amsfonts}
\usepackage{latexsym}
\usepackage{amsmath,amsthm,color, amscd,marginnote}
\usepackage[all]{xy}
 \usepackage{mathrsfs}
\usepackage{MnSymbol}
\usepackage{hyperref}
\hypersetup{colorlinks, linkcolor=blue,citecolor=blue}

\usepackage{graphicx,psfrag,epsfig}

\theoremstyle{plain}
\newtheorem{theorem}{Theorem}

\newtheorem{lemma}[theorem]{Lemma}
\newtheorem{corollary}[theorem]{Corollary}
\newtheorem{proposition}[theorem]{Proposition}
\theoremstyle{definition}
\newtheorem{definition}[theorem]{Definition}
\newtheorem{remark}[theorem]{Remark}
\newtheorem{assumption}[theorem]{Assumption}
\newtheorem{example}[theorem]{Example}
\newtheorem{amplification}[theorem]{Amplification}

\def\eps{\varepsilon}
\def\ep{\varepsilon}
\def\epsi{\delta}

\def\bdef{\begin{definition}}
\def\endef{\end{definition}}
\def\bthm{\begin{theorem}}
\def\ethm{\end{theorem}}
\def\blm{\begin{lemma}}
\def\elm{\end{lemma}}
\def\brm{\begin{remark}}
\def\erm{\end{remark}}
\def\bprop{\begin{proposition}}
\def\eprop{\end{proposition}}
\def\bcor{\begin{corollary}}
\def\ecor{\end{corollary}}
\def\be{\begin{eqnarray}}
\def\ee{\end{eqnarray}}
\def\beal{\begin{aligned}}
\def\enal{\end{aligned}}

\def\om{\omega}
\def\Om{\Omega}

\def\eps{\varepsilon}

\def\phi{\varphi}
\def\f{\varphi}

\def\R{\mathbb R}

\def\Nn{\mathbb N}
\def\T{\mathbb T}

\def\Z{\mathbb Z}

\def\cC{\mathcal C}
\def\cP{\mathcal P}

\def\cF{\mathcal F}
\def\~{\tilde}

\def\fB{\mathfrak B}

\def\cB{\mathcal B}
\def\cC{\mathcal C}

\def\cD{\mathcal D}

\def\cS{\mathcal S}

\def\PP{\mathbf{P}}
\def\EE{\mathbf{E}}

\def\p{\partial}
 \newcommand{\strela}{\rightharpoonup}
 
  \def\BR{\bar B_R}
  
  \def\hve{\hat v^\eps }
   \def\vst{v^{st} }
  \def\vep{v^\eps }

  \def\Iet{I^\eps_\tau}
   \def\Ist{I^{st}_\tau}

  \def\tilI{\tilde I^\varepsilon_\tau}
   \def\tilst{\tilde I^{st}_\tau}

\def\lan{\langle}
\def\ran{\rangle}

\def\llan{\langle\!\langle}
\def\rran{\rangle\!\rangle}

\def\be{\begin{equation}}
\def\ee{\end{equation}}
\def\bdef{\begin{definition}}
\def\endef{\end{definition}}
\def\blm{\begin{lemma}}
\def\elm{\end{lemma}}
\def\beal{\begin{aligned}}
\def\enal{\end{aligned}}

\newtheorem*{Pf}{Proof}
\renewenvironment{proof}{\begin{Pf} \begin{upshape}} {\end{upshape} \qed\end{Pf}}

\def\bv{\mathbf{v}}
\def\bbeta{\boldsymbol\beta}

\def\bV{\mathbf{V}}
\def\bP{\mathbf{P}}
\def\bQ{\mathbf{Q}}
\def\bE{\mathbf{E}}
\def\bR{\mathbf{R}}
\def\bA{\mathbf{A}}
\def\b1{\mathbf{1}}

\title[Averaging for stochastic perturbations  ]{Averaging for stochastic perturbations of  integrable  systems}
\begin{document}
\author{Guan Huang}
\address{School of Mathematics and Statistics, Beijing Institute of Technology, Beijing, China and  Peoples' Friendship University of Russia (RUDN University),  Moscow,  Russia}
\email{huangguan@bit.edu.cn}

\author{Sergei Kuksin}
\address%{Institut de Math\'emathiques de Jussieu--Paris Rive Gauche, CNRS, Universit\'e Paris Diderot, UMR 7586, Sorbonne Paris Cit\'e, F-75013, Paris, France \& School of Mathematics, Shandong University, Jinan, Shandong,  China}
{Universit\'e  Paris Cit\'e and Sorbonne Universit\'e, CNRS, IMJ-PRG, F-75013 Paris, France, and
%{Institut de Math\'emathiques de Jussieu--Paris Rive Gauche, CNRS, Universit\'e Paris Diderot, UMR 7586, Sorbonne Paris Cit\'e, Paris, France, and
 Peoples' Friendship University of Russia,  %(RUDN University),
 Moscow,  Russia, and Steklov Mathematical Institute of
  Russian Academy of Sciences, Moscow, Russia}
\email{sergei.kuksin@imj-prg.fr}

\author{Andrey Piatnitski}
\address{The Arctic University of Norway, campus Narvik, Norway,
and
Institute for Information Transmission Problems of RAS, Moscow, Russia, and  Peoples' Friendship University of Russia (RUDN University),  Moscow,  Russia}
\email{apiatnitski@gmail.com}

\begin{abstract}
We are concerned with averaging theorems for $\eps$-small stochastic perturbations of integrable equations in
$\R^d \times \T^n =\{(I,\f)\}$
$$
\dot I(t) =0, \quad \dot \f(t) = \theta(I),
  \eqno{ (1)}
$$
and in $\R^{2n} = \{ v=(\bv_1, \dots, \bv_n), \ \bv_j \in \R^2\}$,
$$
\dot\bv_k(t) =W_k(I) \bv_k^\perp, \quad k=1, \dots, n,
  \eqno{ (2)}
$$
where $I=(I_1, \dots, I_n)$ is the vector of actions, $I_j = \frac12 \| \bv_j\|^2$.
The vector-functions $\theta$ and $W$
 are locally Lipschitz and non-degenerate. Perturbations of these
  equations are assumed to be locally Lipschitz and
 such that some few first  moments of the norms of
  their  solutions are bounded uniformly in $\eps$, for $0\le t\le \eps^{-1} T$. For $I$-components
 of solutions for perturbations of (1) we establish their convergence in law to solutions of the corresponding
 averaged  $I$-equations, when $0\le \tau := \eps t\le T$ and $\eps\to0$.
 Then we  show that if the system of averaged $I$-equations
 is mixing, then the convergence is uniform in the slow time $\tau=\eps t\ge0$.

 Next using these results, for   $\eps$-perturbed
  equations (2) we construct  well posed {\it effective stochastic equations} for $v(\tau)\in \R^{2n}$ (independent of $\eps$)
  such that when $\eps\to0$, the actions of solutions for the perturbed equations   with $t:=  \tau/\eps$
   converge in distribution to actions
  of solutions for the effective equations. Again, if the effective system is mixing, this convergence is
  uniform in the slow time $\tau  \ge0$.

  We provide easy sufficient conditions on the perturbed equations which ensure that our results apply to their solutions.
 \end{abstract}

\maketitle
\numberwithin{equation}{section}
\numberwithin{theorem}{section}
%\tableofcontents
%

\section{Introduction}

\subsection
{Setting and problems.}
The goal of this paper is to present an averaging theory for   stochastic differential equations, obtained by
 diffusive perturbations of integrable  deterministic   differential equations. All equations in our work have locally
 Lipschitz coefficients.

  In previous publication
 \cite{GHKu}, written by two of us, an easier problem of stochastic perturbations of linear systems with imaginary
 spectrum was considered. This work may be regarded as a natural continuation of  \cite{GHKu}.
  Both papers are based on the Khasminski approach to the averaging in stochastic systems and    use some
 technical ideas, developed by the authors and A.~Maiocchi in their work on stochastic PDEs
 \cite{KP, K, HKM, HGK22}.
  We consider two classes of problems as above.    Firstly,
   given an  integrable system in $  \R^d \times \T^n$,
   \footnote{ If $n=d$ and $\theta(I) = \nabla h(I)$ for some $C^1$-function $h$ on $\R^d$,
   then \eqref{int_syst} is an integrable Hamiltonian system on the symplectic space
   $(\R^n \times \T^n, dI \wedge d\f)$
  with the Hamiltonian function $h$. Otherwise \eqref{int_syst}  is integrable in some general sense.   }
   \be\label{int_syst}
   \dot I =0, \quad \dot \f = \theta(I), \qquad (I,\f) \in \R^d \times \T^n, \quad \T^n = \R^n/ 2\pi \Z^n,
   \ee
   we study the asymptotic properties, as $\varepsilon\to0$,  of solutions for
    perturbed systems of the form
   \begin{equation}\label{2}
\begin{split}
 d I&=\eps P^I(I ,\varphi )dt+ \sqrt{{\ep}}\,\Psi^I(I ,\varphi ) d\beta(t), %\quad I\in \R^d,
 \\
d\varphi &=\big[\theta(I )+\eps P^\varphi(I ,\varphi )\big]dt+\sqrt{{\ep}}\,\Psi^\varphi(I ,\varphi )d\beta(t),
% \quad \phi\in \T^n,
\end{split}
\end{equation}
where  $\beta(t)$ is the standard Wiener process in $\mathbb{R}^{{d_1}}$ and
the matrices $\Psi^I$ and $\Psi^\varphi$ have  corresponding dimensions.
We also consider stochastic
 perturbations of  integrable equations
  in $\R^{2n}$\,\footnote{  If $W(I)= \nabla h(I)$ for some $C^1$-function $h$, then \eqref{Birk} is a Hamiltonian equation in $\R^{2n}$,
  integrable in the sense of Birkhoff.
    }
 \be\label{Birk}
 \dot{\mathbf{v}}_k=W_k(I) \mathbf{v}_k^\perp, \qquad k=1,2,\ldots,n,
 \ee
where $\mathbf{v}_k=(v_k,v_{-k})^t\in\mathbb R^2$, $ \mathbf{v}_k^\perp= (-v_{-k}, v_{k})^t$
 and $I_k$ is the $k$-th action
 $\frac12\|\mathbf{v}_k\|^2$.  Similar to the above  we are interested in  stochastic perturbations of  the form
\begin{equation}\label{3}
\begin{split}
\displaystyle
d\mathbf{v}_k &=\big[ W_k(I)\mathbf{v}_k^\perp+
\eps {\mathbf{P}}_k(v)\big]dt+\sqrt{{\ep}}\,\sum\limits_{j=1}^{n_1}B_{kj}(v)d\beta_j(t), \quad k=1, \dots, n.
\end{split}
\end{equation}
For convenience of language, below we call systems of equations like \eqref{Birk} and \eqref{3} just ``equations".

  Our goal is to study the asymptotic behaviour of  solutions for systems  \eqref{2} and   \eqref{3}
   as $\eps\to0$ on time-interval of length $\sim\eps^{-1}$, and the limiting behaviour of the solutions for $0\le t\le\infty$
    when $\eps\to0$, provided that suitable  mixing assumptions hold.
 To study this behaviour it
 is convenient to pass to the slow time $\tau = \eps t$, and this is what we do below in the main text. In the
 introduction we discuss our results using both the original fast time $t$ and slow time $\tau$.

\subsection{The results.}
Systems \eqref{2} are examined in  the first part of our work (Sections~\ref{s_ps}-\ref{s_lp}).
The first main result of this analysis is
Theorem~\ref{t_main02} which  describes the statistical
behaviour of components  $I^\eps(t)$ of solutions for   \eqref{2} on time-intervals of order  $\eps^{-1}$.
There we assume that  the mapping $\theta$ is  nondegenerate in sense of Anosov (i.e.  for almost all $I\in\R^d$
components of the vector $\theta(I)$ are rationally independent),
that the diffusion in the $I$-equation in \eqref{2} is bounded non-degenerate, that the coefficients of  equation \eqref{2}
satisfy some mild regularity assumptions and that certain a-priori estimates for solutions of \eqref{2} hold uniformly in $\eps$.
The theorem states that the vector of  $I$-components $I^\eps(\eps^{-1}\tau)$
 of a solution  for $0\le \tau \le T$ converges in
law to a solution $I^0(\tau)$  of an averaged stochastic differential equation, obtained by means of an
appropriate averaging in $\f$ of  coefficients of  the $I$-equation in \eqref{2}.  More precisely, the theorem states
that
\be\label{conv1}
\cD I^\eps (\eps^{-1}\tau) \strela \cD I^0(\tau),
\ee
where $\cD$ signifies the distribution of a  r.v. (random variable), the arrow stands for weak convergence of
 measures and $I^0(\tau)$ is a solution of a stochastic equation
 \be\label{I_aver}
 dI(\tau) = \lan P^I\ran(I) d\tau + \llan \Psi^I\rran (I) d\beta(\tau).
 \ee
In  \eqref{I_aver}  $ \lan P^I\ran$ is the averaging of $P^I$ in angles and matrix
 $ \llan \Psi^I\rran (I)$ is obtained from $\Psi^I$ using
the rules of stochastic averaging, see in Subsection~\ref{ss_3.2}.

In Section \ref{s_stationary_nonl}, assuming the mixing, we examine the asymptotic in time
 behaviour of components  $I^\eps(\tau)$. More precisely, Proposition~\ref{p_5.4_sec10} states that if in addition to the
 above-mentioned assumptions system \eqref{2} and the averaged equation \eqref{I_aver} both are mixing with stationary
 measures $\nu^\eps$ and $\mu^0$ respectively,\footnote{Concerning the mixing in SDEs and its basic properties which
 we use see  \cite{khasminskii, GHKu} and references in \cite{GHKu}.}
  then the $I$-projection of  $\nu^\eps$ weakly converges to $\mu^0$ as $\eps\to0$. So for any solution
  $(I^\eps, \f^\eps)$ of \eqref{2} we have
  \be\label{conv2}
\lim_{\eps\to0} \lim_{\tau\to\infty} \cD I^\eps(\tau) =\mu^0.
\ee
Convergence \eqref{conv1}, valid for any finite  time $T$, jointly with $\eqref{conv2}$ suggest that,
 in fact, convergence  \eqref{conv1}
holds uniformly in   $\tau\ge0$. Later in Section~\ref{s_stationary_nonl} we show that indeed this is the case.
Namely, let $\| \mu_1- \mu_2\|_L^*$ be the dual-Lipschitz distance between measures $\mu_1, \mu_2$,
see Definition~\ref{d_dual-Lip}.\footnote{Also known as the Kantorovich--Rubinstein distance.} In
Theorem~\ref{t_mixAP}, assuming that equation \eqref{I_aver} is mixing
 with a mild quantitative property of the rate of mixing,\footnote{and without assuming that system
 \eqref{2}  is mixing.}
 we prove that
 \be\label{conv3}
\lim_{\eps\to0} \sup_{\tau\ge0}\| \cD I^\eps(\eps^{-1}\tau) - \cD I^0(\tau)\|_L^*=0.
\ee

The  behaviour of  $I$-components of solutions for  deterministically  perturbed
equation \eqref{2} with $\Psi^I =0$ and $\Psi^\f=0$  for typical initial data
on time-intervals of order $\eps^{-1}$ is the subject of the classical averaging
theory due to Anosov--Kosuga--Neishtadt, see  \cite[Sec.~6.1]{AKN} and  \cite{LM}.
Moreover, due to Neishtadt and Bakhtin the set of atypical initial data for which the averaging fails
 has measure $\lesssim\eps^\gamma$,
where $\gamma>0$ depends on ``how non-degenerate mapping $\theta$ is".
 Theorem~\ref{t_main02} is a natural counterpart
of that theory for stochastic perturbations of integrable equations. But  stochastic
convergences \eqref{conv2} and \eqref{conv3},
established in  Proposition~\ref{p_5.4_sec10} and Theorem~\ref{t_mixAP} and valid for any initial data,
seem to have no deterministic analogy.
 \smallskip

 Equation \eqref{2} is a fast-slow stochastic system with fast variable $\f$ and slow variable $I$. The averaging in
  such systems (which allows to approximatively describe the law of the slow variable $I$ on time-integrals of
  order $\eps^{-1}$)  is a well developed topic of stochastic analysis, e.g. see papers \cite{Khas68, Ver, Kif}, books
   \cite{ FW, Krs} and references therein.    But in big
  majority of the corresponding works the fast motion of  angles $\f$
   is given not by an ODE, but by a stochastic differential equation  with
  $I$-depending coefficients and non-degenerate diffusion. Even if the fast motion  of $\f$
  was given by a non-random ODE,   then usually strong restrictions were imposed on  ergodic properties of the latter
  (e.g. see \cite[Section~II.3]{Skor}),  and system \eqref{Birk} fails to meet them.     Still, a  local version of system
   \eqref{2} was considered  in \cite{FrWe03}, where   for processes $(I,\f)(t)$
   whose $I$-components are defined in a  bounded domain the authors proved   convergence \eqref{conv1}
    till the time of exit of the vector $I(t)$  from the bounded domain above.
   However,  the techniques used in  \cite{FrWe03}  differ essentially from those in our proof.
  The exact statement of Theorem~\ref{t_main02}  and its proof, given in Sections~\ref{s_co_ac}-\ref{s_lp}, are crucially used in
  Sections~\ref{s_Birk_int}-\ref{s_9} of our work, which form its main second part (discussed below).
   The uniform in time convergence \eqref{conv3} seems to be a completely
  new result.
  \medskip

  Second part of the paper, made by Sections \ref{s_Birk_int}-\ref{s_9}, is dedicated to equation \eqref{3}, where we again assume
  that the frequency mapping $W$ is Anosov--nondegenerate and that equation \eqref{3} satisfies some other restrictions, similar
  to those, imposed on system \eqref{2}.
   It is convenient to rewrite the equation  in terms of the action-angle variables
$
(I,\f) \in \R_+^n \times \T^n,
$
where $I_k =\frac12\|\mathbf{v}_k\|^2$ and $\f_k = \arctan (v_k/ v_{-k})$. Then the corresponding  $I$-equation reads
\be \label{I_eq}
 dI_k =\eps {{P}}_k^I(I,\varphi)dt+\sqrt{\eps}B_{kj}^I(I,\varphi)d\beta_j(t), \quad k=1,\ldots,n,
 \ee
%where $P_k^I(v)=\mathbf{v}_k^tP_k(v)+\frac12\sum\limits_{j=1}^{n_1}$
where $P^I$ and $B^I$ are quadratic functions of $v, {\mathbf{P}}_k(v)$ and $B_{kj}$, while the $\f$-equation is
\be\label{f_eq}
 d \f_k(t) = W_k(I) dt + [ \dots],
\ee
where $ [ \dots]$ stands for the   drift and dispersion terms of order $\eps$ and $\sqrt\eps$, respectively, which are singular
when $I_k=0$. See equations \eqref{I-equation1} and \eqref{phi-equation1}.
 We wish to apply the results from the first part   to this  system of $(I, \f)$-equations. This is not at all
 straightforward since the coefficients of the  $(I,\f)$-system have singularities at the locus
 $$
 %{\aleph}
\aleph = \{v\in\mathbb R^{2n}: \bv_j=0 \text{ for some } j \},
 $$
 and the dispersion in  \eqref{I_eq} vanishes on $\aleph$.   This difficulty is overcome with the help of
 two additional groups of
 results. Firstly, in the rather technical Lemma~\ref{l_key_ham} we establishe that trajectories  of eq.~\eqref{3} stay
  in the vicinity of  locus   ${\aleph}$ only a short time, uniformly in $\eps$.
  This result allows us to prove in Theorem~\ref{th_ham_main}
  that if we write solutions $I^\eps $ of \eqref{I_eq}  as
  $
  I^\eps(\eps^{-1}\tau)$, $0\le \tau\le T,
  $
  then their laws are pre-compact in the space of  measures on
   $C(0,T; \R^n)$  and
   every limiting point of this family is a weak solution $I(\tau)$ of the corresponding averaged $I$-equation (see
  equation \eqref{averaged-1}). Unfortunately the latter is an equation without uniqueness of a solution (see
  Remark~\ref{remark-degenerate}), so we cannot conclude that the laws $\cD I^\eps(\eps^{-1}\tau)$ converge
  to a limit  when
  $\eps\to0$. To resolve this problem we use the second additional result. Namely, in Section~\ref{construct-effective}
  we return from the $(I,\f)$-system \eqref{I_eq}--\eqref{f_eq} back to the original equation \eqref{3}, remove from them
  the fast terms $W_k \bv_k^\perp dt$ and re-write the obtained equation in the slow time $\tau$ as
  $$
  dv (\tau) = P(v) d\tau + B(v) d\beta(\tau),
  $$
  where $P=(\bP_1, \dots, \bP_n)$, $B=(B_{kj})$ and $\beta= (\bbeta_1, \dots, \bbeta_{n_1})$. Then we formally average this
  equation with respect to the natural action of the torus $\T^n$ on $\R^{2n}$ (see \eqref{Phi_th}), using the rules of
  stochastic averaging. Thus we get a  stochastic equation
  \be\label{ef_eq}
  dv(\tau) = \lan P\ran(v) d\tau + \llan B\rran (v) d\beta(\tau),
  \ee
  which we call  {\it  the effective equation}. The coefficients of \eqref{ef_eq} are locally Lipschitz, so its solution, if exists, is
  unique. The key step of our analysis of the averaging in equation \eqref{3} is Theorem~\ref{lifting-thm},
  stating that the weak solution $I(\tau)$ of the averaged $I$-equation,  which has been
  obtained as a limit in law of some
  sequence $I^{\eps_j}$ of solutions for \eqref{I_eq}, written in the slow time $\tau$, can be lifted to a weak solution
  $v(\tau)$ of \eqref{ef_eq} (i.e. $I(\tau) = I(v(\tau))$). Since a solution $v(\tau)$  of \eqref{ef_eq}
  is unique, then as a consequence we get
  in Theorem~\ref{averaged-unique} that the actions $I^\eps$ of solutions for equation \eqref{3},
   written in  slow time $\tau$, converge in law   to a limit:
  \be\label{convI}
  \cD I^\eps(\eps^{-1}\tau) \strela \cD I(\tau) \quad\text{for}\;\; 0\le\tau\le T \quad\text{as}\;\; \eps\to0,
  \ee
  where $I(\tau)$ is a weak solution of the averaged $I$-equation. The latter
    can be lifted to a unique weak solution $v(\tau)$
  of \eqref{ef_eq}.

  Similarly to the uniform convergence \eqref{conv3},  we show in Theorem \ref{thm-uniform} that if some few moments
  of solutions for equation \eqref{3} are bounded uniformly in time and in $\eps$, and if the effective equation \eqref{ef_eq}
  is mixing, then  convergence \eqref{convI} is uniform in $\tau\ge0$.

  Proposition~\ref{sufficient1} provides an easy
  sufficient condition which implies that all the assumptions, required for the validity of results in
   Sections~\ref{s_Birk_int}-\ref{s_9} are met. In Section~\ref{ss_10.2} we discuss applications of our results
   to damped/driven Hamiltonian systems.

 When we apply Theorem \ref{averaged-unique} (which ensures convergence \eqref{convI}) to stochastic perturbations of
 ``real" integrable Hamiltonian equations,  we arrive at a difficulty due to the fact that often
 integrable Hamiltonian equations which appear in mechanics and physics
  can be put to the Birkhoff normal form \eqref{Birk} (with $W= \nabla h$) not on the
 whole $\R^{2n}$, but only locally.\footnote{But see Example \ref{ex_Fourier} in Section \ref{s_Birk_int} for a class of equations \eqref{3} in $\R^{2n}$  which appear in the
 non-equilibrium statistical physics. }
  Vey's theorem (see \cite{Eliasson},
   \cite[Section~2.3]{HGK}) and references there) provides an
 instrumental  sufficient condition which allows to write an integrable
  Hamiltonian equation in the form \eqref{Birk} in a neighbourhood of some point. Without lost of generality we assume that this point is the origin,
  and assume that the equation is written in the from \eqref{Birk} in the closure $\bar\fB_R$ of  domain
  $
  \fB_R = \{v: \| I(v)\| <R\},
  $
  for some $R>0$.  Accordingly, in Section~\ref{ss_local_version} we discuss a local problem of examining  equation
 \eqref{Birk} in  $\bar\fB_R$.   There   for any  initial data $v_0 \in \fB_R$ we consider the exit time
   $\theta^\eps_R$ of a corresponding trajectory $v^\eps(t)$     of \eqref{3}
    from  $\fB_R$. We show that, firstly, $\theta^\eps_R$ is a random
   variable of order $\eps^{-1}$ and,  secondly, that for the trajectory $v^\eps$, stopped at
   $t=\theta^\eps_R$,  a variation of Theorems~\ref{averaged-unique} applies and implies that the action-vector
   $I(v^\eps(t))$,  written in the slow time $\tau=\eps t$,
   for $\tau \le \eps \theta^\eps_R$ converges in distribution to a
   solution $I(\tau)$ of the averaged $I$-equation, stopped when $\| I \| =  R$.  Similar to  Theorem~\ref{averaged-unique},
   this solution $I(\tau)$ may be lifted to a solution $v(\tau)$
   of the effective equation \eqref{ef_eq},  stopped at $\partial \fB_R$.
   \medskip

   Long time behaviour of deterministic perturbations of integrable systems \eqref{Birk} (in various parts of the phase-space $\R^{2n}$)
    is a classical problem of dynamical systems. If the initial data are allowed to be arbitrarily close to the locus $\aleph$, additional
    difficulty appears. In the case of Hamiltonian perturbations see \cite[Section~6.3.7]{AKN} for corresponding KAM theorems,
    and see \cite{Nieder} for a version of Nekhoroshev's  theorem. It seems that for non-Hamiltonian perturbations of \eqref{Birk}
    no convenient averaging theorem, valid up to $\aleph$, is known. Similarly it seems that no systematical study of the averaging for
    stochastic perturbations of system \eqref{Birk} was performed before our work.
     \medskip

\noindent {\it On the proofs.}
  Our presentation and proofs
  are based on the approach to stochastic averaging,
   originated in the celebrated paper \cite{Khas68} by R.~Khasminskii. In that we partially
   follow our previous works \cite{HGK22, HKM,  KP, K}, dedicated to the averaging   in  stochastic  perturbations of linear
    (the first two paper)     and non-linear (the last two)  PDEs.
           The  Khasminskii approach, as presented in our work, is a flexible form of  stochastic averaging, applicable to
 study stochastic perturbations of integrable equations, linear and nonlinear, in finite and infinite dimensions. In
 particular, to study stochastic perturbations of linear and integrable PDEs. See \cite{HKM, HGK22} and references
 in \cite{HGK22} for perturbations of linear PDEs, including ``linear analogies" for them of Theorems~\ref{averaged-unique}  and
 \ref{thm-uniform}.  Also see \cite{KP} for an analogy of Theorem~\ref{t_main02} for
 perturbations of the KdV equation by dissipation $\eps\Delta$ and a white noise of order $\sqrt\eps$, and see
 \cite{K} for an analogy for that equation of Theorems~\ref{lifting-thm} and \ref{averaged-unique}.\,\footnote{ Also
   see  \cite[Section 4.3]{HGK} for a discussion of the results in
   \cite{KP,K}. }

\bigskip

\noindent{\it Notation.} For a matrix $A$ we denote by $A^t$ its transposed and by $\R^n_+$ denote
 the set of vectors in $\R^n$ with non-negative
components.  For a Banach space $E$ and $R>0$, $B_R(E)$ stands for the open $R$-ball $\{ e\in E: | e|_E < R\}$, and $\bar B_R(E)$ --
 for its  closure $ \{| e|_E \le R\}$;
% $\mathcal{P}(E)$ stands for the space of  probability Borel measures on $E$,
  $C_b(E)$ stands for  the space of bounded continuous function on  $E$, and $C([0,T], E)$ -- for
   the space of continuous curves
$[0,T] \to E$, given the sup-norm.
  By  $\mathcal{D}(\xi)$ we denote the  law of a random variable $\xi$,  by $\strela$ -- the weak convergence of measures, and by $\cP(M)$ -- the
  space of Borel probability  measures on a metric space $M$.  For a measurable mapping $F: M_1\to M_2$ and $\mu\in \cP(M_1)$
  we denote by $F\circ\mu\in \cP(M_2)$ the image of $\mu$ under $F$; i.e. $F\circ\mu(Q) = \mu(F^{-1}(Q))$.

  For $v=(v_1,\dots,v_n)\in\mathbb{R}^n$  we set $\|v\|^2=\sum_{k=1}^n|v_k|^2$ and $|v|=\sum_{k=1}^n|v_k|$.
  %We equip the space $\mathbb{R}^n\times \mathbb{T}^k\ni(I,\varphi)$ with the norm $|(I,\varphi)|=|I|+|\varphi|$.
   We denote $M_{n\times k}$ the space formed by $n\times k$ real matrices with the Hilbert-Schmidt norm $\|\cdot\|_{HS}$,  i.e. the square root of the sum of the squares of all its  elements.
   If     $L= \mathbb{R}^d\times \mathbb{T}^n$, $n\ge0$ (we set  $\T^0:=\{0\}$) and $m\ge0$,
       then $\text{Lip}_m(L, E)$ is the collection of  locally Lipschitz    maps $F:L \to E$ such
 that for any $R>0$ we have
\be\label{nota}
 (1+ |R|)^{-m} \!\! \sup_{\xi\in \BR(\R^d)\times \T^n  }|F(\xi)|_E   =: \cC^m(F) <\infty.
\ee
For a  set $Q$  we denote
 by $\mathbf{1}_Q$  its  indicator  function,  by $Q^c$ -- its complement, and by $\mathcal{L}(Q)$ --  its Lebesgue measure if $Q\subset\mathbb{R}^n$. For a function $f$, depending on angles $\varphi\in\T^n$ (and maybe on some other variables)
 we denote
 \be\label{average}
 \lan f \ran = (2\pi)^{-n}\! \int_{\T^n} f\,d\varphi.
 \ee
Finally, for real numbers $a$ and $b$, $a\vee b$ and $a\wedge b$ indicate their  maximum and minimum.

\section{Problem setup}
\label{s_ps}
Let $d_1\in \Nn$,
 $(\Omega, \mathcal{F}, \mathbf{P})$ be a  probability space,  $\beta(t)$, $t\geq 0$, be  a standard $d_1$-dimensional
Brownian motion  defined on it, and $\{\mathcal{F}_t\}$ be the natural filtration, generated by the process $\beta(s)$, $0\leqslant s\leqslant t$.

We start to examine a diffusive perturbation \eqref{2} of integrable systems \eqref{int_syst}.
There $\eps\in (0,1]$ is a small  parameter,   $P^I$ is a $d$-dimensional vector function, $\theta$ and $P^\varphi$ are
 $n$-dimensional  functions, while  $\Psi^I(\cdot)$ and $\Psi^\varphi(\cdot)$ are  $d\times d_1$\,- and $n\times d_1$-matrix functions.
Our first goal is to study system \eqref{2} for $0\le t\lesssim \eps^{-1}$.
After passing to the slow time $\tau=\eps t$, the  system  takes the form
\begin{equation}\label{ori_scaled}
\begin{split}
d I^\eps&=P^I(I^\eps,\varphi^\eps)d\tau+ \Psi^I(I^\eps,\varphi^\eps)d\beta(\tau),\\
d\varphi^\eps&=\big[\tfrac 1\eps\theta(I^\eps)+ P^\varphi(I^\eps,\varphi^\eps)\big]d\tau+\Psi^\varphi(I^\eps,\varphi^\eps)d\beta(\tau),
\end{split}
\end{equation}
where $0\le \tau \le T$ for some fixed $T>0$. It  is equipped with an initial condition
\begin{equation}\label{ini_cond_det}
I^\eps(0)=I_0,\quad \varphi^\eps(0)=\varphi_0.
\end{equation}
Here
$(I_0,\varphi_0)\in \mathbb R^d\times\mathbb T^n$ is either deterministic, or is a r.v.,  independent of the process $\beta$.
We will mostly dwell on the first case since
 the second  can be directly generalized from the first one (see Remark \ref{remark-main-thm} and Amplification \ref{amp_ran_ini}).
Our goal is to examine  the limiting behaviour of the distribution of  a solution
 $(I^\eps(\cdot), \varphi^\eps(\cdot))$ as $\eps\to0$. In particular, to
show that in the limit the law of $I^\eps(\cdot)$ is a weak solution of a certain  averaged equation,
 independent of     $\varphi$ and  $\eps$.

In what follows  we always assume that the following conditions are fulfilled for  system \eqref{ori_scaled}--\eqref{ini_cond_det}:

\begin{assumption}\label{assumption-1}
(1) The Lebesgue measure of $I\in \mathbb R^d$ for which $\theta(I)$ is rationally dependent equals  zero, that is $\mathcal{L}\big(\bigcup_{k\in\mathbb{Z}^n \setminus \{0\} }\{I\in\mathbb{R}^d: k\cdot\theta(I)=0\}\big)=0$.

%[$\boldsymbol{{\mathcal{C}}2}$.]
(2) The  matrix  $a(I,\varphi)=\big(a_{ij}(I,\varphi)\big):=\Psi^I(I,\varphi)(\Psi^I)^t(I,\varphi)$ satisfies the uniform ellipticity condition, that  is,
there exists $\lambda>0$ such that
$$
\lambda |\xi|^2\leq a(I,\varphi)\xi\cdot\xi\leq \lambda^{-1}|\xi|^2
$$
for all $\xi\in \mathbb R^d$ and all $(I,\varphi)\in \mathbb R^d\times\mathbb T^n$.

(3) There exists $q>0$ such that $\theta\in\text{Lip}_{q}(\mathbb{R}^d,\mathbb{R}^n)$, $P^m\in \text{Lip}_q(\mathbb{R}^d\times\mathbb{T}^n,\mathbb{R}^d)$ and $\Psi^m\in\text{Lip}_q(\mathbb{R}^d\times\mathbb{T}^n,M_{d_m\times d_1})$, $m=I,\varphi$ and $d_I=d, d_\f=n$
(see \eqref{nota}).

(4) There exists $T>0$ such that  for every $(I_0,\varphi_0)\in\mathbb{R}^d\times\mathbb{T}^n$,  system
 \eqref{ori_scaled}--\eqref{ini_cond_det}  has a unique strong solution $\big(I^{\eps}(\tau),\varphi^{\eps}(\tau)\big):=\big(I^{\eps},\varphi^{\eps}\big)(\tau; I_0,\varphi_0)$, $\tau\in [0,T]$, equal $(I_0,\varphi_0)$ at $\tau=0$. Moreover, there exists $q_0>(q\vee2)$ such that
\begin{equation}
\label{apriori-1}
\mathbf{E}\sup_{\tau\in[0, T]}| I^\eps(\tau)|^{2q_0}\leq C(|I_0|,T)\qquad \forall \eps\in(0,1],
\end{equation}
where $C(\cdot)$ is a non-negative continuous function on $\mathbb{R}_+^2$ non-decreasing in both arguments.
\end{assumption}

Throughout the text except the last Section \ref{s_9},
the time $T>0$ is fixed and the dependence on it usually is not indicated. The process $(I^{\eps}(\tau),\varphi^{\eps}(\tau))$, $\tau\in[0,T]$ is always understood as  a unique strong solution of the system \eqref{ori_scaled}--\eqref{ini_cond_det}.
\smallskip

\begin{remark}\label{r_Anosov}
Item (1) in Assumption~\ref{assumption-1} is called {\it the Anosov condition}\,\footnote{``The set of slow variables $I$ for which
the motion of fast variable $\f$ is not ergodic has zero measure," see \cite{LM}, p.~12, assumption iii).}
and is rather mild. E.g. it clearly holds if  the mapping $\theta$ is $C^1$-smooth, $n\ge d$ and the set of $I$'s for which
rank$\;\partial_I\theta(I)<d$ has zero measure. But it also may hold for systems \eqref{int_syst} with $n<d$ (even for systems with $n=1$).
In particular, it holds for any $n,d$ if the mapping $\theta$ is analytic and satisfied {\it R\"ussmann's condition}: there exists $N\in\Nn$
such that for every $I\in \R^d$ the vectors
$$
\frac{\p^{|q|} \theta(I)}{\p_1^{q_1}\dots \p_d^{q_d}}, \quad q \in \Z_+^d, \quad |q| \le N,
$$
jointly span $\R^n$. Indeed, if this condition holds, then for any non-zero vector $s\in \R^n$  the function $\theta(I) \cdot s$ does not
vanish identically (see \cite[Section~6.3.2]{AKN}, item $7^\circ$). Then by analyticity the set in (1) has zero measure.
\end{remark}

Since item (4) of Assumption~\ref{assumption-1} is formulated not in terms of  coefficients of equations \eqref{ori_scaled}, we provide here a sufficient
condition that ensures its validity.
\begin{proposition}\label{l_suff_assum_i4}
  Let all  coefficients of equation  \eqref{ori_scaled} be globally Lipschitz continuous. Then for any $T>0$ and
  $(I_0,\varphi_0)\in \mathbb R^{d}\times \mathbb T^n$ problem \eqref{ori_scaled}--\eqref{ini_cond_det} has a unique solution, and inequality \eqref{apriori-1} holds for  every $q_0\in\mathbb{N}$.
\end{proposition}

\begin{proof}
Under the above assumptions, for any $(I_0,\varphi_0)$ system  \eqref{ori_scaled}-\eqref{ini_cond_det} has a
unique solution.\footnote{{ To prove this we regard \eqref{ori_scaled} as an equation on $\R^{d}\times\R^{n}$ with periodic
in $\f$ coefficients, evoke the usual theorem on stochastic equations in $\R^{d+n}$ with Lipschitz coefficients
(see for instance \cite[Theorem~2.9]{brownianbook}) to get a solution for this equation, and next  apply the projection
 $ \R^{d}\times\R^{n} \to \R^{d}\times\T^{n}$ to obtain a solution in $ \R^{d}\times\T^{n}$.}}
 For  $q_0\in\mathbb{N}$ we need to show that
estimate  \eqref{apriori-1} holds with a constant $C(q_0,I_0,T)$  that does not depend on $\eps$.  Let us
 fix an arbitrary $R>\|I_0\|^{2q_0}$  and introduce the stopping time $\tau^\eps_R=\inf\{\tau>0\,:\,\|I^\eps(\tau)\|^{2q_0}>R\}$.
%=\inf\{\tau>0\,:\,I(v^\eps(\tau))>\frac12R\}$.
By It\^o's formula,
 the process $\|I^\varepsilon(\tau\wedge\tau_R^\varepsilon)\|^{2q_0}$ satisfies the equation
\begin{equation}\label{q0-I-norm}
\begin{split}
&\|I^\eps(\tau\wedge\tau_R^\eps)\|^{2q_0}\\
=&\|I_0\|^{2q_0}+2q_0\int_0^{\tau\wedge\tau_R^\eps}\|I^\eps(s)\|^{2(q_0-1)}\Big(I^\eps(s), P^I(I^\eps(s),\varphi^\eps(s))\Big)ds\\
&+\frac{1}{2}\int_0^{\tau\wedge \tau_R^\eps}\text{Trace}\big(\Psi^I(I^\eps,\varphi^\eps)\big)^t\nabla_I^2\|I^\eps\|^{2q_0}\Psi^I(I^\eps,\varphi^\eps)ds\\
&+\int_0^{\tau\wedge\tau_R^\eps}2q_0\|I^\eps(s)\|^{2(q_0-1)}\Big(I^\eps(s), \Psi^I(I^\eps(s),\varphi^\eps(s))d\beta(s)\Big).
\end{split}
\end{equation}
Taking the expectation and using the  global Lipschitz continuity of coefficients we conclude that
$$
\mathbf{E}\|I^\eps(\tau\wedge\tau_R^\eps)\|^{2q_0}\leqslant C_1+C_2\int_0^{\tau}\mathbf{E}\|I^\eps(s\wedge\tau_R^\eps)\|^{2q_0}ds
$$
with constants $C_1$ and $C_2$ that depend continuously  on $q_0$, $\|I_0\|$ and  the global Lipschitz constants of the coefficients and do not depend on $\eps$ and on $R$.
%Summing up these inequalities over $k$ we obtain
%$$
%\mathbf{E}\|I^{\eps,R}(\tau)\|^2 \leqslant C'_1+C'_2\int_{0}^{\tau} \mathbf{E}\|I^{\eps,R}(s)\|^2ds.
%$$
By  Gronwall's  lemma,
$$\mathbf{E}\|I^\eps(\tau\wedge\tau_R^\eps)\|^{2q_0}\leqslant C_1\exp(C_2\tau).$$
 Since for any fixed $\tau$
the sequence ${\tau\wedge\tau^{\eps}_R}$ a.s. converges to $\tau$ as $R\to\infty$,  then by the Fatou lemma,
\begin{equation}\label{q0-I-norm2}\mathbf{E}\|I^\eps(\tau)\|^{2q_0}  \leqslant C_1\exp(C_2\tau),\;\forall \tau\geqslant0.\end{equation}

Now consider equation \eqref{q0-I-norm} without the stopping time $\tau_R^\eps$. For any $T>0$, we have
\[\mathbf{E}\sup_{0\leqslant\tau\leqslant T}\|I^\eps(\tau)\|^{2q_0}\leqslant \|I_0\|^{2q_0}+C\int_0^T\mathbf{E}\|I^\eps(\tau)\|^{2q_0}d\tau+C\mathbf{E}\sup_{\tau\in[0,T]}M(\tau),\]
where $M(\tau)$ is the martingale term in \eqref{q0-I-norm}, and the constant $C$ does not depend on $\eps$.  Due to \eqref{q0-I-norm2}, the estimate \eqref{apriori-1} follows by applying  Doob's inequality to $M(\tau)$.
\end{proof}

Assumption~\ref{assumption-1}.(4) also holds for systems \eqref{2} where the coefficients are not globally Lipschitz, but items (2) and (3)
of the assumption are valid
and the vector field $P$ is coercive. Cf. below Proposition~\ref{sufficient1}, where this is discussed for system
\eqref{3}.

\section{Tightness and  averaged equation}\label{s_co_ac}

In this section we first show the collection of the laws of  the  processes
$I^{\eps}(\tau), \tau\in[0,T]$ with  ${\eps}\in(0,1]$ is tight. Then we introduce the averaged equation
 and finally prove  key technical statements.

Since $T>0$ is fixed, then dependence of constants on it usually is not indicated.

\subsection{Tightness}\label{ss_3.1}
{
Notice that under condition (3) of Assumption \ref{assumption-1}  there exists
an increasing function $\nu(M)\,:\:\R_+\to\R_+$ such that  $\nu(M)\to\infty$ as $M\to\infty$,
and
\be\label{nu_Lip}
\begin{split}
\text{
 in the set $\{(I,\varphi)\,:\,|I|\leq \nu(M)\}$ norms and  Lipschitz constants}\\
\text{  of all coefficients in equations \eqref{ori_scaled}
do not exceed $M$}
 \end{split}
\ee
(the coefficient $\frac1{\eps} \theta$ should be  taken without the  factor $\eps^{-1}$).
This a bit unusual way to write the locally Lipschitz property is convenient for the calculation below.
}

\begin{lemma}\label{l_compa}
For 	any fixed $(I_0, \f_0)$ the family of laws of processes $I^\eps(\tau)$ $(\tau\in[0,T], {\eps}\in(0,1])$ which are the $I$-components of
solutions $(I^\eps, \f^\eps)(\tau; I_0, \f_0)$, is tight
%$$
%I^\eps(\tau)=I_0+\int\limits_0^\tau P^I(I^\eps(s),\varphi^\eps(s))ds+\int\limits_0^\tau \Psi^I(I^\eps(s),\varphi^\eps(s))d\beta(s),$$
%is tight
in the  Banach space $C([0,T];\mathbb R^d)$.
\end{lemma}
\begin{proof}
  According to (4) of Assumption \ref{assumption-1}, for any $\delta>0$ there exists $R=R(\delta)>0$ such that
  $\mathbf{P}\{\sup\limits_{\tau\in[0, T]}|I^\eps(\tau)|>R\}<\delta$.  Denoting by $\tau^\eps_R$ the exit time
  $\tau^\eps_R=\inf\{\tau>0\,:\, |I^\eps(\tau)|>R\}$ and denoting   by $(I^\eps_R,\varphi^\eps_R)$ the stopped  process
  $(I^\eps(\tau\wedge\tau^\eps_R),\varphi(\tau\wedge\tau^\eps_R))$, we have
   $\mathbf{P}\{(I^\eps(\cdot),\varphi^\eps(\cdot))\not=(I^\eps_R(\cdot),\varphi^\eps_R(\cdot)) \ \hbox{on }[0,T]\}<\delta$.

  By (3) of Assumption \ref{assumption-1}  the functions  $P^I$ and $\Psi^I$ are bounded in  the ball $\{|I|\leq R(\delta)\}$. Therefore
for any moments $0\leq \tau_1\leq \tau_2\leq T$ we have
\[\begin{split}
\mathbf{E}\big\{|I_R^\eps(\tau_2)-I_R^\eps(\tau_1)|^4\big\}&=
\mathbf{E}\bigg\{\Big|\int\limits_{\tau_1\wedge\tau^\eps_R}^{\tau_2\wedge\tau^\eps_R}P^I(I^\eps,\varphi^\eps)ds
+\int\limits_{\tau_1\wedge\tau^\eps_R}^{\tau_2\wedge\tau^\eps_R}\Psi^I(I^\eps,\varphi^\eps)d\beta(s)\Big|^4\bigg\}\\
&\leq C_R(|\tau_2-\tau_1|^4+|\tau_2-\tau_1|^2),
\end{split}\]
and the assertion follows by a direct application of the Prokhorov theorem, cf. \cite[Lemma~2.2]{GHKu}.
\end{proof}

\subsection{The averaged equation }\label{ss_3.2}
The vector field  $P^I(I,\f)$  in the $I$-equation of \eqref{ori_scaled} depends explicitly both on $I$-variables and $\varphi$-variables.
Since we are  interested in the evolution of the $I$-component of the system as ${\eps}\to0$, then we introduce
in consideration the vector field $\langle P^I\rangle$, obtained by  the averaging of $ P^I$  in angles  $\f$,
%\begin{equation}\label{eff_driftt}
$$
\langle P^I\rangle(I)=\int_{\mathbb T^n}P^I(I,\varphi)\,d\varphi.
$$
If an $I$ is such that the vector $\theta(I)$ is rationally independent, then
\be\label{another_def}
 \lim_{N\to\infty} \frac1{N} \int_0^N P^I\big(I, \f + t\theta(I) \big)\, dt
=\langle P^I\rangle(I),
\ee
and the convergence is uniform in $\f$.
Moreover, $\langle P^I\rangle \in \text{Lip}_q(\mathbb{R}^d,\mathbb{R}^d)$ with the same $q$ as  for $P^I$,
 $\mathcal{C}^q(\langle P^I\rangle)=\mathcal{C}^q(P^I)$,  and for $\langle P^I\rangle $ the function $\nu(N)$ (as in \eqref{nu_Lip}
 is the same as for $P^I$;  see  \cite[Section~3.2]{GHKu}.
The convergence to the limit on the l.h.s. of \eqref{another_def}  is the faster, the more diophantine is the vector $\theta(I)$ (see \cite[Section~6.1.5]{AKN} for a related
discussion).

Similarly we set
%\begin{equation}\label{eff_diffu}
$$
\langle a^I\rangle(I)=\int_{\mathbb T^n}
\Psi^I(I,\varphi)\big(\Psi^I(I,\varphi)\big)^t\,d\varphi.
$$
By (2) of Assumption \ref{assumption-1},
$$\lambda|\xi|^2 \leqslant \langle a^I\rangle(I)\xi\cdot\xi\leqslant \lambda^{-1}|\xi|^2, \quad\forall \xi, I\in\mathbb{R}^d.$$
 Let $\llangle \Psi^I \rrangle(I)$ be the principal square
 root of $\langle a^I\rangle(I)$.\footnote{I.e. $\llangle \Psi^I \rrangle(I)$ is a non-negative self-adjoint matrix such that
 $\big(\llangle \Psi^I \rrangle\big)^2 =\langle a^I\rangle$. }
 By  the estimates above and
  \cite[Theorem~5.2.1]{SV}, $\llangle \Psi^I \rrangle (\cdot)$ belongs to
$ \text{Lip}_{q}(\mathbb{R}^d,M_{d\times d})$ and satisfies the   uniform ellipticity condition as in (2) of Assumption \ref{assumption-1} (with the same $\lambda>0$).

We introduce the averaged equation for the limiting as ${\eps}\to0$ evolution of $I$-variables,
\begin{equation}\label{averaged-eq-1}
dI(\tau)=\langle P^I\rangle(I(\tau))d\tau+\llangle \Psi^I \rrangle(I(\tau))dW(\tau), \quad I(0)=I_0,
\end{equation}
where $W(\tau)$ is a standard $d$-dimensional Wiener process. As we  discussed above,  coefficients
 of \eqref{averaged-eq-1} are locally  Lipschitz, therefore for each $I_0\in\mathbb{R}^d$ the
 equation  has at most one  solution $I(\tau;I_0)$. We will show in the next section that as $\eps\to0$, the
 $I$-component $I^\eps(\tau)$ of a solution for \eqref{ori_scaled}  with $(I(0), \f(0))=(I_0, \f_0)$, for any $\f_0\in \T^n$
 converges in law to  solution  $I(\tau;I_0)$ (in particular, the latter exists).
The key technical results, needed  to establish this  convergence, are proved in the next subsection.

\subsection{Main lemmas}\label{ss_3.3}
\begin{lemma}\label{lemma-key-1} For the process $(I^{\eps}(\tau),\varphi^{\eps}(\tau))$, $\tau\in [0,T]$, we have
\begin{equation}\label{lim_rel_first}
\Upsilon^\eps :=
 \EE\, \max\limits_{0\leq \tau\leq T}
  \Big|\int_0^\tau P^I(I^\eps(s),\varphi^\eps(s))\,ds-\int_0^\tau \langle P^I\rangle(I^\eps(s))\,ds\Big|\to 0 \quad \text{as}\;\; \eps\to0.
\end{equation}
\end{lemma}
 A proof of this relation relies on a number of auxiliary
  statements, given  below.

Let us first define a family of sets $ \mathcal{A}^\delta_{N,R}$, where $N,R>0$ and $0<\delta\le1$, made by vectors  $I$ such that
 $\theta(I)$ has poor diophantine properties, and so  the rate of convergence in \eqref{another_def} is slow:
%}
\begin{equation}\label{def_AN}
  \mathcal{A}^\delta_{N,R}=\Big\{I\in \mathbb R^d\,:\, |I|<R,\
  \max\limits_{\varphi\in\mathbb T^n}\Big|\frac1N\int_0^N \!\! P^I(I,\varphi+t \theta(I))\,dt
   -\langle P^I\rangle(I)  \Big|>\delta \Big\}.
\end{equation}
We set $\mathcal{A}^\delta_{N}:=  \mathcal{A}^\delta_{N,\infty} = \cup_{R>0}  \mathcal{A}^\delta_{N,R}$.
\begin{lemma}\label{l_small_meaA}
For any $\delta>0$ and  $R>0$ we have
$$
\lim\limits_{N\to\infty}\mathcal{L}(\mathcal{A}^\delta_{N,R})=0
$$
\end{lemma}

\begin{proof}
For $N>0$ denote
$
b_N(I) =  \max\limits_{\varphi\in\mathbb T^n}\Big|\frac1N\int_0^N
\!\! P^I(I,\varphi+t \theta(I))
dt  -   \langle P^I\rangle(I) \Big|.  %}
$
If the vector $\theta(I)$ is non-resonant, then by \eqref{another_def}  $b_N(I) \to0$ as $N\to\infty$,
So, by (1) of Assumption~\ref{assumption-1},
$b_N(I) \to  0$ a.s., and the assertion follows since the a.s. convergence implies the convergence in measure.
\end{proof}

\begin{lemma}\label{prop_nondege}
For any $R, N,\delta>0$ the  probability $\mathbf{P}\{I^\eps(\tau)\in \mathcal{A}_{N,R}^\delta\}$ admits the upper bound
  \begin{equation}\label{int_est_probb}
 \mathbf{E}
  \int_0^{T\wedge\tau^\eps_R} \mathbf{1}_{ \mathcal{A}_{N,R}^\delta}(I^\eps(\tau))\,d\tau
    \leq C(R)\big(\mathcal{L}( \mathcal{A}_{N,R}^\delta)\big)^{\frac1d},
  \end{equation}
  where $\tau^\eps_R$ is the exit time of $I^\eps$ from the ball $B_R=\{I\in\mathbb R^d\,:\,|I|<R\}$.
  The constant $C(R)$ does not depend on $\eps$,  $N$ and $\delta$.
\end{lemma}
\begin{proof}
Since the process $I^\eps(\tau)$ satisfies  It\^o's equation
$
dI =P^I d\tau + \Psi^I d\beta(\tau),
$
where the diffusion matrix $ \Psi^I(\Psi^I)^t$ is uniformly non-degenerate by (2) of Assumption~\ref{assumption-1}, then the assertion follows from  \cite[Theorem~2.2.4]{Kry77}.
Indeed, by the latter result, choosing there $f(\tau, \cdot)$  to be  the characteristic function
  of the set $\mathcal{A}^\delta_{N,R}$, we obtain
  $$
% \int_0^{T\wedge\tau^\eps_G}\mathbf{P}\{z^\eps(t)\in \mathcal{A}_{N,G}^\delta\}\,dt =
  \mathbf{E}
  \int_0^{T\wedge\tau^\eps_R}\mathbf{1}_{ \mathcal{A}_{N,R}^\delta}(I^\eps(\tau))\,d\tau
  \leq C(R) \|\mathbf{1}_{ \mathcal{A}_{N,R}^\delta}\|\big._{L^d(B_R)}=
  C(R)\big(\mathcal{L}( \mathcal{A}_{N,R}^\delta)\big)^{\frac1d},
  $$
  and the assertion of the lemma is proved.
\end{proof}

Due to \eqref{int_est_probb},
  \begin{equation}\label{int_unbou_est_probb}
  \int_0^{T} \mathbf{P}\{I^\eps(\tau)\in \mathcal{A}_{N}^\delta\}\,d\tau
    \leq C(R)\big(\mathcal{L}( \mathcal{A}_{N,R}^\delta)\big)^{\frac1d}+T\mathbf{P}\{\tau_R^\eps<T\}.
  \end{equation}
Since $\{\omega\,:\,\tau_R^\eps<T\} {\, \subset \,}
 \{\omega\,:\, \sup_{0\le \tau \le T} | I^\eps (\tau) | \ge R\}$, then by (4) of Assumption \ref{assumption-1}
  and Chebyshev's inequality we get that  probability $\mathbf{P}\{\tau_R^\eps<T\}$ tends to zero as $R\to\infty$, uniformly in $\eps\in(0,1]$.
Combining  this fact with  Lemma~\ref{l_small_meaA} to estimate the r.h.s. of \eqref{int_unbou_est_probb}
 we conclude that for any  $\delta>0$  there exists a positive function
 $
 N\mapsto  {\alpha_N^\delta},
 $
 $N>0$, converging to zero as $N\to\infty$, such that
 \begin{equation}\label{est_small_gen}
   \int_0^{T} \mathbf{P}\{I^\eps(\tau)\in \mathcal{A}_{N}^\delta\}\,d\tau\leq  {\alpha_N^\delta},\quad \forall\,\eps\in(0,1].
\end{equation}

We are now in  position to prove Lemma~\ref{lemma-key-1}

\begin{proof}[of Lemma~\ref{lemma-key-1}] {\bf Step 1:} Take  $N>0$ and $\tau_0  \in[0,\eps N)$ which will be specified later
and consider a partition of  $[\tau_0,T)$
to  intervals $[\tau_j, \tau_{j+1})=:\Delta_j$, $0\le j \le {j_N}$, where
$ \tau_j=\tau_0+j \eps  N$, \,${j_N}$ is  the biggest $j$ such that $\tau_j <T$, and $\tau_{j_{N}+1} :=T$.
We assume that $\eps N\leq \frac13 (1\wedge T)$ -- below
we deal with $N$'s such that  $\eps N\ll 1$, so  this assumption is not a restriction.  Then
$$
j_N \le {2}/{(\eps N)}.
$$

 Next let us  introduce the random variable
$$
L^{\eps,\delta}_{N,\tau_0} (\omega)=
\#\big\{j\in [0, {j_N}]
\,:\, I^\eps(\tau_j)\in \mathcal{A}_N^\delta\big\},
$$
 which counts the moments $\tau_j$'s for which  the frequency vector $\theta(I^\eps(\tau_j))$ has poor
diophantine properties.  Since
$$
L^{\eps,\delta}_{N,\tau_0}=\sum_{j=0}^{{j_N}}\mathbf{1}_{\mathcal{A}_N^\delta}(I^\eps(\tau_j)),
$$
we obtain
$$
\mathbf{E} L^{\eps,\delta}_{N,\tau_0}=\sum_{j=0}^{{j_N}}\mathbf{P}\{I^\eps(\tau_j)\in \mathcal{A}_N^\delta\}.
$$
Integration of this equality in the variable $\tau_0$ over  interval $[0, \eps N)$ yields
\[\begin{split}
&\int_0^{\eps N}\mathbf{E}L^{\eps,\delta}_{N,\tau_0}\,d\tau_0=\sum\limits_{j=0}^{{j_N}}\int_0^{\eps N}
\mathbf{P}\{I^\eps(\tau_j(\tau_0))\in \mathcal{A}_N^\delta\}\,d\tau_0   \\
& =
\int_0^{T}
\mathbf{P}\{I^\eps(\tau)\in \mathcal{A}_N^\delta\}\,d\tau\leq \alpha_N^\delta,\end{split}
\]
where
the last estimate  follows from \eqref{est_small_gen}.  Therefore there exists a {\it non-random}  number
 $\tau_0^*\in[0,\eps N)$ such that
\begin{equation}\label{choice_point}
  \mathbf{E}L^{\eps,\delta}_{N,\tau_0^*}\leq  \alpha_N^\delta(\eps N)^{-1}.
\end{equation}
From now on we fix this $\tau_0^*$ for  the choice of $\tau_0$ in the definition of  the partition $\{\Delta_j\}$
 of interval $[\tau_0,T)$. So from now on
$
\tau_j := \tau_0^* + j\eps N
$
for all $j$, and   below we write $L^{\eps,\delta}_{N,\tau_0}$ simply  as $L^{\eps,\delta}_{N}$.

{\bf Step 2:}
We define the first  good event for our argument (there will be three of them) as a collection $\mathcal{E}_1$
 of all $\omega$ such that $L^{\eps,\delta}_{N} (\omega)$ is relatively small:
$$
 \mathcal{E}_1 =
  \big\{L^{\eps,\delta}_{N}\leq \ (\alpha_N^\delta)^\frac12 (\eps N)^{-1}\big\}.
$$
In view of \eqref{choice_point} and %}
 Chebyshev's  inequality, for the complement  $\mathcal{E}_1^c=\Omega\setminus \mathcal{E}_1$ we have
\begin{equation}\label{after_cheb}
\mathbf{P}( \mathcal{E}^c_1)
%\mathbf{P}\big\{L^{\eps,\delta}_{N,t_0}\geq \ (\alpha_N^\delta(T+1))^\frac12 (\eps N)^{-1}\big\}
 \leq  (\alpha_N^\delta)^\frac12.
\end{equation}

Due  to (3) and (4) of 	Assumption \ref{assumption-1}, for any $j$
  the stochastic %and the absolutely continuous
terms on the right-hand sides of \eqref{ori_scaled} admit the following upper bounds:
\[ \begin{split}
&\mathbf{E}\Big(\Big|\int_{\tau_j}^{\tau_{j+1}}\Psi^m(I^\eps(s),\varphi^\eps(s))\,d\beta(s)\Big|^2\Big)
\leq \mathbf{E}\int_{\tau_j}^{\tau_{j+1}}|\Psi^m(I^\eps(s),\varphi^\eps(s))|^2\,ds \\
&\leq C\int_{\tau_j}^{\tau_{j+1}}\mathbf{E}\big((|I^\eps|+1)^{2q}\big)d\tau\leq C_1\eps N,\qquad m=I,\,\varphi.
\end{split}
\]
Therefore, by  Doob's  inequality we  have
\begin{equation}\label{auxi_one}
\mathbf{E}\Big(\sup\limits_{\tau \in \Delta_j }\Big|\int_{\tau_j}^{\tau}\Psi^m(I^\eps(s),\varphi^\eps(s))\,d\beta(s)\Big|^2\Big)\leq C_2\eps N,\qquad m=I,\,\varphi.
\end{equation}
Similar,  by (3) and (4) of 	Assumption \ref{assumption-1} we have
\be\label{after_8_10}
\begin{split}
\mathbf{E}\Big(\sup\limits_{\tau\in\Delta_j}\Big|\int_{\tau_j}^{\tau}P^m(I^\eps(s),\varphi^\eps(s))\,ds\Big|^2\Big)
%\mathbf{E}\Big(\sup\limits_{\tau\in\Delta_j}\Big[(\tau-\tau_j)\int_{\tau_j}^{\tau}\big|P^m(I^\eps(s),\varphi^\eps(s))\big|^2\,ds\Big]\Big)\\
%\leq
%\eps N \mathbf{E}\Big(\sup\limits_{\tau\in\Delta_j}\int_{\tau_j}^{\tau}\big|P^m(I^\eps(s),\varphi^\eps(s))\big|^2\, \Big)
\leq
\eps N \mathbf{E}\Big(\int_{\tau_j}^{\tau_{j+1}}\big|P^m(I^\eps(s),\varphi^\eps(s))\big|^2\,ds  \Big)\\
\leq %\eps N\mathbf{E}\int_{\tau_j}^{\tau_{j+1}}(|I^\eps|+1)^{2q}d\tau
\eps NC\int_{\tau_j}^{\tau_{j+1}}\mathbf{E}\big((|I^\eps|+1)^{2q}\big)d\tau \leq
C_1(\eps N)^2,\qquad \; m=I,\,\varphi.
\end{split}
\ee

{\bf Step 3:}
For $j=0, \dots, {j_N}$ we measure the size of the perturbative part in eq.~\eqref{ori_scaled} on interval $\Delta_j$ by the random variable
$$
\zeta_j=\sum_{m=I,\varphi}\Big(  \ \sup\limits_{\tau\in\Delta_j}\Big|\int_{\tau_j}^\tau
\Psi^m(I^\eps(s),\varphi^\eps(s))\,d\beta(s)\Big|+
%\sum_{m=I,\varphi}\ \
 \sup\limits_{\tau\in\Delta_j}\Big|\int_{\tau_j}^\tau
P^m(I^\eps(s),\varphi^\eps(s))\,ds\Big| \Big)\,,
$$
and introduce another random variable which counts the number of  large $\zeta_j$'s:
$$
\widetilde{L}^{\eps,\delta}_{N}=\#\big\{j\,: \zeta_j  \geq \eps^\frac14 \big\}.
$$
%Then $\widetilde{L}^{\eps,\delta}_{N}=\#\{\tau_j\,:\, \zeta_j \geq \eps^\frac14\}$.
We have  $\mathbf{E}\zeta_j\leq C_3({\eps N})^\frac12$. Indeed,
 $ \mathbf{E}\zeta_j \leq  \big( \mathbf{E}\zeta^2_j  \big)^{\frac12}$ and due to \eqref{auxi_one} and \eqref{after_8_10}
\[
\begin{split}
 \mathbf{E}\zeta^2_j
 \leq
 C\big(\eps N+(\eps N)^2\big) \leq C \eps N
\end{split}
\]
(we recall that $\eps N \le 1/3$).
By  Chebyshev's  inequality
$\mathbf{P}\{\zeta_j\geq \eps^\frac14\}\leq C_3N^\frac12\eps^\frac14$, then
$$
\mathbf{E}\widetilde{L}^{\eps,\delta}_{N}=\mathbf{E}\Big(\sum_j \mathbf{1}_{\{\zeta_j\geq \eps^\frac14\}}\Big)
=\sum_j \mathbf{E}\big(\mathbf{1}_{\{\zeta_j\geq \eps^\frac14\}}\big)=\sum_j \mathbf{P}\{\zeta_j\geq \eps^\frac14\}\leq
\frac{T C_3 N^\frac12\eps^\frac14}{{\eps N}}=C_3TN^{-\frac12}\eps^{-\frac34}.
$$
Now we define the second good event for our argument as the set, where $ \widetilde{L}^{\eps,\delta}_{N}$ is not too big:
$$
\mathcal{E}_2=\big\{\omega\in\Omega\,:\, \widetilde{L}^{\eps,\delta}_{N}\leq \eps^{-\frac78}\big\}.
$$
Again due to  Chebyshev's  inequality,  the probability of  its complement satisfies
\begin{equation}\label{est_tildeL}
  \mathbf{P}(\mathcal{E}^c_2)\leq C_3TN^{-\frac12}\eps^\frac18.
\end{equation}
{
Denote
$$
\mathcal{M}_j:=\{\omega\in\Omega\,:\,\zeta_j\leq \eps^\frac14\quad \text{and} \quad
\sup\limits_{0\leqslant\tau\leqslant T} |I^\eps(\tau)|\leq \nu(M)\},
$$
where $\nu(M)$ is defined in \eqref{nu_Lip}.  In what follows we assume that $M\leqslant N$. Then for each $\omega \in \mathcal{M}_j$ on the interval $\Delta_j$
the curve $I(\tau)$ is close to $I(\tau_j)$:
}
\be \label{first_interv}
\begin{split}&\sup\limits_{\tau\in\Delta_j}\big|I^\eps(\tau)-I^\eps(\tau_j)\big|\\
&\leq \sup\limits_{\tau\in\Delta_j}\Big|\int_{\tau_j}^{\tau}P^I(I^\eps(s),\varphi^\eps(s))\,ds\Big|+
\sup\limits_{\tau\in\Delta_j}\Big|\int_{\tau_j}^{\tau}\Psi^I(I^\eps(s),\varphi^\eps(s))\,d\beta(s)\Big|
\leq
\eps^\frac14.
\end{split}\ee
Due to the definition of $\nu(M)$, for  $\omega\in \mathcal{M}_j$
$$
\Big|\int_{\tau_j}^{\tau}(\theta(I^\eps(s)-\theta(I^\eps(\tau_j))\,ds\Big|\leqslant N\int_{\tau_j}^{\tau}
\sup\limits_{\tau\in\Delta_j}\big|I^\eps(s)-I^\eps(\tau_j)\big|\,ds\leqslant \eps^\frac54 N^2.
$$
Therefore, for $\omega$ from the event $\mathcal{M}_j$,  on the interval $\Delta_j$
 the curve  $\f(\tau)$ is close to its linear approximation:
\be \begin{split} \label{second_interv}
&\sup\limits_{\tau\in\Delta_j}\big|\varphi^\eps(\tau)-\varphi^\eps(\tau_j)-\frac{\tau-\tau_j}\eps\theta(I^\eps(\tau_j))\big|\\
&\leq
\sup\limits_{\tau\in\Delta_j}\Big|\int_{\tau_j}^{\tau}P^\varphi(I^\eps(s),\varphi^\eps(s))\,ds\Big|
+
\sup\limits_{\tau\in\Delta_j}\Big|\frac1\eps\int_{\tau_j}^{\tau}(\theta(I^\eps(s)-\theta(I^\eps(\tau_j))\,ds\Big|\\
&\quad+
\sup\limits_{\tau\in\Delta_j}\Big|\int_{\tau_j}^{\tau}\Psi^\varphi(I^\eps(s),\varphi^\eps(s))\,d\beta(s)\Big|
\leq N^2\eps^\frac14+2 \eps^\frac14. %\leq CN^2\eps^\frac14.
\end{split}
\ee
{\bf Step 4: } Consider now a random collection of  intervals $\Delta_j= [\tau_j, \tau_{j+1})$
 such that
 \be\label{typ}
 \text{
 $I^\eps(\tau_j)\not\in \mathcal{A}_N^\delta$ (cf. event $\mathcal{E}_1$),\;
  $\zeta_j\leq \eps^\frac14$  (cf. event $\mathcal{E}_2$)
{
   and $ \sup\limits_{\tau\in\Delta_j}|I^\eps(\tau)|\leq \nu(M)$.}}
   \ee
We call these intervals {\it typical}.
 Clearly,
\be\label{om_typ}
\text{ if $\Delta_j^\omega$ is typical, then $\omega \in \mathcal{M}_j$.   }
\ee
For a typical interval $\Delta_j$  we will estimate the following quantity:
\[ \begin{split}
%\sup\limits_{t_j\leq t\leq t_{j+1}}
&\Big|\int_{\tau_j}^{\tau_{j+1}}\big(P^I(I^\eps(s),\varphi^\eps(s))-\langle P^I\rangle(I^\eps(s))\big)\,ds\Big| \\
&\leq\Big|\int_{\tau_j}^{\tau_{j+1}}\big(P^I(I^\eps(s),\varphi^\eps(s))-P^I(I^\eps(\tau_j),
\varphi^\eps(\tau_j)+{\textstyle \frac{s-\tau_j}\eps}\theta(I^\eps(\tau_j)))\big)\,ds\Big|\\
&+\Big|\int_{\tau_j}^{\tau_{j+1}}\big(P^I(I^\eps(\tau_j),
\varphi^\eps(\tau_j)+{\textstyle \frac{s-\tau_j}\eps}\theta(I^\eps(\tau_j)))-
\langle P^I\rangle(I^\eps(\tau_j))\big)\,ds\Big| \\
&+
\Big|\int_{\tau_j}^{\tau_{j+1}}\big(
\langle P^I\rangle(I^\eps(\tau_j))-\langle P^I\rangle(I^\eps(s))\big)\,ds\Big|=J_1+J_2+J_3.
\end{split}
\]
Due to \eqref{first_interv}, \eqref{second_interv} and \eqref{om_typ},
$$
J_1\leq CN^2\eps^\frac54+CN^3\eps^\frac54=C N^3\eps^\frac 54,\qquad J_3\leq C N^2\eps^\frac 54.
 $$
 Considering \eqref{def_AN} and Lemma~\ref{l_small_meaA} we find that
 \[ \begin{split}
 J_2&=\Big|\int_{\tau_j}^{\tau_{j+1}}\big(P^I(I^\eps(\tau_j),
\varphi^\eps(\tau_j)+{\textstyle \frac{s-\tau_j}\eps}\theta(I^\eps(\tau_j)))-
\langle P^I\rangle(I^\eps(\tau_j))\big)\,ds\Big| \\
&=\Big|\eps N \frac1N\int_0^{N}\big(P^I(I^\eps(\tau_j),
\varphi^\eps(\tau_j)+\xi \theta(I^\eps(\tau_j)))-
\langle P^I\rangle(I^\eps(\tau_j))\big)\,d\xi\Big|\leq \delta\eps N.
\end{split}
\]
Thus for a typical interval $\Delta_j$ we have:
\begin{equation}\label{esti_typica}
\Big|\int_{\tau_j}^{\tau_{j+1}}\big(P^I(I^\eps(s),\varphi^\eps(s))-\langle P^I\rangle(I^\eps(s))\big)\,ds\Big|
\leq CN^3\eps^\frac54+\delta {\eps N}.
\end{equation}

{\bf Step 5: }
We
 introduce the third good event {
 $$\mathcal{E}_3=\{\omega\in\Omega\,:\, \sup\limits_{0\leq\tau\leq T} |I^\eps(\tau)|\leq\nu(M)\} $$
 }
  (corresponding to the third condition in the definition \eqref{typ}   of a typical interval).   Due to (4) of Assumption \ref{assumption-1},
  {
$$
  \mathbf{P} (\mathcal{E}^c_3) \to 0 \quad \hbox{as }M\to\infty,
$$
}
and the rate of convergence is independent of $\eps$.
Finally we set
$$
 \mathcal{E} =  \mathcal{E}_1 \cap  \mathcal{E}_2\cap  \mathcal{E}_3.
$$
Then by \eqref{after_cheb} and  \eqref{est_tildeL}
{
$$
\PP  (\mathcal{E}^c )\le  \PP  (\mathcal{E}^c_1) +\PP  (\mathcal{E}^c_2) +\PP  (\mathcal{E}^c_3)  \le
(\alpha_N^\delta)^{1/2} +C_3 TN^{-1/2} +\PP  (\mathcal{E}^c_3) =: \beta^\delta_{N,M},
$$
where $\beta^\delta_{N,M} \to 0$ as $N, M\to\infty$, for each $\delta>0$.
Note that although the sets $\mathcal E$ and  $\mathcal E_j$ depend on $\eps$, $\delta$ and $N$,  the upper bound
$\beta_{N,M}^\delta$ for $\PP  (\mathcal{E}^c )$ is independent of $\eps$.
}

Denote
$$
\mathcal{I}^\eps:=\sup\limits_{0\leq \tau\leq T}\Big| \int_0^\tau \big(P^I(I^\eps(s),\varphi^\eps(s))-
\langle P^I\rangle(I^\eps(s))\big)\,ds\Big|
$$
and %  for any $j$, $0\leq j\leq {j_N}$,
$$
\mathcal{I}_j^\eps:=\Big| \int_{\tau_j}^{\tau_{j+1}} \big(P^I(I^\eps(s),\varphi^\eps(s))-
\langle P^I\rangle(I^\eps(s))\big)\,ds\Big|, \qquad 0\leq j\leq {j_N}.
$$
Apart from $\mathcal{I}^\eps$ let us consider $\mathcal{I}^\eps_N$, defined as follows:
$$
\mathcal{I}^\eps_N:=\sup\limits_{\tau= \tau_1, \dots, \tau_{j_N}}\Big| \int_{\tau_0}^\tau \big(P^I(I^\eps(s),\varphi^\eps(s))-
\langle P^I\rangle(I^\eps(s))\big)\,ds\Big|.
$$
Since  $\mathcal E \subset  \mathcal E_3$, then for $\omega \in  \mathcal E$ the norms of vector fields $P^I(I^\eps(s), \f^\eps(s))$ and $\lan P^I\ran(I^\eps(s))$ are
bounded by $M$ (we recall \eqref{nu_Lip}). Therefore,
\be\label{67}
| \mathcal{I}^\eps - \mathcal{I}^\eps_N| \le C N^2 \eps \quad \text{for} \quad \omega \in \mathcal{E}.
\ee

Conditions (3) and (4) of Assumption \ref{assumption-1} imply that
{
\begin{equation}\label{determ_esttt}
\mathbf{E} (\mathcal{I}^\eps_N)^2, \
\mathbf{E} (\mathcal{I}^\eps)^2   \leq CT^2 \quad\hbox{and}\quad \mathbf{E} (\mathcal{I}_j^\eps)^2   \leq C\eps^2 N^3.
\end{equation}
Due to \eqref{determ_esttt} and \eqref{67},
\be\label{66}
\Upsilon^\eps = \mathbf{E}\,\mathcal{I}^\eps =
 \mathbf{E}\big(\mathbf{1}_{\mathcal{E} }\,\mathcal{I}^\eps \big)+
\mathbf{E}\big(\mathbf{1}_{
\mathcal{E}^c }\,{\mathcal{I}^\eps} \big)
\leq \mathbf{E} \big( \mathbf{1}_{\mathcal{E} }\,{\mathcal{I}^\eps_N} \big)+CN^2 \eps+(CT^2)^{\frac12} (\beta_{N,M}^\delta)^{1/2}.
\ee
}
{\bf Step 6: }  Let us estimate $\mathbf{E} \big( \mathbf{1}_{\mathcal{E} }\,{\mathcal{I}^\eps_N} \big)$.  Denoting
$$
\mathcal J =\mathcal J^\eps = \{ j: \Delta_j\ \text{ is typical} \}
$$
we have that
\be\label{666}
\mathbf{E} \big( \mathbf{1}_{\mathcal{E} }\,{\mathcal{I}^\eps_N} \big) \le
\mathbf{E} \big(  \sum_{j \in \mathcal J}
\mathcal{I}^\eps_j \big) + \mathbf{E} \big( \mathbf{1}_{\mathcal{E}}\, \sum_{j \notin \mathcal J}
\mathcal{I}^\eps_j \big) =: S_1 +S_2.
\ee
From the definitions of sets $\mathcal{E}_1, \mathcal{E}_2$ we see that for any $\omega\in \mathcal{E}$,
$$
\# \mathcal{J}^c \le (\alpha^\delta_N)^{1/2} (\eps N)^{-1} + \eps^{-7/8}.
$$
Since for each $\omega \in \mathcal{E}_3$ all $\mathcal{I}^\eps_j$  are trivially bounded by {
 $C\eps NM$, } then
by the above estimate
{
$$
S_2 \le C M \big( (\alpha^\delta_N)^{1/2} +N \eps ^{1/8}\big).
$$
}
\smallskip

Due to \eqref{esti_typica}, for each $\omega$
$$
 \sum_{j \in \mathcal J} \mathcal{I}^\eps_j  \le j_N \big( C N^3 \eps^{5/4} + \delta\eps N\big) \le
 C_1 \big(  N^2 \eps^{1/4} + \delta\big) .
$$
So  $$S_1\le C_1 \big(  N^2 \eps^{1/4} + \delta\big).$$

By \eqref{66}, \eqref{666}  and the estimates for $S_1$ and $S_2$,
{
$$
\Upsilon^\eps \le  C'\beta_{N,M}^\delta  +  C_1 \big(  N^2 \eps^{1/4} + \delta\big) +
CM \big( (\alpha^\delta_N)^{1/2} +N \eps ^{1/8}\big) + CN^2 \eps.
$$
Now for any $\delta>0$ we choose $M=M(\delta)$ so big that
$
C'\beta_{N,M'}^\delta  \leqslant \delta$ for all $M'\geqslant M$.
Then for this $M$ we choose sufficiently large $N$ so that
$ C  M(\alpha^\delta_N)^{1/2} \le \delta.$
For these $M$ and $N$ we have
$$
\Upsilon^\eps \le \delta  +  C_1   N^2 \eps^{1/4} +  C_1 \delta  + C N \eps ^{1/8}  +  CN^2 \eps.
$$
}
Taking $\eps$ small enough we achieve that $\Upsilon^\eps \le (C_1+2)\delta $.
Since $\delta$ is an arbitrary positive number, then \eqref{lim_rel_first} follows.
\end{proof}
With exactly the same proof  we also have
\begin{lemma}\label{lemma-key-2}
The following convergence  holds as ${\eps}\to0$:
\[
\mathbf{E}\Big\{\sup\limits_{0\leq \tau\leq T}\Big| \int_0^\tau \big( \Psi^I(I^\eps(s),\varphi^\eps(s))
(\Psi^I(I^\eps(s),\varphi^\eps(s)))^t-
\langle a^I\rangle(I^\eps(s))\big)\,ds\Big|  \Big\}\to0.
\]\end{lemma}

\section{The averaging theorem}
\label{s_lp}
In the section we show that  the limiting laws of the family  $\{\cD(I^{\eps}(\cdot))\}$, as ${\eps}\to0$,
 are solutions of the martingale problem for  the averaged equation \eqref{averaged-eq-1}, thus are weak
 solutions of the later. We begin with the corresponding definition.
Let us introduce a natural filtered measurable space for the problem we consider
\begin{equation}\label{natural-m-s}(\tilde\Omega,\mathcal{B}_T,\{\mathcal{B}_\tau,0\leqslant \tau\leqslant T\}),\end{equation}
where $\tilde\Omega$ is  the Banach space $C([0,T];\mathbb R^d)=\{a(\tau),\tau\in[0,T]\}$, $\mathcal{B}_T$ is its Borel $\sigma$-algebra and a $\sigma$-algebra
  $\mathcal{B}_\tau \subset  \mathcal{B}_T$ is  generated by  the restriction mapping
 $ %C([0,T];\mathbb R^d)
 \tilde\Omega \ni
a(\cdot) \mapsto a(\cdot) \mid_{[0,\tau]}.$
  Consider the process
% \begin{equation}\label{drift-martingale}
$$
 N^I(\tau; a):=a(\tau)-\int_0^\tau\langle P^I\rangle(a(s))ds,\; a\in\tilde\Omega,\tau\in[0,T],
 $$
 and
 note that for any $0\leqslant \tau\leqslant T$, $N^I(\tau; \cdot)$ is a $\cB_\tau$-measurable continuous functional on~$\tilde\Omega$.

 \begin{definition} (see \cite{SV, brownianbook}).
 A measure $Q$ on the space \eqref{natural-m-s} is called a solution of the martingale problem for  equation \eqref{averaged-eq-1}, if
$a(0)=I_0$ $Q$-a.s. and

1) the process $\{N^I(\tau; a)\in\mathbb{R}^d,\tau\in [0,T]\}$ is a vector-martingale on the filtered space \eqref{natural-m-s} with respect to the measure $Q$;

2) for any $k,j=1\dots,d$ the process $N^I_k(\tau; a)N_j^I(\tau; a)-\int_0^\tau \langle a^I\rangle_{kj}(a(s))ds$ is a martingale on the space \eqref{natural-m-s} with respect to the measure $Q$.
 \end{definition}
For each $\eps>0$ we define a probability measure $Q^\eps$ on $(\tilde\Omega,\mathcal{B}_T)$ as the law of $\{I^\eps(\cdot)\}$ and denote by
$\mathbf{E}^{Q^\eps}$  the corresponding expectation.
According to Lemma~\ref{l_compa} the family $\{Q^\eps\}$
is tight in $\mathcal{P}(C([0,T],\mathbb{R}^d)$. Take a  sequence ${\eps}_j\to0$ such that
\begin{equation}\label{weak-limit-1}
Q^{{\eps}_j}\strela Q^0,\; \text{ as }{\eps}_j\to0\;  \text { in } \mathcal{P}(C([0,T],\mathbb{R}^d)).
\end{equation}
\begin{lemma} The probability measure $Q^0$ above is a solution of the martingale problem for the averaged equation \eqref{averaged-eq-1}.
\end{lemma}
\begin{proof} For any $s\in[0,T]$,
let $\Phi(\cdot)$ be a bounded continuous
functional defined on $C([0,s];\mathbb R^d)$. Then, for any $\tau\in[0,T]$ such that $0\leq s\leq \tau\leq T$, we have
\[ \begin{split}
&\mathbf{E}^{Q^0}\Big\{\big(N^I(\tau;a)-N^I(s;a)\big)\Phi(a([0,s]))\Big\}\\
&=\lim\limits_{\eps_j\to0}
\mathbf{E}^{Q^{{\eps}_j}}\Big\{\big(N^I(\tau;a)-N^I(s;a)\big)\Phi(a([0,s]))\Big\} \\
&=\lim\limits_{\eps_j\to0}\mathbf{E}\Big\{\Big(I^{\eps_j}(\tau)-I^{\eps_j}(s)-\int_s^\tau  \langle P^I\rangle(I^{\eps_j}(u))\,du\Big)\Phi(I^{\eps_j}([0,s]))\Big\}\\
&=\lim\limits_{\eps_j\to0}\mathbf{E}\Big\{\Big(\int_s^\tau  \big(P^I(I^{\eps_j}(u),\varphi^{\eps_j}(u))-\langle P^I\rangle(I^{\eps_j}(u))\big)\,du\Big)\Phi(I^{\eps_j}([0,s]))\Big\} \\
&+
\lim\limits_{\eps_j\to0}\mathbf{E}\Big\{\Big(I^{\eps_j}(\tau)-I^{\eps_j}(s)-\int_s^\tau P^I(I^{\eps_j}(u),\varphi^{\eps_j}(u))\,du\Big)\Phi(I^{\eps_j}([0,s]))\Big\},
 \end{split}
 \]
 where to get the first equality we
have  used the fact  that the r.v.'s  $\|I^\eps\|^{q}_{C([0,T])}$
are uniformly integrable, which is guaranteed by (3) and (4) of Assumption \ref{assumption-1} (cf. \cite[Lemma~4.4]{GHKu}).
The first limit on the r.h.s of the
last equality  vanishes due to Lemma~\ref{lemma-key-1} and the second one vanishes because
$I^\eps(\tau)-\int_0^\tau P^I(I^\eps(u),\varphi^\eps(u))\,du$ is a martingale. Therefore,
 $$
\mathbf{E}^{Q^0}\Big\{\Big(N^I(\tau;a)-N^I(s;a)\Big)\Phi(a([0,s]))\Big\}=0,
$$
for each $\Phi$ as above. So
the process $N^I(\tau;a),\; \tau\in[0,T] $, is a martingale with respect to the measure
 $Q^0$ and  filtration $\{\mathcal{B}_\tau\}$.

Arguing in the same way and evoking Lemma~\ref{lemma-key-2} we conclude that the processes  $N^I_k(\tau; a)N_j^I(\tau; a)-\int_0^\tau \langle a^I\rangle_{kj}(a(s))ds$, $\tau\in[0,T]$, $k,j=1,\dots,d$, also
are  $(Q^0,\mathcal{B}_\tau)$-martin\-gales. Hence, the assertion of the lemma is proved.
\end{proof}

 As we have discussed before, since
  the drift term $\langle P^I\rangle$ and the dispersion matrix $\llangle \Psi^I \rrangle$ in \eqref{averaged-eq-1} are locally Lipschitz with respect $I$,  then by the Yamada-Watanabe theorem (see \cite[Section~5.3.D]{brownianbook})   the  just obtained  solution of the
  martingale problem for \eqref{averaged-eq-1} is unique.
  So,  in particular,  the limit $Q^0$ in \eqref{weak-limit-1} does not depend on the choice of ${\eps}_j$,
   and the whole family $Q^\eps$ converges as $\eps\to0$ to the measure
 $Q^0$. Thus we have obtained the following result.

 \begin{theorem}\label{t_main02} Under Assumption \ref{assumption-1}, for any $(I_0,\varphi_0)\in\mathbb{R}^d\times\mathbb{T}^n$ we have
 $$
 \cD(I^\eps(\cdot))\strela Q^0\; \text{ as }{\eps}\to0\;  \text { in } \mathcal{P}(C([0,T],\mathbb{R}^d)),
 $$
 where $Q^0$ is the law of a weak solution $I^0(\tau)$ for problem  \eqref{averaged-eq-1}.
  Moreover, for  $q_0$ and $C(\cdot)$  as in \eqref{apriori-1} we have
   $\mathbf{E}\sup_{\tau\in[0,T]}|I^0(\tau)|^{2q_0}\leqslant C(|I_0|,T)$.
   \end{theorem}
 The last assertion follows directly from the Skorokhod representation theorem (see \cite[Section 6 ]{Bill} and Fatou's Lemma, cf. \cite[Remark 4.8]{GHKu}.
 \begin{remark}\label{remark-main-thm}
1) It is straightforward to see that the statement of Theorem~\ref{t_main02}  remains true with the same proof if in \eqref{ini_cond_det} the initial data $(I_{0,{\eps}}, \varphi_{0,{\eps}})$ depends
on $\eps$ and converges to a limit  $(I_0,\varphi_0)$ as $\eps\to0$.

2) The result in Theorem~\ref{t_main02} admits an immediate generalization to the case when the initial data $(I_0,\varphi_0)$ is a random variable. Cf.  Amplification \ref{amp_ran_ini}.

3) A local in $I$  version of Theorem \ref{t_main02} was earlier proved in \cite{FrWe03} by another method.

\end{remark}

  \section{Stationary solutions}\label{s_stationary_nonl}
  The goal of this section is to characterize the asymptotic behaviour  of   stationary solutions of  equations
   \eqref{ori_scaled}   and to find out their relation with a stationary solution of  averaged equation \eqref{averaged-eq-1}.
   We recall that  a solution $(I^\eps(\tau),\varphi^\eps(\tau))$, $\tau\geqslant0$, of  equation \eqref{ori_scaled}  is called stationary if
   there exists $\nu^\eps\in\mathcal{P}(\mathbb{R}^d\times\mathbb T^n)$ such that  $\mathcal{D}\big((I^\eps(\tau),\varphi^\eps(\tau))\big)=\nu^\eps$ for all $\tau\geqslant0$.
   Then the   measure $\nu^\eps$ is called a stationary measure for equation \eqref{ori_scaled}.
   A stationary solution and  stationary measure of the averaged equation \eqref{averaged-eq-1} are defined in the same way.

  Throughout  this section  a strengthened  version of Assumption \ref{assumption-1} is imposed
  on the system \eqref{ori_scaled}.  Namely, we assume that
  \begin{assumption}\label{assumption-mixing-1} i) Items (1)-(3) of Assumption \ref{assumption-1} hold true.

  ii) For each $\eps\in (0,1]$ and any $(I_0,\varphi_0)\in\mathbb{R}^d\times\mathbb T^n$  problem
   \eqref{ori_scaled}-\eqref{ini_cond_det} has a unique strong solution $\big(I^{\ep},\varphi^\eps\big)(\tau;I_0,\varphi_0)$,
 $\tau\in[0,+\infty)$,    satisfying
  \begin{equation}\label{v-upper-mix_sec10}
  {\EE}\sup_{ {T'}\leqslant\tau\leqslant  {T'}+1}|I^{\ep}(\tau;I_0,\varphi_0)|^{2q_0}\leqslant C_{q_0}(|I_0|),
  \end{equation}
for each $ {T'}\geqslant0$ and some number $q_0>(q\vee 2)$.

iii)
  Equation \eqref{ori_scaled} is  mixing. So it has a stationary weak solution
$\big(I^\eps_{\rm st}(\tau),\varphi^\eps_{\rm st}(\tau)\big)$ such that $\cD \big(I^\eps_{\rm st}(\tau),\varphi^\eps_{\rm st}(\tau)\big) \equiv   \nu^\eps\in \cP(\mathbb R^d\times\mathbb T^n)$,  and
\be\label{mixing_sec10}
\mathcal{D}\big(I^{\ep}(\tau;I_0,\varphi_0),\varphi^{\ep}(\tau;I_0,\varphi_0)\big)\rightharpoonup\nu^{\ep}\text{ in }\mathcal{P}(\mathbb{R}^d\times\mathbb T^n)\; \;\text{ as }\tau\to+\infty,
\ee
for each $(I_0,\varphi_0)$.
  \end{assumption}

  Under   Assumption \ref{assumption-mixing-1} equation \eqref{ori_scaled} defines in $\R^d\times\T^n$ a Markov process
  with the transition probability
  $
  \Sigma_\tau^{\eps}(I,\varphi) \in \cP(\R^d\times\T^n)$, $\tau \ge0$, $(I,\varphi)\in \R^d\times\T^n,
  $
  where
  $
  \Sigma^{\eps}_\tau(I,\varphi) = \cD \big(I^\eps(\tau; I,\varphi),\varphi^\eps(\tau; I,\varphi)\big);
  $
  e.g. see \cite[Section~5.4.C]{brownianbook}.

  As a consequence of  \eqref{v-upper-mix_sec10} and \eqref{mixing_sec10}
  by a straightforward argument (e.g. see in \cite{GHKu}
  Lemma~5.3 and its proof),   we obtain
   that the stationary solution $(I^\ep_{\rm st},\varphi^\ep_{\rm st})$ satisfies the following estimate
   \begin{equation}\label{v-upper-stat_sec10}
  {\EE}\sup_{ {T'}\leqslant\tau\leqslant  {T'}+1}|I_{\rm st}^{\ep}(\tau)|^{2q_0}\leqslant C_{q_0}(0).
  \end{equation}

Using  \eqref{v-upper-stat_sec10} we derive from the first equation in \eqref{ori_scaled}
that for any $N\in\Nn$,
     the collection of measures
  $\{\cD(I_{\rm st}^{\ep}|_{[0,N]}),0<{\ep}\leqslant1\}$ is tight. Choosing for each $N$
 a sequence ${\ep}^{(N)}_l\to0$ %(depending on $N$)
 such that
  \[\cD(I_{\rm st}^{{\ep}^{(N)}_l}|_{[0,N]})\rightharpoonup Q^0\text{ in } \mathcal{P}\big(C([0,N],\mathbb{R}^d)\big)\]
  and  applying the diagonal procedure  we conclude that for a subsequence $\{{\ep}_l\}$ the  relation
  $
  \cD (I_ {\rm st}^{{\ep}_l} )\strela Q^0
  $
holds in $\cP(X)$,  % (see \eqref{X} for the definition of $X$).
where $X$ is the complete separable metric space
  $
  X=C([0,\infty), \R^d)
  $
  with the distance
%  \be\label{X}
$$
  \text{dist} (a_1, a_2) = \sum_{N=1}^\infty 2^{-N} \frac{ \| a_1 -a_2\|_{C([0,N], \R^d)}} {1+ \| a_1 -a_2\|_{C([0,N], \R^d)}},\qquad a_1, a_2\in X\,.
  $$

Denote $\mu^\ep_I(\tau)=\cD (I^{{\ep}}_{st}(\tau)) = I\circ \nu^\eps$.
%   Since $a^{\ep}(0)=v_{st}^{\ep}(0)$,
   Then $\mu_I^{{\ep}_l}(0)\rightharpoonup\mu^0:=Q^0|_{\tau=0}$.
   Let $I^0(\tau)$ be a solution of equation \eqref{averaged-eq-1} with an initial condition $I_0$, distributed as
   $\mu^0$. Then, by Remark~\ref{remark-main-thm},
   \be\label{hrum_sec10}
   \cD (I_{\rm st}^{{\ep}_l} (\cdot) )\strela  \cD (I^{0} (\cdot) )\quad \hbox{in } X,
   \ee
   and for any $\tau\ge0$
   we have $Q^0(\tau)= \cD (I^{0} (\tau))=\lim\limits_{\eps_l\to0}  \cD I(_{\rm st}^{{\ep}_l} (\tau)) =
   \lim\limits_{\eps_l\to0}  \cD (I_{\rm st}^{{\ep}_l} (0))=\mu^0$.
%   In particular, $\cD I^0(0) = \mu^0$.
%By the same arguments as those in the proof of Proposition \ref{p_5.4} in \cite{GHKu}
We obtained the following statement:

  \begin{proposition}\label{p_5.4_sec10}
   The process $I^0$ is a stationary weak solution of the averaged equation \eqref{averaged-eq-1}
 %   \begin{equation}\label{lim_proce2_sec10}
 %$$
% dI^{0}(\tau)=\langle P^I\rangle(I^{0}(\tau))\,dt+q^I(I^{0}(\tau))\,d\tilde W(\tau),
% ,\qquad I^{\rm h}(0)=I_0,
%  \end{equation}
%$$
%   \eqref{effective-equation},
and  $\cD(I^0(\tau))\equiv\mu^0$, $\tau\in[0,\infty)$. In particular,
  any limiting point of the collection of measures
  $\{\mu_I^{\ep}:=\cD (I^{{\ep}}_{\rm st}(\tau))\}$ as ${\ep}\to0$ is a stationary measure of  the
  averaged  equation \eqref{averaged-eq-1}.
  If  the latter  is mixing, then  its stationary measure is unique, and so
   convergence \eqref{hrum_sec10}    holds as $\eps\to0$.
  \end{proposition}

 % The last assertion holds that if equation \eqref{averaged-eq-1} is mixing, then  the limiting measure $\mu^0$ is
 % defined in a unique way.
  \medskip

  \subsection{Uniform convergence in the averaging theorem} \label{ss_average_uniform}
  To describe quantitively   the weak convergence of measures in Theorem~\ref{t_main02}
   we introduce the dual-Lipschitz distance.
\begin{definition}\label{d_dual-Lip}
Let $M$ be a      complete and separable metric space. For any two measures $ \mu_1,\mu_2\in\cP(M)$ we define
the dual-Lipschitz distance between them as
\[
\|\mu_1-\mu_2\|_{L,M}^*:=\sup_{f\in C_b (M),\,  |f|_L\leqslant1}\Big|\langle f ,\mu_1\rangle-\langle f ,\mu_2\rangle\Big|
=\sup_{f\in C_b (M),\,  |f|_L\leqslant1}\Big(\langle f ,\mu_1\rangle-\langle f ,\mu_2\rangle\Big)
\le2,
 \]
   where  $|f|_L= |f|_{L,M}=
    \,$Lip$\,f+\|f\|_{C(M)}$.
 \end{definition}

 In this definition and below we denote
    \be\label{recall}
   \langle f ,\mu\rangle :=\int_M f(m)\mu(dm).
   \ee
 The dual-Lipschitz  distance converts  $\mathcal{P}(M)$ to a complete metric space and
    induces on it  a  topology, equivalent to the  weak convergence of measures, e.g. see \cite[Sec\-tion ~1.3]{Dud}.

  If equation  \eqref{averaged-eq-1}   is mixing with  a  quantitive  property as in the following assumption,  then the convergence in Theorem~\ref{t_main02} is uniform in time with respect to the dual-Lipschitz distance.
\begin{assumption}\label{assumption-mixing-2}
1) i) and (ii) of Assumption \ref{assumption-mixing-1} hold true.

2) The averaged equation \eqref{averaged-eq-1} is mixing with a unique stationary measure $\mu^0$. Moreover, for each $M>0$ and any $I_0, I_1\in \bar B_M(\mathbb{R}^d)$ the laws of solutions for  \eqref{averaged-eq-1} with these initial data satisfy
\[\|\cD(I^0(\tau; I_0))-\cD(I^0(\tau; I_1))\|_{L,\mathbb{R}^d}^*\leqslant \frak f_M(\tau),\]
where $\frak f_M$ is a positive
continuous function of $(M,\tau)$ which goes to zero when $\tau\to\infty$ and is non-decreasing in $M$.
\end{assumption}

\begin{theorem}\label{t_mixAP}
 Under Assumption \ref{assumption-mixing-2}, in the setting of Theorem~\ref{t_main02}   for any initial data
 $(I_0,\varphi_0)\in\mathbb{R}^d\times\mathbb{T}^n$ we have
\[\lim_{\eps\to0}\sup_{\tau\geqslant0}\|\cD(I^\eps(\tau; I_0,\varphi_0))-\cD(I^0(\tau; I_0))\|_{L,\mathbb{R}^d}^*=0.\]
\end{theorem}

Concerning the theorem's proof see  Subsection \ref{ss_10.1}.  Assumptions \ref{assumption-mixing-1} and \ref{assumption-mixing-2} allow
for an instrumental sufficient condition, cf. below Proposition~\ref{sufficient1}.

\section{Perturbations of integrable equations in $\R^{2n}$}\label{s_Birk_int}

In this section we study  diffusive perturbations of integrable equations in $\R^{2n}$
in the framework of previous sections.
By bold italic letters we  denote various vectors in $\R^2$ (regarded as column-vectors).

Let us consider a perturbed   integrable equation \eqref{3}
 for a vector $v=(\bv_k,\; k=1,\dots,n)\in\mathbb{R}^{2n}$ and write it  in the slow time $\tau = \eps t$:
\begin{equation}\label{v-equation1}
\begin{cases}
d\bv_k=\ep^{-1}W_k(I)\bv_k^{\bot}d\tau+{\mathbf{P}}_k(v)d\tau+\sum_{j=1}^{n_1}B_{kj}(v)d\bbeta_j(\tau),\quad
k=1,\dots,n,\\
v(0)=v_0\in\mathbb{R}^{2n}.
\end{cases}
\end{equation}
Here $\bv_k=\left(\begin{array}{c}v_k\\v_{-k}\end{array}\right)\in\mathbb{R}^2$,
$\bv_k^{\bot}=\left(\begin{array}{c}-v_{-k}\\v_{k}\end{array}\right)$,
 $B_{kj}(v)$ are $2\times2$ matrix functions,  $I=(I_1,\dots,I_n)$ with $I_k=\frac{1}{2}\|\bv_k\|^2$ and
  $\bbeta_j(\tau)=\left(\begin{array}{c}\beta_j(\tau)\\\beta_{-j}(\tau)\end{array}\right)$ are
 independent standard Brownian motions in $\mathbb{R}^2$, defined on a filtered probability space
 \be\label{OmOm}
 \big( \Omega, \cF, \{ \cF_\tau\}, \PP\big).
 \ee
  The unperturbed part  \eqref{Birk}  of  equation  \eqref{v-equation1}
   is integrable, and  the functions $I_k(v), k=1,\dots, n$, are its integrals of motion.
  If $W(I) = \nabla h(I)$ for some $C^1$-function $h$, then  equation \eqref{Birk}   is
  Hamiltonian with the Hamiltonian function  $h(I(v))$, and then  it is  is called
  integrable  in the sense of   Birkhoff.

   As in  previous sections, the initial data $v_0$ can be deterministic or random, and a solution of problem
   \eqref{v-equation1} will be denoted
 $
 v^\ep(\tau;v_0) = (\bv_k^\ep(\tau;v_0), k=1, \dots n),
 $
 or simply $v(\tau)$.  We will  focus on the deterministic case,
 always assuming  that equation \eqref{v-equation1} satisfies the following assumption.

\begin{assumption}\label{assumption-v-3}
(1) The Lebesgue measure of $I\in \mathbb R^n_+$ for which  components of the vector
 $W(I)=:(W_k(I),k=1,\dots,n)$ are  rationally dependent is equal to zero. That is,
 $\mathcal{L}\big(\cup_{k\in\mathbb{Z}^n\setminus \{0\} }\{I\in\mathbb{R}_+^n: W(I)\cdot k=0\})=0$.

(2) The $2n\times 2n$ diffusion matrix $S(v)=B(v)B(v)^t$, where $B(v)=(B_{kj}(v))$,
satisfies the uniform ellipticity condition. That is, there exists $\lambda>0$ such that
\be\label{ellipt}
\lambda \|\xi\|^2\leqslant S(v)\xi\cdot\xi\leqslant\lambda^{-1} \|\xi\|^2,\quad \forall v, \xi\in\mathbb{R}^{2n}.
\ee

(3) There exists $q>0$ such that $W(I)\in \text{Lip}_q(\mathbb{R}^n,\mathbb{R}^n)$, $P(v):=({\mathbf{P}}_k(v),k=1,\dots,n)\in\text{Lip}_{q}(\mathbb{R}^{2n},\mathbb{R}^{2n})$ and $B(v)\in \text{Lip}_q(\mathbb{R}^{2n},M_{2n\times 2n_1})$ (we recall \eqref{nota}).

(4) There exists $T>0$ such that for every $v_0\in\mathbb{R}^{2n}$
 equation \eqref{v-equation1} has a unique strong solution $v^{\eps}(\tau; v_0)$, $\tau\in[0,T]$, equal $v_0$ at $\tau=0$. Moreover, there exists $q_0>(q\vee4)$ such that
\begin{equation}\label{apriori-v-2}
\mathbf{E}\sup_{\tau\in[0,T]}\|v^{\eps}(\tau;v_0)\|^{2q_0}\leqslant C_{q_0}(\|v_0\|,T),\quad \forall {\eps}\in(0,1],
\end{equation}
where $C_{q_0}(\cdot)$ is a non-negative continuous function on $\mathbb{R}_+^2$, non-decreasing in both arguments.
\end{assumption}

\begin{remark}\label{r_111} 1) In the assumption above, (1) is  Anosov's condition  (see Remark~\ref{r_Anosov}), and it holds,
in particular,   if $W(I)$ is a constant vector with rationally independent coefficients.
  Equations \eqref{v-equation1} with constant frequency vectors $W$ and without
assuming that its components are  rationally independent  are examined in \cite{GHKu}, and for the case $W=\,$const the results,
given below in Sections~\ref{s_Birk_int}-\ref{s_9}, are special cases of more general theorems in that work. But equations
\eqref{v-equation1} with non-constant frequency vectors $W(I)$ are significantly more complicated then those with
$W=\,$const.

2) Item (4) holds if assumptions
(2), (3) are valid and if the coefficients of equation \eqref{v-equation1} are globally Lipschitz
 (cf. Proposition~\ref{l_suff_assum_i4}),  or if the vector field $P$ is coercive, see below Proposition~\ref{sufficient1}.
\end{remark}

\begin{example}\label{ex_Fourier}
In statistical physics they  often examine stochastic perturbations of chains of nonlinear oscillators
\be\label{system}
\ddot q_k = -Q(q_k),
% -(\alpha^2 q_k + \beta q_k^3),
 \qquad k=1, \dots, n,
\ee
where $Q$ is a polynomial $Q^0(q) = \alpha q+ \beta q^3$,  $\alpha, \beta>0$, or more generally,  is a smooth function of the form
\be\label{Q}
Q(q) = Q^0(q) + O(q^4).
\ee
System \eqref{system} may be re-written in the Hamiltonian form as
\be\label{s1}
\dot q_k =  (\p/\p p_k)H, \quad \dot p_k =-  (\p/\p q_k)H,  \qquad k=1, \dots, n,
\ee
where
$
H= \sum\big( \frac12  p_k^2 + \int_0^{q_k} Q(l) dl).
$
E.g. see \cite{Ek} and \cite[Section 4]{BLR}. Under a suitable diagonal symplectomorphism
$$
\R^{2n} \to \R^{2n} , \quad (q,p) \mapsto v=(\bv_1,\dots, \bv_n),  \;\; \bv_j ={\mathbf F}(q_j, p_j),
$$
system \eqref{s1}  may be re-written in the form \eqref{Birk}, where $W_k=\omega(I_k)$, $k=1,\dots,n$,
and $\omega$ is a smooth function. See below
Appendix~\ref{app_B} concerning canonical transformations ${\mathbf F}: \R^2 \to \R^2$ with the required properties.
In the  just cited papers the
  vector field $P=(\bP_1, \dots , \bP_n)$ in \eqref{v-equation1}    is of the
type ``next neighbour interaction", i.e. $\bP_k= \bP_k(\bv_{k-1},  \bv_{k}, \bv_{k+1})$.
Our results apply to  corresponding equations \eqref{v-equation1}
with non-degenerate dispersion matrices $B$ if assumptions (1), (3) and (4) hold. The first two are easy to verify, while for a sufficient
condition for (4) see  Proposition~\ref{sufficient1}.  Below in Section~\ref{ss_10.2} we discuss a class of damped/driven
Hamiltonian systems which includes the equations we are now discussing if there the vector field $P$ is Hamiltonian with added damping and
the stochastic term is a diagonal random force.

 Usually  the systems in  works of
physicists correspond to  equations \eqref{v-equation1} with degenerate diffusion such that
only the $2\times2$-matrices $B_{11}$ and $B_{nn}$ are non-zero. To treat them the methods of our work
should be developed further. Still we mention that stochastic perturbations of linear equations \eqref{system}
with $Q(q) = \alpha q$, $\alpha>0$,
 and with degenerate (or non-degenerate) diffusion may be examined using the results of \cite{GHKu}. \qed
\end{example}
\smallskip

Let us consider  the action-angles
 mapping $\mathbb{R}^{2n}\to\mathbb{R}^{n}_+\times\mathbb{T}^{n}, \  v\mapsto (I,\varphi)(v)$, where
\be\label{ac_an}
\begin{split}
&I(v)=(I_k(v),k=1,\dots,n),\; \;\; \;\; I_k(v)=\frac{1}{2}\|\bv_k\|^2, \\
 &\varphi(v)=(\varphi_k(v), k=1,\dots,n), \; \;\;\;\;\varphi_k(v)=\text{Arg}(\bv_k)=\arctan\frac{v_{k}}{v_{-k}}
  \end{split}
  \ee
  and $\varphi_k:=0\text{ if }\bv_k=0$.  Then
  \be\label{inverse}
  \bv_k = \sqrt{2I_k}\,(\cos \varphi_k,  \sin \varphi_k), \quad k=1, \dots, n.
  \ee

 By It\^o's formula, if $v(\tau)$ is a solution of \eqref{v-equation1}, then the
 equations for  the actions $I_k(v)$ read
\begin{equation}\label{I-equation1}
dI_k=\bv_k^t{\mathbf{P}}_k(v)d\tau+\frac{1}{2}\sum_{j=1}^{n_1}\|B_{kj}(v)\|^2_{HS}d\tau+\sum_{j=1}^{n_1}\bv_k^tB_{kj}(v)d\bbeta_j(\tau), \;\;\; k=1,\dots,n,
\end{equation}
where  $\|\cdot \|_{HS}$ is the  Hilbert-Schmidt norm.
 Equations for $ \varphi_k(v)$, $k=1,\dots,n$,  hold if all $\bv_k$ are non-zero and read
\begin{equation}\label{phi-equation1}
d\varphi_k(\tau)={\eps}^{-1}W_k(I)d\tau+\Phi^1_k(v)d\tau+\sum_{j=1}^{n_1} \Phi^2_{kj}(v)d\bbeta_j(\tau),\; \;\;k=1,\dots,n,
\end{equation}
where
\[\Phi_k^1(v)=\Big(\nabla_{\bv_k}\arctan\big(\frac{v_k}{v_{-k}}\big)\Big)\cdot {\mathbf{P}}_k(v)+\frac{1}{2}\sum_{j=1}^{n_1}\text{Trace}\Big( B_{kj}(v)\Big(\nabla_{\bv_k}^2\arctan\big(\frac{v_k}{v_{-k}}\big)\Big) B_{kj}(v)^t\Big),\]
and
$\Phi^2_{kj}(v)=\Big(\nabla_{\bv_k}\arctan\big(\frac{v_k}{v_{-k}}\big)\Big)\cdot  B_{kj}(v).$
\smallskip

\begin{remark}
Note that in view of \eqref{inverse}, for $I_k$ near zero
the r.h.s of \eqref{I-equation1} is a H\"older-$\tfrac12$  function of $I_k$  and the r.h.s.
of \eqref{phi-equation1} has  a  strong singularity  when $I_k$ vanishes.
 Moreover, the dispersion part of \eqref{I-equation1} vanishes with $I_k$. Hence,
  system \eqref{I-equation1}+\eqref{phi-equation1} is singular and degenerated at the set
  $\bigcup_{k=1}^n\{(I,\varphi): I_k=0\}$.  %=\bigcup_{k=1}^d\{v: \bv_k=0\}$.
\end{remark}

As in  previous sections, we introduce
the averaged equations for $I_k(\tau)$, $k=1,\dots,n$, as
\begin{equation}\label{averaged-1}
\begin{split}
dI_k(\tau)&=F_k(I)
%\langle\bv_k^t{\mathbf{P}}_k\rangle(I)d\tau+\frac{1}{2}\langle\sum_j\|B_{kj}\|_{HS}^2\rangle(I)
d\tau+\sum_{j=1}^nK_{kj}(I)d\beta_j(\tau),\\
F_k(I) &= \langle\bv_k^t{\mathbf{P}}_k\rangle(I)
+\tfrac{1}{2}\langle\sum_j\|B_{kj}\|_{HS}^2\rangle(I),
\end{split}
\end{equation}
with the initial condition
\begin{equation}\label{initial-a1}
I(0)=I_0=I(v_0).
\end{equation}
The brackets $\langle\cdot \rangle$ signify the 	averaging in $\varphi$,  see \eqref{average},
and the dispersion matrix ${K}(I)=\big(K_{kj}(I)\big)_{1\leqslant k,j\leqslant n}$ is chosen to be the principal
symmetric square root  of the averaged diffusion matrix ${S}(I)=(S_{km}(I))_{1\leqslant k,m\leqslant n}$ of equation \eqref{I-equation1},
\begin{equation}\label{diffusion-a1}
S_{km}:=\big\langle \sum_{j=1}^{n_1}\bv_k^tB_{kj}B_{mj}^t\bv_{m}\big\rangle(I), \; k,m=1,\dots,n.
\end{equation}
So $K=K^t\ge0$ and
%\begin{equation}\label{dispersion-a1}
$$
\sum_{j=1}^nK_{kj}(I)K_{mj}(I)=S_{km}(I), \; k,m=1,\dots,n.
$$

\begin{remark}\label{remark-degenerate}Under (3) of  Assumption \ref{assumption-v-3}, the drift  and dispersion terms
 of \eqref{averaged-1} are only H\"older-$\tfrac12$ smooth with respect to $I_k$ (and Lipschitz with respect to $\sqrt{I_k}$), $k=1,\dots,n$.
 Moreover, the dispersion term vanishes at $I_k=0$. If we strengthen the assumption by assuming
  that the vector-field $P$  and the dispersion matrix $B$ are
   $C^2$-smooth, then by a straightforward application of Whitney's theorem (see in \cite{Whit} Theorem~1 and
   the last remark of the paper),  the drift term in \eqref{averaged-1}
   will be $C^1$-smooth  in  $I$. However, this will not   improve
   the fact that the  dispersion term vanishes with $I_k$ (and may have there a square-root singularity, see  \cite[Proposition 6.2]{GHKu} for an
   example). Thus the  well-posedness of equation  \eqref{averaged-1} is a delicate question.
\end{remark}

\begin{remark}\label{r_66}
Dispersion matrix $K$ of equation  \eqref{averaged-1} should  not   necessarily  be neither  symmetric nor  square, and
 if we replace it by another (possibly non-square)
 matrix $\bar{{K}}$ such that   $\bar{{K}}\bar{{K}}^t={S}$,  we would get a new  equation with  the same set of weak solutions. See \cite[Section~5.3]{SV} and \cite[Section~5.4.B]{brownianbook}.  This fact concerning equation \eqref{averaged-1} 	and other SDEs
 is systematically used below.
\end{remark}

The following analogue of Theorem~\ref{t_main02} holds  for solutions of equation~\eqref{v-equation1}.

\begin{theorem}\label{th_ham_main}Under  Assumption \ref{assumption-v-3}, for any $v_{0}\in\mathbb{R}^{2n}$
 the collection of  laws of the processes $\{I(v^{\eps}(\tau;v_{0})),\tau\in[0,T]\}$, $0<\eps \le1$,
  is tight in $\mathcal{P}(C([0,T],\mathbb{R}_+^n))$.  If
 we take  any sequence ${\eps}_j\to0$ such that
$\cD(I(v^{{\eps}_j}(\cdot;v_{0}))\strela Q^0\in\mathcal{P}(C([0,T],\mathbb{R}_+^n))$, then $Q^0$ is the law of a weak solution $I^0(\tau)$, $\tau\in[0,T]$, of the averaged equation \eqref{averaged-1}, equal $I_0=I(v_0)$ at $\tau=0$. Moreover,
\begin{equation}\label{post-bound}
\mathbf{E}\sup_{\tau\in[0,T]}|I^0(\tau)|^{q_0}\leqslant C_{q_0}(\|v_0\|,T),
\end{equation}
and for $k=1,\dots,n$,
\begin{equation}\label{small-delta}
\mathbf{E}\int_0^T\mathbf{1}_{\{I_k\in[0,\delta]\}}(I^0(\tau))d\tau\to0\text{ as }\delta\to0.
\end{equation}
\end{theorem}

As we  discussed in Remark \ref{remark-degenerate},
the uniqueness of a solution to the averaged equation   \eqref{averaged-1}, \eqref{initial-a1} is a delicate  issue.
 Therefore in the above theorem we state
the convergence only for a subsequence ${\eps}_j\to0$. In Section \ref{construct-effective} we construct an
{\it effective equation} in  $v$-variables such that for its solution $v(\cdot)$ with the same initial data $v_0$
the process of actions  $I(v(\cdot))$ is exactly the weak solution $I^0(\tau)$ from the theorem
above. This equation has locally Lipschitz coefficients,
 so its solution is unique, and thus the    convergence in Theorem~\ref{th_ham_main}  holds as $\eps\to0$.

Let us denote
$$
I^{\eps}(\tau)=I(v^{{\eps}}(\tau;v_{0}))\;\text{ and }\;\varphi^{\eps}(\tau)=\varphi(v^{{\eps}}(\tau;v_{0})), \;\tau\in[0,T].
$$
Denote also $\mathcal{S}_\delta=\{I\in\mathbb R^n_+\,:\, \min\limits_{1\leq j\leq n}I_j\leqslant\delta\}$ and $\bar A_{N,\delta}=A_{N}^\delta\setminus\mathcal{S}_\delta$, where $A_{N}^\delta$ is defined as in  \eqref{def_AN}, but
 with $P^I$ replaced by the drift term in \eqref{I-equation1}.
The following lemma  is a key technical result for a proof of Theorem~\ref{th_ham_main}. It states
  that with high probability the process $I^{\eps}(\tau)$, $\tau\in[0,T]$ stay away
from the locus  $\mathcal{S}_0=\cup_{k=1}^n\{I\in\mathbb{R}_+^n:I_k=0\}$, for  most of the time, uniformly in $\eps$.

\begin{lemma}\label{l_key_ham}
Under Assumption \ref{assumption-v-3}, there exist  a function $\kappa(\delta)$ and 	a function
$\alpha(\delta,N): \R_+\times\R_+\mapsto\R_+$
such that $\kappa(\delta)\to0$ as $\delta\to0$, $\alpha(\delta,N)\to0$ as $N\to\infty$  for each $\delta>0$, and
  \be\label{estim_ogo2} %(6.18)
  \;\mathbf{E}\int_0^T\mathbf{1}_{\bar A_{N,\delta}}(I^{\eps}(\tau))d\tau\leqslant\alpha(\delta,N),
\end{equation}
\begin{equation}\label{estim_ogo1} %(6.17)
  \mathbf{E}\int_0^T \mathbf{1}_{\mathcal{S}_\delta}(I^\eps(\tau))\,d\tau\leq \kappa(\delta)
  \ee
for any $\delta>0$,
uniformly in  $0<\eps\le1$.
\end{lemma}

A proof of   \eqref{estim_ogo2} follows from the same argument as in the demonstration of Lemma~\ref{prop_nondege} since
 the diffusion in \eqref{I-equation1} is non-degenerate outside $\cS_\delta$. That of   \eqref{estim_ogo1}
  is rather technical due to the  degeneracy at the locus
  $\cS_0$ and  presence  of the  ${\eps}^{-1}$-term in \eqref{v-equation1}.
The proof of this inequality  is  based on an argument, similar to that used in \cite{KP,K} for an  infinite-dimensional
stochastic equation.
 It relies on introducing a family of auxiliary processes which are suitable rotations of a solution $v^\eps(\tau; v_0)$.
 They  are constructed as It\^o processes  such that
 their actions coincide with $I^\eps$'s, while  presence of the rotations allows to remove terms of order  ${\eps}^{-1}$
from the  equations for these processes.
Detail of the proof of \eqref{estim_ogo1} is provided in Appendix \ref{app_proof_lem_not0}.

With the help of Lemma~\ref{l_key_ham}  the proof of Theorem~\ref{t_main02} can be adapted to demonstrate
 Theorem~\ref{th_ham_main}. We only give here a sketch of the corresponding argument,
emphasizing the differences and leaving  detail to the reader.
\medskip

{\it Sketch of the proof of Theorem~\ref{th_ham_main}}:
We introduce the starting time $\tau_0$,  number
$N$ and  intervals $\Delta_j$ in the same way as at Step 1   in the proof of Lemma~\ref{lemma-key-1}.
In particular, $\Delta_j=[\tau_j,\tau_{j+1})=[\tau_0+j\eps N, \tau_0+(j+1)\eps N)$.
By the argument from the proof of Lemma~\ref{l_compa} and Chebyshev's inequality, for any $\delta>0$ and $N>0$ it  holds that
\begin{equation}\label{lem_incre}
\mathbf{P}\big\{\ \max\limits_{0\leqslant\tau'<\tau''\leqslant\tau'+\eps N\leqslant T}\ |I^\eps(\tau')-I^\eps(\tau'')|\geqslant\tfrac\delta2\big\}\longrightarrow 0,\quad\hbox{\rm as }\eps\to0.
\end{equation}
By  Prokhorov's  theorem this relation implies the tightness of the family $\{I^\eps(\cdot)\}$ in the space
$C([0,T]; \mathbb R^n)$. So any sequence $\eps_j\to 0$ has a subsequence such that along it
the laws $\mathcal{L}(I^\eps(\cdot))$ converge weakly in $C([0,T]; \mathbb R^n)$  to a limit probability
measure $Q^0$.
From Lemma~\ref{l_key_ham} and convergence  \eqref{lem_incre} we derive that  for any $\delta>0$ and $N>0$
$$
\lim\limits_{\eps\to0}\ \mathbf{E}\,\frac{\#\big\{j\,:\, \{I^\eps(\tau)\,:\,
\tau\in\Delta_j\big\}\cap \mathcal{S}_{\frac\delta2}\not=\emptyset\big\}}{ {T}/({\eps N})}
=0.
$$
For this proof  we  call an interval $\Delta_j$ {\it typical}, if in addition to the properties, listed in \eqref{typ},
the curve $\{I^\eps(\tau)\,:\,\tau\in\Delta_j\}$  does not intersect  the set
$ \mathcal{S}_{\frac\delta2}$.
  Then direct analogies of
 Lemmas~\ref{lemma-key-1} and \ref{lemma-key-2} hold for the process
 $(I^{\eps}(\tau),\varphi^{\eps}(\tau))$ due to the  argument, used in Section~\ref{s_co_ac} and
 enriched by Lemma~\ref{l_key_ham}. Next,  arguing  in the same way as in Section \ref{s_lp}, we show that the limiting
  measure $Q^0$ is a solution of the martingale problem for equation \eqref{averaged-1},
  satisfying $Q^0\{I(0)=I_0\}=1$.  Relation \eqref{small-delta} follows from \eqref{estim_ogo2}.

\section{The effective equation}\label{construct-effective}

In this section we construct an
 effective equation for  \eqref{v-equation1} with  small $\eps$.
 This is a $v$-equation such that under the mapping $v\mapsto I(v)$ its weak
solutions go to  weak solutions of \eqref{averaged-1}. Thus by Theorem~\ref{th_ham_main}
the equation   controls the behaviour of actions $I^k(v^\eps(\tau;v_0))$ as $\ep\to0$.
The construction of the effective  equation is a finite-dimensional
modification of the infinite-dimensional construction, used in \cite{K} for purposes of averaging a stochastic PDE
with analytic nonlinearity.

To get the effective equation we firstly remove from equation \eqref{v-equation1} the fast rotating terms
$
\ep^{-1}W_k \bv_k^{\bot},
$
and then average the resulting equation with respect to the action of the $n$-torus on $\R^{2n}$, using the rules of the stochastic averaging,
similar to how we earlier got the averaged $I$-equation \eqref{averaged-eq-1} from the $I$-equation in \eqref{ori_scaled}. The action of
$\T^n =\{\theta\}$ on $\R^{2n}$ is given by the block-diagonal matrix
\be \label{Phi_th}
 \Phi_\theta=
 \text{diag}\{\Phi_\theta^1,\dots,\Phi_\theta^n\}, \quad
 \Phi_\theta^k= \left(\begin{array}{cc}\cos\theta_k&-\sin\theta_k\\
\sin\theta_k&\cos\theta_k\end{array}\right),\quad 1\leqslant k \leqslant n.
\ee

The drift in the effective equation is the  averaging
$
\lan P\ran = (\lan \bP\ran_1, \dots, \lan\bP\ran_n)
$
 of the vector field $P$  in \eqref{v-equation1}
 with respect to the action $\Phi_\theta$. For  convenience of  future calculation  we abbreviate
$
\lan P\ran(v) =  R(v) = ( \mathbf{R}_1(v), \dots  \mathbf{R}_n)(v).
$
Then
$$
\bR_k(v):=
\lan \bP\ran_k(v)  =\int_{\mathbb{T}^n}\Phi^k_{-\theta}{\mathbf{P}}_k(\Phi_\theta v)d\theta,\quad k=1,\dots,n
$$
(cf. \cite[Section 3]{HGK}).

To obtain the dispersion matrix $\llangle {B}\rrangle(v)$ of the effective equation we start with
   $2n\times 2n$ matrix ${X}(v)=\big(X_{km}(v)\big)_{1\leqslant k,m\leqslant n}$,
   formed by  $2\times 2$-blocks $X_{km}(v)$,
$$
X_{km}(v): =\sum_{j=1}^{n_1}\int_{\mathbb{T}^n}\Phi^k_{-\theta}B_{kj}(\Phi_\theta v)
\big(B_{mj}(\Phi_\theta v)\big)^t\Phi^m_{\theta}d\theta.
$$
That is,
$\ {X}(v)=\int_{\mathbb{T}^n}\Phi_{-\theta}{B}(\Phi_\theta v)( {B}(\Phi_\theta v))^t\Phi_{\theta}d\theta$, where we denoted by $B$  the
block matrix ${B}=\big(B_{kj}\big)$.
The dispersion matrix in question   $\llangle {B}\rrangle(v)=\big(\llangle B_{kj}\rrangle(v)\big)_{1\leqslant k,j\leqslant n}$ is
  the principal  square root of  ${X}(v)$ (see the 10-th footnote). So
\be\label{formula}
\sum_{j=1}^n\llangle B_{kj}\rrangle(v) \llangle B_{mj}\rrangle^t(v)=
\sum_{j=1}^{n_1}\int_{\mathbb{T}^n}\Phi^k_{-\theta}B_{kj}(\Phi_\theta v)
\big(B_{mj}(\Phi_\theta v)\big)^t\Phi^m_{\theta}d\theta
\ee
for $1\leqslant k, m\leqslant n.$

Then the  effective equation for \eqref{v-equation1} is the following one:
\begin{equation}\label{effective1}
d\bv_k(\tau)={\mathbf{R}}_k(v)d\tau+\sum_{j=1}^n \llangle B_{kj}\rrangle(v) d\bbeta_j(\tau),\quad k=1,\dots,n,
\end{equation}
or
$$
d v(\tau) = R(v) d\tau + \llangle B \rrangle(v) d\beta(\tau),
$$
where $\beta(\tau) = (\bbeta_1,\dots,\bbeta_n)(\tau)$ is a standard Wiener process in $\R^{2n}$, defined on the space \eqref{OmOm}.

\begin{proposition}\label{square-root} Under (2) and (3) of Assumption \ref{assumption-v-3},

i) the vector-function $R(v)$ and
 matrix functions ${X}(v)$ and  $\llangle {B}\rrangle(v)$ are  locally Lipschitz in $v$;

ii) for any $\tilde\theta\in\mathbb{T}^n$, $R(\Phi_{-\tilde\theta} v) = \Phi_{-\tilde\theta} R(v)$, while
 ${X}(\Phi_{-\tilde\theta} v)=\Phi_{-\tilde\theta} {X}(v)\Phi_{\tilde\theta}$ and
 $\llangle {B}\rrangle(\Phi_{-\tilde\theta} v)=\Phi_{-\tilde\theta} \llangle {B}\rrangle (v)\Phi_{\tilde\theta}$.
\end{proposition}

\begin{proof} i) The local
 Lipschitz continuity of $R(v)$ and  ${X}(v)$ follow  from  the relations which define them. Since the operator $BB^t(v)$
is uniformly elliptic (see \eqref{ellipt}), then $X(v)$ is uniformly elliptic as well, so the Lipschitz continuity of
$\llangle {B}\rrangle(v) = (X(v))^{1/2}$ is a consequence of  \cite[Lemma~5.2.1 and Theorem~5.2.2]{SV} and the Lipschitz continuity of
${X}(v)$.

 ii) The  relations for $R(v)$ and $X(v)$ are direct consequences of their definitions, and the relation for  $\llangle {B}\rrangle$
 follows from that for ${X}$.
\end{proof}

Since the coefficients  of equation  \eqref{effective1}  are  locally Lipschitz, then its
  strong solution,  if exists, is unique. So by
 the Yamada-Watanbe theorem (see in \cite{brownianbook})
we  have
\begin{proposition}\label{effective-unique-weak}
If $v^1(\tau)$ and $v^2(\tau)$, $0\le\tau\le T$,
are weak solutions of equation \eqref{effective1} such that $\cD(v^1(0))=\cD(v^2(0))$,
then $\cD(v^1(\cdot))=\cD(v^2(\cdot))$,  and a strong solution  with the initial data $v^1(0)$
 exists for $0\le  \tau \le T$.
\end{proposition}

%Due to ii) of Proposition \ref{square-root} (and the assertion, mentioned in Remark \ref{r_66})
 %the set of weak solutions is invariant under the rotations  $\Phi_\theta$, $ \theta\in\mathbb{T}^n$:

%\begin{proposition}\label{phi-invariant}
%If $v(\tau)$ is a weak solution of equation  \eqref{effective1} then for any $\theta\in\mathbb{T}^n$,
 %$\Phi_\theta\big( v(\tau)\big)$  also is its weak solution.
%\end{proposition}

Since for a matrix $Q$ we have
$
\| Q\|^2_{HS} = $tr$\,QQ^t,
$
then  taking the trace of relation  \eqref{formula} with $k=m$ we get that
$\
\sum_{j=1}^n \| \llan B_{kj}\rran (v)\|^2_{HS} = \sum_{j=1}^{n_1} \int_{\T^n} \| \Phi^k_\theta B_{kj}(\Phi_{\theta}v)\|^2_{HS}.
$   %}
Using this equality we write  It\^o's formula for the actions
 $I_k(v(\tau))$, $k=1, \dots, n$, of a solution  $v(\tau)$ for
\eqref{effective1}, as
\begin{equation}\label{I-effective}
dI_k=\bv_k^t{\mathbf{R}}_k(v)d\tau+ \frac{1}{2} \Big(\sum_{j=1}^n\int_{\mathbb{T}^n}\|\Phi^k_\theta B_{kj}(\Phi_\theta v)\|_{HS}^2d\theta
\Big)d\tau +\sum_{j=1}^n\bv_k^t\llangle B_{kj}\rrangle(v)d\bbeta_j(\tau).
\end{equation}
The first term of  the drift in this equation may be re-written as
\[\bv_k^t{\mathbf{R}}_k(v)=\bv_k^t\int_{\mathbb{T}^n}\Phi_{-\theta}^k{\mathbf{P}}_k(\Phi_\theta v)d\theta=\int_{\mathbb{T}^n}(\bv_k^t{\mathbf{P}}_k)(\Phi_\theta v)d\theta=\langle\bv_k^t{\mathbf{P}}_k\rangle(I).
\]
Since $\|\Phi^k_\theta B_{kj}(\Phi_\theta v)\|_{HS}^2=\|B_{kj}(\Phi_\theta v)\|_{HS}^2,$ then in the second term of the drift we have
 $\int_{\mathbb{T}^n}\|\Phi^k_\theta B_{kj}(\Phi_\theta v)\|_{HS}^2d\theta=\langle \|B_{kj}\|_{HS}^2\rangle(I).$
 Therefore  the drift in \eqref{I-effective} is $F_k(I)$, i.e. is
 the same as that in \eqref{averaged-1}.

Using once again \eqref{formula} we see that the diffusion matrix in \eqref{I-effective}
 is $\bar{{S}}=(\bar S_{km})$ with
\[ \begin{split}
\bar S_{km}&:=\sum_j\mathbf{v}_k^t\llangle B_{kj}\rrangle(v)(\llangle B_{mj}\rrangle(v))^t \mathbf{v}_m \\
%&=\sum_{j}\int_{\mathbb{T}^n}\bv_k^t\Phi_{-\theta}^kB_{kj}(\Phi_\theta v)(B_{mj}(\Phi_\theta v))^t\Phi_\theta^m\bv_md\theta\\
&=\sum_{j}\int_{\mathbb{T}^n}(\Phi_\theta^{k}\bv_k)^tB_{kj}(\Phi_\theta v)(B_{mj}(\Phi_\theta v))^t\Phi_\theta^m\bv_md\theta\\
&=\big\langle\sum_j \bv_k^tB_{kj}(v)B_{mj}^t(v)\bv_m\big\rangle(I)=S_{km}(I),
\end{split}
\]
where $S_{km}(I)$ is as in \eqref{diffusion-a1}.  We conclude that the diffusion matrices of equations
 \eqref{I-effective} and  \eqref{averaged-1} also coincide.  Hence,

\begin{proposition} \label{p_7.4} Let the process $v(\tau)\in\mathbb{R}^{2n}$, $0\leqslant \tau\leqslant T$ be a weak solution of \eqref{effective1}.

i)
 Then the process $\big(I(v(\tau)), v(\tau)\big)\in\mathbb{R}_+^n\times\mathbb{R}^{2n}$, $0\leqslant \tau\leqslant T$ is a strong solution of      system  \eqref{I-effective}$+$\eqref{effective1}, driven by the  set of Brownian motions, corresponding to the weak solution $v$.

 ii) The drift and diffusion matrix in \eqref{I-effective} are functions of the actions $\{I_k\}$ and are the same as in equation \eqref{averaged-1}. So
  the process $I(v(\tau))$, $0\leqslant\tau\leqslant T$ is a weak solution of \eqref{averaged-1}.
\end{proposition}

A disadvantage of  system  \eqref{I-effective}$+$\eqref{effective1} is that dispersion matrix in  \eqref{I-effective} depends both on the actions $I$ and angles $\f$,
despite the corresponding diffusion matrix depends only on the actions. To fix this let us denote by
 $(\sqrt{2I(v)},0)$ the $2n$-vector with components
 $(\sqrt{2I_i(v)},0)^t$, $i=1,\dots,n$. Then
 $(\sqrt{2I(v)},0) = \Phi_{-\varphi(v)}v $, so by Proposition~\ref{square-root}.ii),
 \[\begin{split}
 \bv_k^t\llangle B_{kj}\rrangle (v)d\bbeta_j(\tau)&=(\sqrt{2I_k(v)},0)\Phi^k_{-\varphi(v)}\llangle B_{kj}\rrangle(\Phi_{\varphi(v)}(\sqrt{2I(v)},0))d\bbeta_j(\tau)\\
&=(\sqrt{I_k(v)},0)\llangle B_{kj}\rrangle((\sqrt{2I(v)},0))\Phi^j_{-\varphi(v)}d\bbeta_j(\tau)\\
&=M_{kj}(I)d\tilde{\bbeta}_j(\tau).
\end{split}\]
Here $d\tilde{\bbeta}_j(\tau)=\Phi^j_{-\varphi(v)}d\bbeta_j(\tau)$, $j=1,\dots,n$, are differentials of
 independent  standard Brownian motions in
 $\mathbb{R}^2$ and $M_{kj}(I)$ is the 2-vector
$(\sqrt{I_k},0)\llangle B_{kj}\rrangle\big((\sqrt{I},0)\big)$.
 Then equation  \eqref{I-effective} can be re-written as
\begin{equation}\label{I-E2}dI_k=F_k(I)d\tau+\sum_{j=1}^n\big( M_{kj}(I),d\tilde{\bbeta}_j(\tau)\big),\; k=1,\dots,n.
\end{equation}
When driven by the  set of Brownian motions $\{\tilde\bbeta_j(\tau)\}$, the effective equation \eqref{effective1} reads
\begin{equation}\label{effective2}
d\bv_k ={\mathbf{R}}_k(v)d\tau+\sum_{j=1}^n \tilde B_{kj}(v)d\tilde{\bbeta}_j(\tau),\; k=1,\dots, n,
\end{equation}
where $\tilde B_{kj}(v)=\llangle B_{kj}\rrangle (v)\Phi^j_{\varphi(v)}.$
The system   \eqref{I-E2}+\eqref{effective2} is just the system  \eqref{I-effective}+ \eqref{effective1}, written using another
standard Wiener process in $\R^{2n}$, so the two systems have the same sets of weak solutions.
 In difference with equation \eqref{effective1},
dispersion matrix $\tilde B(v)$ in \eqref{effective2} is  not locally Lipschitz. But for any $N\ge1$ and $\delta>0$ it is Lipschitz in domain
$\{ v: \|v\| \le N, \| \bv_j\| >\delta \ \forall\,j\}$.

 The  drift and diffusion in equation \eqref{I-E2}  are the same as in \eqref{I-effective}, so by Proposition~\ref{p_7.4}\,ii) they are the same as in
 equation \eqref{averaged-1}. Thus equations \eqref{I-E2} and  \eqref{averaged-1} have the same set of weak solutions. We have established

\begin{lemma}\label{same-weak2}
Systems   \eqref{I-E2}+\eqref{effective2}  and \eqref{I-effective}+ \eqref{effective1}
have the same set of weak solutions. So do equations
 \eqref{I-E2} and \eqref{averaged-1}.
\end{lemma}

\section{Lifting of  solutions}\label{lifting}
In this section we prove that a weak solution $I(\tau)$
 of the averaged equation \eqref{averaged-1}, constructed in Theorem~\ref{th_ham_main}, is distributed as $I(v(\tau))$,
 where $v(\tau)$ is   some weak solution of the effective equation \eqref{effective1}. That is,  $I(\tau)$  can be
   lifted to a weak solution of \eqref{effective1}.  We follow a strategy from \cite{K}, where such a lifting is constructed
    for an infinite dimensional equation.
Since we work with  weak solutions of the equations, then using  Lemma~\ref{same-weak2} we replace averaged
equation \eqref{averaged-1}
by equation \eqref{I-E2}, and effective equation \eqref{effective1} --  by equation \eqref{effective2}. If
 $v(\tau)$ solves \eqref{effective2}, then by It\^o's formula the equation for the vector of actions  is
  \eqref{I-E2}, while  equations for $\varphi_k$'s read
\begin{equation}\label{phi-E2}
d\varphi_k(\tau)=R_k^\varphi(v)d\tau+\sum_{j=1}^n\frak R_{kj}^\varphi(v)d\tilde\bbeta_j(\tau),\quad k=1,\dots,n,
\end{equation}
where
\[R_k^\varphi(v)=\!\big(\nabla_{\bv_k}\!\arctan\big(\frac{v_k}{v_{-k}}\big)\big)\cdot {\mathbf{R}}_k(v)+
\frac{1}{2}\sum_{j=1}^n\text{Trace}\Big(\tilde B_{kj}(v)\Big(\nabla_{\bv_k}^2\!\arctan\big(\frac{v_k}{v_{-k}}\big)\Big)\big(\tilde B_{kj}(v)\big)^t\Big),\]
and
$\
\frak R_{kj}^\varphi(v)=\Big(\nabla_{\bv_k}\arctan\big(\frac{v_k}{v_{-k}}\big)\Big)\cdot \tilde B_{kj}(v). $

For any $\theta=(\theta_1,\dots,\theta_n)\in \mathbb{T}^n$ and any vector $I=(I_1,\dots,I_n)\in \mathbb{R}^n_+$ we set
\be\label{V_theta}
V_\theta(I)=(\bV_{\theta_1}(I_1),\dots,\bV_{\theta_n}(I_n))),
\ee
 where
 $$
 \bV_{\alpha}(J)=(\sqrt{2J}\cos\alpha,\sqrt{2J}\sin\alpha)^t\in\mathbb{R}^2.
 $$
Then $\varphi(V_\theta(I))=\theta$, $I(V_\theta(I))=I$ and  the mapping $\mathbb{R}_{+}^n\times\mathbb{T}^n: (I,\varphi)\mapsto V_\varphi(I)$
 is a left-inverse for the mapping   $v\mapsto(I,\varphi)$.
   For any $I=(I_1,\dots,I_n)\in\mathbb{R}_+^n$ let us denote
\be\label{[I]}
[I]=\min_{1\leqslant k\leqslant n}\{I_k\}.
\ee
Then for a  $\delta>0$
 the mapping $v\mapsto (I,\varphi)$ defines a diffeomorphism  of  domain
 $\mathbb R^{2n}_\delta:=\{v\in\mathbb{R}^{2n}: [I(v)]>\delta\}$ and domain
$\{I\in \mathbb{R}_+^n\,:\,[I]>\delta\}\times \mathbb{T}^n$.
Therefore we have

\begin{lemma} \label{v-I-phi} For any $\delta>0$, on domain $\mathbb{R}_\delta^{2n}$
 equation  \eqref{effective2} is equivalent to system
  \eqref{I-E2}+\eqref{phi-E2} in the following sense:
  let $\tau_1\leqslant \tau_2$ be two stopping times with respect to the natural filtration. Then

i) if for $\tau_1 \le \tau \le \tau_2$ a process $v(\tau)$ lies in  $ \mathbb{R}_\delta^{2n}$ and is a weak
solution of \eqref{effective2}, then for such $\tau$'s
 $(I(v),\varphi(v))$  is a weak solution of \eqref{I-E2}+\eqref{phi-E2};

ii)   if for $\tau_1 \le \tau \le \tau_2$ a process  $(I(\tau),\varphi(\tau))$ satisfies $[ I(\tau)]>\delta$ and is a
 weak solution of  \eqref{I-E2}+\eqref{phi-E2}, where $v(\tau)=V_{\varphi(\tau)}(I(\tau))$,
 then  for such $\tau$'s
 the process $v(\tau)=V_{\varphi(\tau)}(I(\tau))$ is a weak solution of \eqref{effective2}.
\end{lemma}

The following statement is the main result of this section.

\begin{theorem}\label{lifting-thm} If the process $I^0(\tau)=(I_k^0(\tau), k=1,\dots,n)\in\mathbb{R}^n_+$, $0\leqslant \tau\leqslant T$,
 is a weak solution of equation \eqref{averaged-1}, constructed in Theorem~\ref{th_ham_main} by taking  the limit along
 a sequence
 ${\eps}_i\to0$, then, for any vector $\theta\in \mathbb{T}^n$, there exists a
   weak solution $v^\theta(\tau)\in\mathbb{R}^{2n}$,
 $0\leqslant \tau\leqslant T$ of the  effective equation  \eqref{effective1} such that

(i) the law of $I(v^\theta(\cdot))$ coincides with that of $I^0(\cdot)$;

(ii) $v^\theta(0)=V_\theta(I_0)$ a.s.
\end{theorem}

  The properties of solutions of equation \eqref{effective1} are  important for  analysis in the following sections.
  They are  the subject of Theorem~\ref{averaged-unique} below.

We  begin with explaining the key ideas of the theorem's proof.

{\it Strategy of the proof:} By Lemma \ref{same-weak2} we may regard $I^0(\tau)$ as a weak solution of equation \eqref{I-E2}.
Since the uniqueness of a solution $v(\tau)$ is claimed in Proposition~\ref{effective-unique-weak},
then only its existence and properties (i), (ii) should be established.
For any $\delta>0$, we will divide $[0,T)$ into a finite or
countable  set of random closed intervals
 (see Figure \ref{pic1})
 $\Lambda_j$, $j\geqslant0$, and $\Delta_j$, $j\geqslant1$, such that

(1) $\Lambda_0\leqslant \Delta_1\leqslant \Lambda_1\leqslant \Delta_2\leqslant\dots$,

(2) $[I^0]\geqslant\delta$ on each $\Lambda_j$, and $[I^0]\leqslant 2\delta$ on each $\Delta_j$.\\
For definiteness we assume that $[I_0]\ge\delta$.

Next we iteratively construct  on these intervals a process $v^\delta(\tau)$  such that
\begin{equation}\label{equal-i-v}
I(v^\delta(\tau)) \equiv I^0(\tau), \; \; a.e.
\end{equation}
Suppose that we  already know $v^\delta$
 at the left end point of some $\Lambda_j$. To construct $v^\delta(\tau)$ on $\Lambda_j$ we note that since on every $\Lambda_j$ we have
 $[I(v^\delta(\cdot))]=[I^0(\cdot)]\geqslant \delta$ , then by Lemma~\ref{v-I-phi} there
  equation \eqref{effective2} is equivalent to \eqref{I-E2}+\eqref{phi-E2}.
   As $I(\tau)=I^0(\tau)$ is known, it remains to solve \eqref{phi-E2}, regarded as a stochastic equation with progressively measurable
   coefficients for  $n$-vector~$\varphi(\tau)$. Since the initial value of $\varphi$ is given on the left end point of $\Lambda_j$,
   such a solution $\varphi$ is uniquely determined. Then for $\tau\in \Lambda_j$ we set  $v^\delta(\tau):=V_{\varphi(\tau)}(I^0(\tau))$.
    By Lemma~\ref{v-I-phi} this $v^\delta$  solves \eqref{effective2} weakly on $\Lambda_j$. Clearly, \eqref{equal-i-v} is satisfied.

On the next  interval $\Delta_{j+1}$ we have $[I^0]\leqslant 2\delta$. We want to extend $v^\delta(\tau)$ to $\Delta_{j+1}$,
 keeping the property \eqref{equal-i-v},  so that when eventually
  $v^\delta$ is constructed on all $\Lambda_j$'s and $\Delta_k$'s,
  we may obtain the desired weak solution of \eqref{effective2} by taking a limit as $\delta\to0$. Such a task turns out to be not easy
     since by \eqref{equal-i-v} we have  $\bv_k=\big(\sqrt{2I_k(\tau)}\cos\varphi_k(\tau),\sqrt{2I_k(\tau)}\sin\varphi_k(\tau)\big)$ with
     some phase function  $\varphi_k$, and so on $\Delta_{j+1}$ a-priori $|\dot{\bv}_k|\sim\delta^{-1/2}$.
     Hence, a naive extension may fail to guarantee the existence of a limit as $\delta\to0$.
     To construct a right lifting of $I^0(\tau)$ when it is small, we use the fact the $I^0(\tau)$ is a limit of the process $I^{{\eps}_i}(\tau):=I(v^{{\eps}_i}(\tau))$, where $v^{{\eps}_i}$ solves equation \eqref{v-equation1} with ${\eps}={\eps}_i$. The process $v^{{\eps}_i}(\tau)$ is a lifting of $I^{{\eps}_i}(\tau)$ which is singular as ${\eps}_i\to0$. In Appendix~\ref{app_proof_lem_not0}
      for each ${\eps}_i>0$ a modified process $\bar v^{{\eps}_i}(\tau)$ is constructed such that $I(\bar v^{{\eps}_i}(\tau))=I(v^{{\eps}_i}(\tau))$ and $|\frac{d}{d\tau}\bar v^{{\eps}_i}|\sim1$ as ${\eps}_i\to0$. A limit in law of processes
      $\bar v^{{\eps}_i}(\cdot)$ as ${\eps}_i\to0$ provides  a right lifting of $I^0$ on  interval $\Delta_{j+1}$. We thus extend
       the process $v^\delta(\tau)$ to $\Delta_{j+1}$ in such a way that \eqref{equal-i-v} holds  and $|\dot{v}^\delta|\sim1$.

By iterating the two  constructions above we obtain a process $v^\delta(\tau)$
 which solves \eqref{effective2} for $\tau\in \cup\Lambda_j$ and satisfies  good estimates on the complementary set $\cup \Delta_j$.
 By Theorem~\ref{th_ham_main} the Lebesgue measure of $\cup \Delta_j$ becomes small with $\delta$. This allows us to
 show that any limit distribution of the process $v^\delta(\cdot)$ as $\delta\to0$ gives a weak solution $v(\cdot)$ for~\eqref{effective2}.

To run  this construction we need  good upper bounds for  the numbers of  intervals $\Lambda_j$ and $\Delta_j$,
which could be large when the norm of the process $I^0(\tau)$ is large. To get the bounds
 we begin the proof by introducing the stoping times $\tau_N=\inf_\tau\{|I^0(\tau)|=N/2\}$ and replacing $I^0(\tau)$ by
 a trivial modification for $\tau\geqslant \tau_N$. Then we first obtain a weak solution $v_N(\tau)$, corresponding
 to the modified process $I_N^0(\tau)$,  and next take the limit as $N\to\infty$ to get a real weak solution $v(\tau)$ of \eqref{effective2}.

\begin{proof} Let us introduce a natural filtered measurable space $\{\Omega, \cB,\{\cB_\tau\},0\leqslant\tau\leqslant T\}$
 for the problem we consider, where $\Omega$ is the Banach space
\be\label{Om}
\Omega=\Omega_I\times \Omega_v:=C([0,T],\mathbb{R}_+^n)\times C([0,T],\mathbb{R}^{2n}),
\ee
$\cB$ is its Borel $\sigma$-algebra and $\cB_\tau$ is the $\sigma$-algebra generated by the set of random variables
$\{r(s): 0\leqslant s\leqslant \tau \text{ and } r(\cdot)\in \Omega\}$.
Denote by
 $\pi_I$ and $\pi_v$ the natural projections $\pi_I: \Omega\to\Omega_I$ and $\pi_v:\Omega \to\Omega_v$.
We will prove the theorem by  constructing  a probability measure
 $\bQ$ on $\Omega$ such that $\pi_I\circ \bQ=\cD(I^0(\cdot))$,  $\pi_v\circ \bQ$ is the distribution of a weak solution of \eqref{effective2},
 and  $\bQ$-a.s. for $(I'(\cdot),v'(\cdot))\in \Omega$ we have
  $I(v')=I'$. This will be achieved  in four  steps.

{\bf Step 1}. Redefine the equations for large amplitudes.

For any $N\in\mathbb{N}$ consider the stopping time
\[\tau_N=\inf\{\tau\in [0,T]\Big|\|v(\tau)\|^2=2|I(v(\tau))|=N\},\]
(here and in similar situations below $\tau_N=T$ if the set is empty). For $\tau\geqslant \tau_N$ and each $\ep>0$ we redefine equation \eqref{v-equation1} to the trivial equation
\begin{equation}\label{v-trivial}
d\bv_k=d\bbeta_k(\tau), \; k=1,\dots,n,
\end{equation}
and redefine accordingly equations \eqref{I-equation1}, \eqref{I-E2}  \eqref{effective2} and \eqref{phi-E2}.
We  denote the new equations as \eqref{v-equation1}$_N$, \eqref{I-equation1}$_N$ \eqref{I-E2}$_N$, \eqref{effective2}$_N$
and \eqref{phi-E2}$_N$. So if $v_N^\ep(\tau)$ is a solution of \eqref{v-equation1}$_N$, then $I^\ep_N(\tau):=I(v^\ep_N(\tau))$
satisfies \eqref{I-equation1}$_N$.  That is, for $\tau\leqslant \tau_N$,  it satisfies \eqref{I-equation1}, while for $\tau\geqslant \tau_N$ it is a solution of the equation
\[
dI_k=\tfrac{1}{2}d\tau+(v_kd\beta_k+v_{-k}d\beta_{-k})=\tfrac{1}{2}d\tau+\sqrt{2I_k}\,dw_k(\tau),\quad k=1,\dots,n,
\]
where $w_k(\tau)$ is the Wiener process $\int_0^\tau(\cos\varphi_kd\beta_k(\tau)+\sin\varphi_kd\beta_{-k}(\tau))$. So \eqref{I-equation1}$_N$ is the equation
\[dI_k=\mathbf{1}_{\tau\leqslant \tau_N}\cdot\langle \text{r.h.s of \eqref{I-equation1}}\rangle+
\mathbf{1}_{\tau\geqslant \tau_N}(\tfrac{1}{2}d\tau+\sqrt{2I_k}\,dw_k(\tau)),\quad k=1,\dots,n.\]
Accordingly, the averaged equation  \eqref{I-E2}$_N$ reads
\[
dI_k=\mathbf{1}_{\tau\leqslant \tau_N}\Big(F_k(I)d\tau+\sum_{j=1}^n\big( M_{kj}(I),d\tilde\bbeta_j(\tau)\big)\Big)+\mathbf{1}_{\tau\geqslant \tau_N}\Big(\tfrac{1}{2}d\tau+\sqrt{2I_k}d\tilde\beta_k(\tau)\Big), \;k=1,\dots,n.
\]
Here $\tilde\bbeta_j(\tau)$, %=\left(\begin{array}{c}\tilde\beta_j\\
%\tilde\beta_{-j}\end{array}\right)$,
 $j=1,\dots,n$, are independent standard Brownian motions
in $\mathbb{R}^2$.% and we have replaced $w_k(\tau)$ by $\tilde\beta_k$ which makes no
%difference  since we are interested in weak solutions of the system.

For the sequence $\ep_i\to0$, where we have the convergence
$\cD(I^{\ep_i}(\cdot))\rightharpoonup\cD(I^0(\cdot))$, choosing a suitable subsequence (if neccessary) we
 achieve that also  $\cD(I^{\ep_i}_N(\cdot))\rightharpoonup\cD(I_N(\cdot))$ for some process $I_N(\tau)$, for each  $N\in\mathbb{N}$.

\begin{lemma}\label{N-S}
For every  $N\in\mathbb{N}$,

i) the process $I_N(\tau)$, $0\leqslant\tau\leqslant T$, is a weak solution of \eqref{I-E2}$_N$ such that $\cD(I_N)=\cD(I^0)$ for $\tau\leqslant \tau_N$ and $\cD(I_N(\cdot))\rightharpoonup\cD(I^0(\cdot))$ as  $N\to\infty$.

ii)
 the statement in Lemma~\ref{v-I-phi} holds true with \eqref{effective2}$_N$ and \eqref{I-E2}$_N$$+$\eqref{phi-E2}$_N$.
\end{lemma}

\begin{proof}
i) The first part of the statement follows from the same argument as in the proof of Theorem~\ref{th_ham_main}.
Recalling that, for each $\ep>0$, $\cD(I^\ep_N(\cdot))=\cD(I^\ep(\cdot))=: \mathbf{Q}^\ep$ for $\tau\leqslant\tau_N$ and passing to the limit as $\ep_i\to0$ we get the second assertion of the lemma. As $\bQ^\ep\{\tau_N<T\}\leqslant CN^{-1}$ uniformly in $\ep$, then the last assertion also follows.

ii) The assertion  follows by the same argument as  that used to prove  Lemma~\ref{v-I-phi}.
\end{proof}

Now we fix any $N\in\mathbb{N}$. Our goal is to  construct for each $\delta>0$ a measure $\bQ_\delta^N$ on $\Omega$ such that for its natural process $(I^N_\delta(\tau),v^N_\delta(\tau))$, $\tau\in [0,T]$ we have
\begin{enumerate}
\item $\mathcal{D}(I^N_\delta(\cdot))=\mathcal{D}(I_N(\cdot))$;
\item  $I(v^N_\delta(\cdot))\equiv I^N_\delta(\cdot)$, $\bQ_\delta^N$-a.s;
\item  the process $v_N^\delta(\tau)$ is an It\^o process with bounded (in terms of $N$) drift term and diffusion matrix.
Moreover, it solves \eqref{effective2}$_N$ for $\tau$ outside a small  random set, where $[I(\tau)]\lesssim\delta$ (see
\eqref{[I]}).
\end{enumerate}
Next we  will prove the assertion of Theorem~\ref{lifting-thm}
by taking the limits $\delta\to0$ and  $N\to\infty$.
\medskip

{\bf Step 2: } Construction  of the  measure $\mathbf{Q}_\delta:=\bQ^N_\delta$ for every $\delta>0$ and $N$ fixed.

We start with finding an auxiliary
   It\^o process $\bar w(\tau)$ which covers a version of the process $I_N(\tau)$ (but has no relation
with the effective equation \eqref{effective2}). Since $N$ is fixed,  then below in the step the index $N$ usually is dropped.

\begin{lemma} \label{lm-ito-limit}There exists a continuous $\cB_\tau$-adapted continuous
process $(\bar I(\tau), \bar w(\tau))\in\mathbb{R}_+^n\times\mathbb{R}^{2n}$, $0\leqslant\tau\leqslant T$, such that $\cD(\bar I(\cdot))=\cD(I_N(\cdot))$, $\bar I(\cdot)=I(\bar w(\cdot))$ a.s.,
 and  $  \bar w(\tau) $ is an It\^o process in $\R^{2n}$ of the form
\begin{equation}\label{ito-limit}
\begin{split}
d \bar w(\tau)=B(\tau)d\tau+a(\tau)d\beta (\tau), \quad \tau\in [0,T],
\\
 |B(\tau)|\leqslant C,\quad C^{-1} \mathbb{I}\leqslant a(\tau)\big(a(\tau)\big)^t\leqslant C\mathbb{I},\;
\end{split}
\end{equation}
where the constant $C$ depends only on $N$, and $\beta(\tau)$ is a standard Brownian motion in $\mathbb{R}^{2n}$.
\end{lemma}

\begin{proof}
By the construction in Appendix \ref{app_proof_lem_not0} (see there  \eqref{deltaeq}),
for each $\epsi>0$ there exists an It\^o process $\bar w^\epsi(\tau)\in\mathbb{R}^n$, $0\leqslant \tau\leqslant T$,
 such that

$$%\label{ito-limit}
d \bar w^\epsi(\tau)=B^\epsi(\tau)d\tau+a^\epsi(\tau)d\beta (\tau), \quad \tau\in [0,T].
$$
Here
$$
 |B^\epsi(\tau)|\leqslant C,\quad C^{-1} \mathbb{I}\leqslant a^\epsi(\tau)\big(a(\tau)\big)^t\leqslant C\mathbb{I},\;
$$
 with a uniform in $\epsi$   constant $C>0$ depending only on $N$ and
  $I(\bar w^\epsi(\cdot))=I(v^{{\ep}}(\cdot))$ a.s.,
  where $v^{\ep}$ solves equation \eqref{v-equation1}${}_N$.

  Since  the family of the  processes
  $\Big\{\big(I(\bar w^\epsi(\cdot)), \bar w^\epsi(\cdot)\big)\Big\}_{\epsi\in(0,1]}$
  is tight in $C([0,T];  \mathbb R^n_+  \times \mathbb R^{2n})$, then
   taking if necessary a subsequence  $\epsi_i\to 0$ we achieve that  the process
   $\big(I(\bar w^{\epsi_i}(\cdot)), \bar w^{\epsi_i}(\cdot)\big)$ weakly
   converges in law in $C([0,T];  \mathbb R^n_+  \times \mathbb R^{2n})$ to a limit process
   $(\bar I(\tau), \bar w(\tau))$, $0\leqslant \tau\leqslant T$. Using the same arguments as in {\bf Step 3} of Appendix~\ref{app_proof_lem_not0} we conclude that $\bar w(\tau)$  admits    representation \eqref{ito-limit} (cf. equation \eqref{append}).
    So $(\bar I(\tau), \bar w(\tau))$  is a desired process.
\end{proof}

Now let $(\bar I(\tau), \bar w(\tau))$, $0\leqslant\tau\leqslant T$, be a continuous process as in Lemma~\ref{lm-ito-limit}. Fix any
 $\delta>0$. For the process $\bar  I(\tau)$, $0\leqslant\tau\leqslant T$,  we define stopping times
  $\tau_j^\pm\leqslant T$ such that $\cdots<\tau_j^-<\tau_j^+<\tau_{j+1}^-<\dots$, similarly to the stopping times
   in Appendix~\ref{app_proof_lem_not0}. Namely,
\begin{enumerate}
\item If $[\bar I(0)]\leqslant \delta$, then $\tau^-_1=0$; otherwise $\tau_0^+=0$.
\item If $\tau_j^-$ is defined, then $\tau_j^+$ is the first moment after $\tau_j^-$ when $[\bar I(\tau)]\geqslant 2\delta$ (if this never happens, then we set $\tau_j^+=T$; similar in the item below).
\item If $\tau_j^+$ is defined, then $\tau_{j+1}^-$ is the first moment after $\tau_j^+$ when $[\bar I(\tau)]\leqslant\delta$.
\end{enumerate}
We denote $\Delta_j=[\tau_j^-,\tau_j^+]$, $\Lambda_j=[\tau_j^+,\tau_{j+1}^-]$ and set
\[
\Delta^\delta=\cup \Delta_j, \;\;\; \Lambda^\delta=\cup \Lambda_j.
\]
See Figure \ref{pic1}.
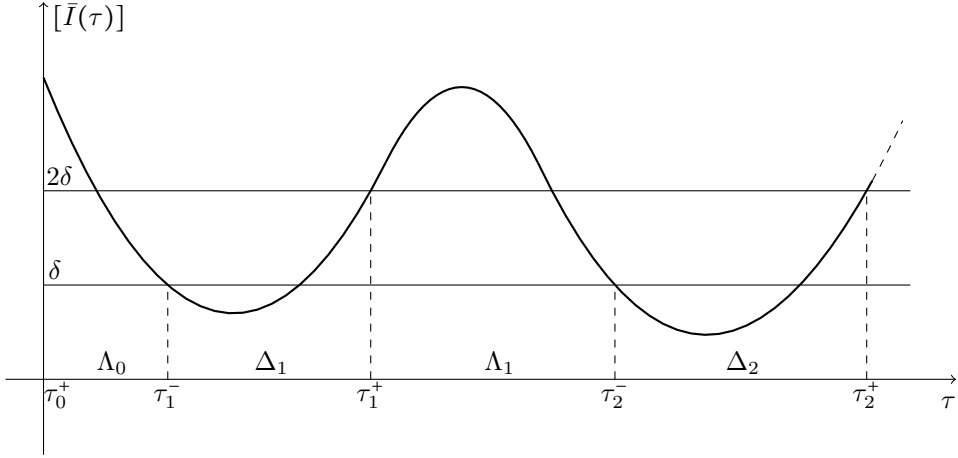
\begin{figure}[h]
{$\!\!\!\!\!\!\!\!\!\!\!\!\!\!\!\!\!$
\begin{tikzpicture}
\draw[, ->, black] (0,1) -- (12.5,1);
\draw[, ->, black] (0.5,0) -- (0.5,6);
%\draw [thick] (1.5,1.25) -- (3,1.25);
\draw[thin] (0.5,2.25) -- (11.9,2.25);
\draw[thin] (0.5,3.5) -- (11.9,3.5);
\draw[thin, dashed] (2.134,1) -- (2.134,2.25);
\draw[thin, dashed] (4.803,1) -- (4.803,3.5);
\draw[thin, dashed] (8.015,1) -- (8.015,2.25);
\draw[thin, dashed] (11.326,1) -- (11.326,3.5);
%\draw[thick, black, domain=0:1.5] plot (\x, {(\x-2)*(\x-2) +1});
%\draw[thick, black, domain=3:4.262] plot (\x, {2*(\x-3)*(\x-3) +0.5*(\x-3)+1.25});
\node at (12.4,0.7) {\color{black} $\tau$};
\node at (1.1,5.8) {\color{black} $[\bar I(\tau)]$
};
\node at (2.134,0.8) {\color{black} $\tau_1^-$};
\node at (4.803,0.8) {\color{black} $\tau_1^+$};
%\node at (0.65,1.2) {\color{black} $0$};
\node at (0.65,2.45) {\color{black} $\delta$};
\node at (0.71,3.7) {\color{black} $2\delta$};
\node at (0.69,0.8) {\color{black} $\tau_0^+$};
\node at (1.4,1.22) {\color{black} $\Lambda_0$};
\node at (3.5,1.22) {\color{black} $\Delta_1$};
\node at (6.5,1.22) {\color{black} $\Lambda_1$};
\node at (9.7,1.22) {\color{black} $\Delta_2$};
\node at (8.015,0.8) {\color{black} $\tau_2^-$};
\node at (11.326,0.8) {\color{black} $\tau_2^+$};
%\draw[, ->, black] (12,0) -- (12,5);
%\draw[, ->, black] (7.5,2) -- (14,2);
%\draw [thick] (11,2.25) -- (12,2.75);
%\draw [thick] (12,2.75) -- (13,3.75);
%\draw[thin, dashed, black, domain=0.1:0.5] plot (\x, {0.5*(\x-0.5)*(\x-0.5)-2.5*(\x-0.5) +5});
\draw[thick, black, domain=0.5:5] plot (\x, {0.5*(\x-0.5)*(\x-0.5)-2.5*(\x-0.5) +5});
%\draw[thin, dashed, black, domain=5:5.5] plot (\x, {0.5*(\x-0.5)*(\x-0.5)-2.5*(\x-0.5) +5});
\draw[thick,  black, domain=5:7] plot (\x, {-(\x-5)*(\x-5)+2*(\x-5) +3.875});
\draw[thick,  black, domain=7:11.4] plot (\x, {-0.008*(\x-7)*(\x-7)*(\x-7)+0.5*(\x-10)*(\x-10)+ 0.9*(\x-7) -0.625});
\draw[dashed,  black, domain=11.4:11.8] plot (\x, {-0.008*(\x-7)*(\x-7)*(\x-7)+0.5*(\x-10)*(\x-10)+ 0.9*(\x-7) -0.625});
%\draw[thick, black, domain=7.5:11] plot (\x, {0.25*(\x-10)*(\x-10) +2});
\end{tikzpicture}
}
%  \centering
 % \includegraphics[width=]{}
  \caption{A typical behaviour of the stopping times  $\tau_j^\pm$}\label{pic1}
\end{figure}
%}

By Theorem~\ref{th_ham_main},
\begin{equation}\label{short-delta}
\bE\,\mathcal{L}(\Delta^\delta)=o_{\delta}(1), \quad  \text{so} \;\;  \bE\,\mathcal{L}(\Lambda^\delta)= T- o_{\delta}(1),
\end{equation}
where $\mathcal{L}$ is the Lebesgue measure and   $o_{\delta}(1)$ goes to zero with $\delta$.
Due to the truncation \eqref{v-trivial} for  $\|v\|\geqslant N^{1/2}$, there exists $c(\delta,N)>0$ depending only on $N$ and $\delta$ such that  for each $j$,
\begin{equation}\label{finite-speed}
\mathbf{E}\mathcal{L}(\Lambda_j)\geqslant c(\delta,N)>0.
\end{equation}

For $j=0,1, \dots$ we will iteratively construct on  segments  $[0,\tau_j^\pm]$
 continuous  process $(I(\cdot), v(\cdot))$
 (defined on the space  $\Omega$ as in \eqref{Om})
 such that $I(\tau)=I(v(\tau))$ a.s. and  $\cD(I(\cdot))=\cD(I_N(\cdot))$. Moreover, on each segment $\Lambda_l\subset[0,\tau_j^\pm]$ the process $v(\tau)$, $\tau\in\Lambda_l$, will be a weak solution of \eqref{effective2}$_N$. Next we  will obtain a desirable measure $\mathbf{Q}_\delta$ as a limit of the laws of these processes as  $j\to\infty$.

For the sake of definiteness we assume that $[I(0)] > \delta$, so
 $0=\tau_0^+$.  With suitably chosen $\cB_\tau$-adapted Brownian motions $\{\tilde\bbeta_j,j=1,\dots,n\}$,
 we can assume that the  process $\bar I(\tau)$, $\tau\in[0,T]$ is a strong solution of \eqref{I-E2}$_N$.

a) Let  $\tau\in\Lambda_0$. Substituting  in \eqref{phi-E2}$_N$ \  $v=V_{\varphi(\tau)}(\bar I(\tau))$, we  get  for  the
$n$-vector  $\varphi(\tau)$ an equation, denoted
 by $(S_\varphi)$, with $\cB_\tau$-adapted coefficients and driven by $\{\tilde\bbeta_j\}$ (cf. equation \eqref{phi-E2}).
 Moreover, the drift term and the dispersion matrix are uniformly Lipschitz continuous in $\varphi$, where the Lipschitz constants may
 depend on $\delta$ and $N$. Hence, for any $\theta\in\mathbb{T}^n$
 there exists a unique  solution $\varphi(\tau)$, $\tau\in \Lambda_0$,  of $(S_\varphi)$ with $\varphi(0)=\theta$. Then by Lemma~\ref{N-S} the process $\tilde v^1(\tau)=V_{\varphi(\tau)}(\bar I(\tau))$, $\tau\in\Lambda_0$, is a weak solution of \eqref{effective2}$_N$. Obviously
 $I(\tilde v^1) \equiv \bar I$ 	and $\tilde v^1(0) = V_\theta(I_0)$.

b) Let $\tau \ge \tau_1^-$. For $(\xi_1, \xi_2) \in  (\mathbb{R}^2\setminus\{0\})\times (\mathbb{R}^2\setminus\{0\})$ we denote by
$U(\xi_1,\xi_2) \in \text{SO}(2)$ the  unique rotation of $\mathbb{R}^2$ that maps $\frac{\xi_2}{|\xi_2|}$ to $\frac{\xi_1}{|\xi_1|}$
(cf. Appendix~\ref{app_proof_lem_not0}).  Next for any two vectors
 $v^j=(\bv_1^j,\dots,\bv_n^j)\in  (\mathbb{R}^2\setminus\{0\})^n$,
 $j=1,2$, we set
 $$\mathcal{U}(v^1,v^2):=\text{diag}\{U(\bv_1^1,\bv_1^2),\dots, U(\bv_n^1,\bv_n^2)\}.$$
Consider   the process ${\hat v^1}(\tau):=\mathcal{U}\Big(\tilde v^1(\tau_1^-),\bar w(\tau_1^-)\Big)\bar w(\tau)$, $\tau\ge \tau_1^-$.
This still is an  It\^o process of the form as in  \eqref{ito-limit} with the same constant $C$. Moreover,
$$
{\hat v^1}(\tau_1^-)=\tilde v^1(\tau_1^{-}),
$$
 and  $I({\hat v^1}(\cdot))=\bar I(\cdot)$ a.s.  Now consider the continuous  process
$$
\big(\bar I^1(\tau),\bar v^1(\tau)\big): =\big(\bar I(\tau), \mathbf{1}_{\tau\leqslant\tau_1^-} \tilde v^1(\tau)+\mathbf{1}_{\tau>\tau_1^-}{\hat v^1}(\tau)\big), \; \tau\in[0,T],
$$
and  denote  $\bQ_{1,\delta}:=\cD\Big((\bar I^1(\cdot),\bar v^1(\cdot))\Big)\in \cP(\Omega)$.
Clearly, we have $\pi_I\circ \bQ_{1,\delta}=\cD(\bar I^1(\cdot))=\cD(I_N(\cdot))$, and $I(\bar v^1(\cdot))=\bar I^1(\cdot)$,
$\bQ_{1,\delta}$-a.s. Furthermore the process $\bar v^1(\cdot)$  solves \eqref{effective2}$_N$
weakly on  random interval  $[0,\tau_1^+]$, while on $\Delta_1$ and on the whole
interval $[\tau_1^+,T]$ it is an It\^o process of the form \eqref{ito-limit}.

c) For $\tau\in \Lambda_1$,  by the same method as in a) we construct a weak solution
$\tilde v^2(\tau)$, $\tau\in \Lambda_1$,  of  equation
\eqref{effective2}$_N$, equal  $\bar v^1(\tau_1^+)$  at $\tau = \tau_1^+$,
a.s. Consider the continuous process
$$\big(\bar I^2(\tau),\bar v^2(\tau)\big):=\big(\bar I(\tau),\mathbf{1}_{\tau\leqslant\tau_1^+}\bar v^1(\tau)+\mathbf{1}_{\tau_1^+<\tau\leqslant \tau_2^-}\tilde v^2(\tau)+\mathbf{1}_{\tau\geqslant\tau_2^- }{\hat v^2}(\tau)\big),\;\tau\in[0,T],$$
where ${\hat v^2}(\cdot):=\mathcal{U}\Big(\tilde v^2(\tau_2^-),\bar w(\tau_2^-)\Big)\bar w(\cdot)$.
Denote $\bQ_{2,\delta}=\cD\big((\bar I^2(\cdot), \bar v^2(\cdot)\big)$.
As in b),  $\pi_I\circ \bQ_{2,\delta}=\cD(\bar I^2(\cdot))=\cD(I_N(\cdot))$ and $I(\bar v^2(\cdot))=\bar I^2(\cdot)$, \, $\bQ_{2,\delta}$-a.s.
The process $\bar v^2(\cdot)$ solves \eqref{effective2}$_N$ on
 random intervals $\Lambda_j$, $j=0,1$ and is an It\^o process of the form \eqref{ito-limit} on $[0,T]\setminus \cup_{j=0}^1 \Lambda_j$.

d) Iteratively we construct on the space $\Omega$ measures
$\bQ_{j,\delta}$, $j\in\mathbb{N}$. Due to \eqref{finite-speed} we know  that a.s. the sequence $\tau_j^{\pm}$ stabilizes at $\tau=T$ after a finite (random) number of steps.  Accordingly, as $j\to\infty$,
 the sequence of measure $\bQ_{j,\delta}$ converges to a limiting measure $\bQ_\delta=\bQ_\delta^N$ on $\Omega$.

Let $(I_\delta(\tau), v_\delta(\tau))$, $\tau\in [0,T]$ be the natural process of the measure $\bQ_\delta$. We then have

 i) $\cD(I_\delta(\cdot))=\cD(I_N(\cdot))$;

 ii) $I(v_\delta(\cdot))=I_\delta(\cdot)$ $\bQ_\delta$-a.s.;

  iii) for $\tau\in \Lambda^\delta$ the process $v_\delta$ is a weak solution of \eqref{effective2}$_N$, while for
  $\tau\in \Delta^\delta$   $v_\delta(\tau) $ is  an It\^o process as in \eqref{ito-limit}, where $C$ does not depend on $N$.
  \medskip

{\bf Step 3}. Limit $\delta\to0$.

From the construction we know that the set of measures $\{\bQ^N_\delta, 0\leqslant\delta\leqslant1\}$ is tight.
Let $\bQ^N$ be any limiting measure as $\delta\to0$.  Then
\be\label{convmeas}
\bQ^N_{\delta_j} \strela \bQ^N \quad \text{as} \quad \delta_j\to0,
\ee
for some sequence $\{\delta_j\}$.
Since $\cD(I_\delta(\cdot))=\cD(\bar I(\cdot))=\cD(I_N(\cdot))$ \, $\forall \delta>0$, then $\pi_I\circ \bQ^N=\cD(I_N(\cdot))$. For the
 projection  of $\bQ^N$ to $\Omega_v$  we have

\begin{lemma}\label{lm-N-martingale}
The measure $\bP^N:=\pi_v\circ \bQ^N$ is a solution of  the martingale problem for equation   \eqref{effective2}$_N$.
\end{lemma}
\begin{proof}
Let us denote the drift terms and dispersion terms of \eqref{effective2}$_N$ as ${\mathbf{R}}_k^N(\tau,v)$ and $\tilde B^N_{kj}(\tau,v)$, respectively.
Consider the natural process $v_\delta=(\bv_{1,\delta},\cdots,\bv_{n,\delta})$ of the measure $\bP^N_\delta:=\pi_v\circ \bQ^N_\delta$
on $\Omega_v$. It  satisfies the system of It\^o  equations,
\begin{equation}\begin{split}
d\bv_k(\tau)&=\Big(\mathbf{1}_{\tau\in\Lambda^\delta} \bR^N_k(\tau,v)+\mathbf{1}_{\tau\in\Delta^\delta}B_k(\tau)\Big)d\tau\\
&\quad +\mathbf{1}_{\tau\in\Lambda^\delta}\sum_{j=1}^n\tilde B_{kj}^N(\tau,v)d\bbeta_j(\tau)+\mathbf{1}_{\tau\in\Delta^\delta}\sum_{j=1}^na_{kj}(\tau)d\mathbf{w}_j(\tau)\\
&=: {\bA}_k^\delta(\tau)d\tau+\sum_{j=1}^n\Big(G_{kj}^\delta(\tau,v) d\bbeta_j(\tau)+C_{kj}^\delta(\tau)d\mathbf{w}_j(\tau)\Big),
\end{split}
\quad k=1,\dots,n.
\end{equation}
Here $B(\tau)$ and $a(\tau)$ are the drift and dispersion, corresponding to the It\^o process $v_\delta$ on
 $\Delta^\delta$ (see item d)iii) of Step~2), so
    dispersion matrices $G_{kj}^\delta(\tau)$ and $C_{kj}^\delta(\tau)$ are supported by non-intersecting unions of
 random time-intervals $\Lambda^\delta$ and $\Delta^\delta$.
  Furthermore, for any $\delta>0$ and $k,m=1,\dots,n$, we have

i) the process $\gamma_k^\delta(\tau)=\bv_{k}(\tau)-\int_0^\tau {\bA}_k^\delta(s)ds\in\mathbb{R}^2$ is a $\bP^N_\delta$-martingale;

ii) the process $\Gamma_{km}^\delta=\gamma_k^\delta(\tau)\big(\gamma^\delta_m(\tau)\big)^t-\frac{1}{2}\int_0^\tau \big(X^\delta_{km}(s)+Y^\delta_{km}(s)\big)ds$, where
$$
X^\delta_{km}(s)=\sum_{j=1}^nG_{kj}^\delta(s)\big(G_{mj}^\delta(s)\big)^t, \qquad Y^\delta_{km}(s)=\sum_{j=1}^nC_{kj}^\delta(s)\big(C_{mj}^\delta(s)\big)^t,
$$
is a $\bP^N_\delta$-martingale.

 Note that for any $\delta>0$, by \eqref{short-delta},  we have
\begin{equation}\label{delta-0}
\begin{split}&\mathbf{E}^{\bP^N_\delta}\sup_{0\leqslant\tau\leqslant T}\Big|\int_0^\tau ({\bA}_k^\delta (s)-{\mathbf{R}}_k^N(s,v(s))ds\Big|\\
&\leqslant \bE^{\bP^N_\delta}\int_{\Delta^\delta}\Big(|{\mathbf{R}}_k^N(s,v(s))|+|B_k(s)|)ds\leqslant C(N)o_\delta(1).
\end{split}
\end{equation}

By \eqref{convmeas},   $\bP^N_{\delta_j}\strela\bP^N  := \pi_v(Q^N)$ as $\delta_j\to0$.

(1) We first show that the process $\gamma^0_k(\tau)=\bv_k(\tau)-\int_0^\tau {\mathbf{R}}_k^N(s,v(s))ds$ is a $\bP^N$-martingale.
Let us take any $0\leqslant\tau_1\leqslant\tau_2\leqslant T$ and $\Psi\in C_b(\Omega_v)$
such that $\Psi(\xi(\cdot))$ depends only on $\xi(\tau)$ with $\tau\in[0,\tau_1]$. We have to show that
\begin{equation}\label{martingale1}
\bE^{\bP^N}\Big(\big(\gamma_k^0(\tau_2)-\gamma_k^0(\tau_1)\big)\Psi(\xi)\Big)=0.
\end{equation}
The l.h.s equals
\[\begin{split}
&\lim_{\delta_j\to0}\bE^{\bP^N_{\delta_j}}\big(\gamma_k^0(\tau_2)-\gamma_k^0(\tau_1)\Psi(\xi)\big) \\
=&\lim_{\delta_j\to0}\bE^{\bP^N_{\delta_j}}\Big(\Psi(\xi)\big(\bv_k(\tau_2)-
\bv_k(\tau_1)-\int_{\tau_1}^{\tau_2}{\mathbf{R}}_k^N(s,v(s))ds\big)\Big)\\
=&\lim_{\delta_j\to0}\bE^{\bP^N_{\delta_j}}\Big(\Psi(\xi)\big(\int_{\tau_1}^{\tau_2}({\bA}_k^\delta(s)-{\mathbf{R}}_k^N(v(s)))ds\big)\Big)=0,
\end{split}\]
where in the second equality we use the fact that $\gamma^\delta_k$ is a $\bP^N_\delta$-martingale and in the third we use \eqref{delta-0}. So \eqref{martingale1} is established, and so  $\gamma_k^0(\tau)$ is a $\bP^N$-martingale.

(2) We then show the process $\Gamma_{km}^0(s)=
(\gamma_k^0(\tau))(\gamma^0_m(\tau))^t-\frac{1}{2}\int_0^\tau X_{km}^0(s)ds$ is a $\bP^N$-martingale, where $X_{km}^0(s)=\sum_{j=1}^n\tilde B_{kj}^N(s,v(s))\big(\tilde B_{mj}^N(s,v(s))\big)^t$.  Note that by the definition, for $\delta>0$
we have
\[\bE\sup_{0\leqslant\tau\leqslant T}\Big|\int_0^\tau (X^0_{km}(s)-X^\delta_{km}(s)-Y_{km}^\delta(s))ds\Big|\leqslant C(N)\bE^{\bP^N_\delta}\mathbf{1}_{ s\in \Delta^\delta}\leqslant C(N)o_{\delta}(1).
\]
Then  $\Gamma_{km}^0$ is a $\bP^N$-martingale due to the same reasoning as in (1).
This  finishes the proof of the lemma.
\end{proof}

{\bf Step 4}. Limit $N\to\infty$.

 By the construction and estimate \eqref{post-bound} the set of measures $\{\bQ^N, N\in\mathbb{N}\}$ is tight.
  Consider any limiting measure $\bQ$ for this family as $N\to\infty$.  Repeating in a simpler way the proof of
  Lemma~\ref{lm-N-martingale}, we find that $\pi_v\circ\bQ$ solves the martingale problem of \eqref{effective2}.
Therefore, it is  a weak solution for \eqref{effective2} with $\pi_v\circ \bQ(v(0)=V_\theta(I))=1$ and $I\circ\pi_v\circ\bQ=\cD(I^0(\cdot))$. Hence the  assertion of
Theorem~\ref{lifting-thm} is established.  \end{proof}

Theorem \ref{lifting-thm}, Proposition \ref{effective-unique-weak} and Theorem~\ref{th_ham_main} jointly imply the following.

\begin{theorem}\label{averaged-unique} Under  Assumption \ref{assumption-v-3},

i) For any $v_0\in\mathbb{R}^{2n}$ effective equation \eqref{effective1} has a unique strong solution $v(\tau;v_0)$, $\tau\in[0,T]$, equal $v_0$ at $\tau=0$. It satisfies
\begin{equation}\label{ef-v-bound}\mathbb{E}\sup_{0\leqslant \tau\leqslant T}\|v(\tau; v_0)\|^{2q_0}\leqslant C_{q_0}(|v_0|,T)<+\infty, \end{equation}
where $C_{q_0}(\cdot,\cdot)$ is the same as in Assumption \ref{assumption-v-3}.

ii) For any $v_0\in\mathbb{R}^{2n}$,  solution $v^{\eps}(\tau; v_0)$ of equation \eqref{v-equation1} with $v^{\eps}(0;v_0)=v_0$ satisfies
\begin{equation}\label{converge-i-v}
\cD\big(I(v^\eps(\cdot;v_0))\big)\rightharpoonup\cD\big(I(v(\cdot;v_0))\big)\; \text{ in }\;
\mathcal{P}(C[0,T],\mathbb{R}_+^n)\; \text{ as }{\eps}\to0.
\end{equation}
Moreover, the process $I^0(\tau):=I(v(\cdot;v_0))$, $\tau\in[0,T]$ is a weak solution of \eqref{averaged-1}, equal $I_0=I(v_0)$ at $\tau=0$.
 \end{theorem}

 \begin{remark}A straightforward analysis of the proof
 shows that it goes through without changes if $v^{{\eps}}(\tau; v_{{\eps}0})$ solves \eqref{v-equation1} with an
  initial data $v_{{\eps}0}$ which converges to $v_0$ as ${\eps}\to0$. In this case still
 \begin{equation}\label{conv_new}
 \cD\big(I(v^{\eps}(\cdot;v_{{\eps}0})\big)\rightharpoonup\cD\big(I(v(\cdot;v_0))\big)\; \text{ in }\; \mathcal{P}(C[0,T],\mathbb{R}_+^n)\; \text{ as }{\eps}\to0. \end{equation}
 \end{remark}
The result in Theorem~\ref{averaged-unique} admits an immediate generalization to the case when the initial data $v_0$ in  \eqref{converge-i-v} is a random variable:

\begin{amplification}\label{amp_ran_ini}  Let $v_0$ be a random variable independent of the Wiener processes $\bbeta_j(\tau)$, $j=1,\dots,n_1$. Then  still
convergence \eqref{converge-i-v} holds.
\end{amplification}
\begin{proof}   Let $v^\eps$ be a  weak solution of \eqref{v-equation1}  with $v^\eps(0)=v_0$.   Let
$(\Omega', \mathcal{F}',\mathbf{P}')$ be some probability space and $\xi_0^{\om'}$ be a r.v. on $\Omega'$, distributed
as $v_0$.
Then $v^{\eps\om}(\tau; \xi_0^{\om'})$ is a weak solution of \eqref{v-equation1},
 defined on the probability space $\Omega'\times\Omega = \{(\omega',\omega)\}$.
Take  $f$ to be  a bounded continuous function on $C([0,T], \mathbb{R}_+^n)$. Then by the theorem above,  for each $\omega'\in\Omega'$
%and $v_0^{\omega'}=\xi_0(\omega')$,
 $$
 \lim_{\eps\to0}\mathbf{E}^\Omega f\Big(I\big(v^{\eps \om} (\cdot ;\xi_0^{\omega'})\big)\Big)=\mathbf{E}^\Omega f\Big(I\big(v^{ \om}(\cdot; \xi_0^{\omega'})\big)\Big),
 $$
 where $v^{ \om}(\cdot; \xi_0^{\omega'})\big)$ is a weak solution of \eqref{effective1} with $v^\omega(0)=\xi_0^{\omega'}$.
Since $f$ is bounded, then by the Lebesgue dominated  convergence theorem we have
\[\begin{split}
\lim_{\eps\to0} \mathbf{E} f(I(v^{\eps}(\cdot; v_0)))=
\lim_{\eps\to0}\mathbf{E}^{\Omega'}\mathbf{E}^\Omega  f(I(v^{\eps\om}(\cdot;\xi_0^{\omega'})))\\
=\mathbf{E}^{\Omega'}\mathbf{E}^\Omega f(I(v^{\om}(\cdot ;\xi_0^{\omega'}))) = \mathbf{E} f(I(v(\cdot; v_0))).
\end{split}
\]
This implies the required convergence \eqref{converge-i-v}.
\end{proof}

\begin{proposition}\label{p_inirotate}
1) A weak solution  $v^\theta$, $\theta\in \T^n$,  as in Theorem \ref{lifting-thm}   may be chosen to be
$
v^\theta(\tau) = \Phi_\theta  v(\tau; v_0),
$
where $v(\tau; v_0)$ is the strong solution from Theorem~\ref{averaged-unique}.

2) More generally, if $\tilde\tau$ is a non-negative constant and $\theta\in \T^n$ is a r.v., measurable with respect to
$\cF_{\tilde\tau}$ (see \eqref{OmOm}), then the process
$
\hat v(\tau) = \Phi_\theta v(\tau;v_0)$,  $ \tau\ge \tilde\tau,
$
is a weak solution of equation \eqref{effective1}.
\end{proposition}
\begin{proof}
It suffices to prove 2) since it implies 1) if we choose $\tilde\tau=0$.  Substituting in  \eqref{effective1}
$
v(\tau) = \Phi_{-\theta} \hat v(\tau)
$
we get that
$$
d\hat v(\tau) = \Phi_\theta R\big( \Phi_{-\theta} \hat v(\tau)\big) d\tau+
 \Phi_{\theta} \llan B\rran\big( \Phi_{-\theta} \hat v(\tau)\big) d \beta(\tau).
$$
Or, using Proposition \ref{square-root}.ii), that
$$
d\hat v(\tau) =  R(\hat v(\tau)) d\tau+
  \llan B\rran (\hat v(\tau))  \Phi_{\theta}    d\beta(\tau).
$$
Since the r.v. $\theta$ is $\cF_{\tilde\tau}$-measurable, then the process
$
t\mapsto  \Phi_\theta \big( \beta(t+\tilde\tau) -  \beta(\tilde\tau) \big)
$
is a standard Wiener process in $\R^{2n}$. Thus, $\hat v(\tau)$, $\tau\ge \tilde\tau$, is a weak
solution of \eqref{effective1}.
 \end{proof}

\section{Equations in
 bounded domains}\label{ss_local_version}

In this section we consider problem
\eqref{v-equation1}
in the set
$$
\fB=
\{v \in\mathbb R^{2n}\,:\, I(v)\in {B}\},
$$
 where ${B}= B_R(\R^n)$  for some $R>0$.  We assume  that all  coefficients in
 \eqref{v-equation1} are defined and  Lipschitz continuous on  the set $\bar\fB$ (so $W(I)$ is defined and  Lipschitz continuous on $\bar B$).
We also  assume that  conditions of items (1)--(2)  of Assumption~\ref{assumption-v-3} are fulfilled  on $\bar\fB$.
 Let us consider the effective equation \eqref{effective1} as in Section~\ref{construct-effective} on set $\bar\fB$ (note that to calculate
 coefficients of the equation on $\bar\fB$ it suffices to know $P, B$ and $W$ only where they are defined).

Let $v^\eps(\tau)$ be a solution of  \eqref{v-equation1} with the $v_0$ as above.
Denote by $\tau_{R}^\eps$ its exit time from  domain $\fB$.
%$\{v\in\mathbb R^{2n}\,:\, I(v)\in\bar{B}\}$.
Then $\tau_{R}^\eps = \inf\{\tau>0\,:\, I^\eps(v(\tau))\in \p B\}$,
where $I^\eps(\tau) :=  I^\eps(v(\tau))$ satisfies  \eqref{I-equation1}.
 Similarly, let $\tau_{R}^{\rm h}$ stands for the exit time from $\fB$ of a solution $v(\tau)$ of equation \eqref{effective1}, equal $v_0$ at $\tau=0$.
 Again, it equals the exit time of $I(v(\tau))$ from $B$.

 \begin{theorem}\label{th_ham_loc} Under the above assumptions, for any $T>0$ and any
 $v_0\in \fB$, as $\eps\to0$ the family  of processes
$\{I^\eps(\cdot\wedge\tau_{R}^\eps)\}$  converges in law, weakly  in the space $(C([0,T];\mathbb R^n),\mathcal{B}_T)$, to a weak solution
$I^{\rm h}(\cdot\wedge\tau^{\rm h}_{R})$ of problem \eqref{averaged-1}, \eqref{initial-a1}. The latter solution is
 obtained as the
action-vector for a unique weak solution of the effective equation \eqref{effective1},  equal $v_0$ at $\tau=0$ and
stopped at $\partial {\fB}$.
\end{theorem}

\noindent
\begin{proof} The required statement is essentially a consequence of Theorems~\ref{th_ham_main} and \ref{lifting-thm}. Indeed, by
 Lemma~5.2  in \cite{GHKu}  coefficients ${\mathbf{P}}_k(v)$ and $B_{kj}(v)$
can be extended from the set $\bar\fB$
%$\{v\in\mathbb R^{2n}\,:\, I(v)\in \bar{B}\}$
 to the whole space $\mathbb R^{2n}$ in such a way that
the extensions  are bounded,  Lipschitz continuous, and the extended  matrix $B(v) B(v)^t$ is positive definite.
Using the same lemma we also extend $W(I)$ to a Lipschitz continuous vector-function $\widetilde W$ on
$\mathbb R^d$ with   compact support.

Consider the function $\Gamma$ on $\R^n$,
$$
\Gamma(x) =
\begin{cases}
%\int_R^{|x|}
\exp\big(-\frac1{|x|-R}\big),\qquad &|x| > R,
\\
0, \qquad &|x| \le R.
 \end{cases}
$$
It is smooth, globally Lipschitz, and is flat on $\partial B$.  For $\alpha \in \R$ define vector-functions $W_\alpha$ on
$\R^n$ as  $W_\alpha(x)= \alpha \nabla \Gamma(x)+\widetilde{W}(x)$. All of them are continuations of $W$ to $\R^n$.

\begin{lemma}\label{l_w_extension}
 There is at most countable number of $\alpha$'s for which the components of $W_\alpha$ are rationally dependent on a set of positive measure.
\end{lemma}
\begin{proof}
Since $W_\alpha = W$ on $\bar B_R$, then in view of Assumption~\ref{assumption-v-3}.(1)
 it suffices to examine intersections of the sets of rational dependence
of  components of $W_\alpha$ with $\R^n \setminus B$.

  Assume that there exists two distinct $\alpha_1$,  $\alpha_2$ and a non-zero  vector $m\in \Z^n$
   such that
  $$
  \textstyle
  {\mathcal L} \big(
  \{x\in \mathbb R^n \setminus B\,: m \cdot W_{\alpha_1} =0 \
   \hbox{and }  m \cdot W_{\alpha_2} =0
 \} \big)>0.
  $$
  Then $ m \cdot \nabla \Gamma(x)=0$ on a set of positive measure in $\mathbb R^n\setminus B$. But
  $\Gamma (x)$ may be written as $f(|x|^2)$, where $f(r)$ is a smooth function, vanishing for $r \le  R^2$. Then
    $
  \nabla \Gamma(x) = 2x f'(|x|^2)
  $
  and we see that the set under discussion has zero measure since $f'(r)>0$ for $r>R^2$.
  Therefore for each  $\delta>0$ and any non-zero  vector   $m\in \Z^n$  the number of
  $\alpha$'s for which $ {\mathcal L}  \big( \{x\in\mathbb R^d\setminus B\,: m\cdot \nabla W_\alpha
  =0\}   \big)>\delta$ is at most countable.    This implies the assertion.
\end{proof}

Denote $X_T = C([0,T], \mathbb R^n)$,  take any number $\alpha_0$, different from the countable family in
the lemma above, and choose $W_{\alpha_0}$ for the extentson of $W$.
Then by  Proposition~\ref{l_suff_assum_i4} \ Assumption~\ref{assumption-v-3}.(4) holds, and so
 Theorems~\ref{th_ham_main}, \ref{lifting-thm} and \ref{averaged-unique}
 apply to the obtained stochastic equation
in $\R^{2n}$. Thus    for any $T>0$ the corresponding process $I^\eps$ converges in law
in $X_T$ as $\eps\to0$
 to a solution of the averaged  equation   \eqref{averaged-1} and may be lifted to a solution of the corresponding effective equation.
The initial condition remains  unchanged.

Let $\tau_R(I)=\min\big(T, \inf\{\tau\in (0,+\infty)\,:\,I(\tau)\in\partial B \}\big)$.
For an arbitrary bounded continuous functional $\mathcal{R}$ on $X_T$, consider there another  functional
$$
\mathcal{R}_B(I)=\mathcal{R}\big(I(\tau\wedge\tau_R(I))\big), \quad I \in X_T.
$$
It is not  continuous. However,  the following statement ensures that it is almost surely continuous
with respect to the measure on $X_T$,
generated by the limit process $I^{\rm h}$, constructed in Theorem~\ref{averaged-unique}
 (and called there $I^0$).
\begin{lemma}\label{L_a_s_cont}
Under the standing assumptions,  let $v(\tau)$ be any  solution of the corresponding effective equation in $\R^{2n}$
 and $I^h(\tau)=I(v(\tau))$. Then
$$
\mathbf{P}\big(\inf\{\tau>0\,:\, I^{\rm h}(\tau)\in\partial B\}= \inf\{\tau>0\,:\, I^{\rm h}(\tau)\not\in \bar{B}\}\big)=1.
$$
\end{lemma}
\begin{proof}
The desired inequality is an immediate consequence of the fact  that   the diffusion coefficient of the process $v(\tau)$ at $\partial B$
 does not degenerate in the direction of a normal vector to $\partial B$.
\end{proof}
Combining this with the statements of Theorems~\ref{th_ham_main}, \ref{lifting-thm}
 and \cite[Theorem~5.2]{Bill} we conclude that the law of
$I^\eps(\tau\wedge \tau^\eps_R)$ converges to that of $I^{\rm h}(\tau\wedge\tau^{\rm h}_R)$.
\end{proof}

It should be noted that the distribution of $\tau^\eps_R$ does not concentrate in the vicinity of zero, as $\eps\to0$. So
Theorem~\ref{th_ham_loc} with high probability describes the behaviour of solutions for \eqref{v-equation1} on $t$-time intervals
of order $\eps^{-1}$.
More precisely, %as a consequence of the tightness of the family $\{I^\eps(v(\cdot))\}$ 
we have
%\begin{proposition}\label{l_loubou}
%Let $\|v_0\|<R$. Then for any $\delta>0$ there exists $\bar s=\bar s(\delta, R-\|v_0\|)>0$ such that
%$$
%\mathbf{P}\{\tau_R^\eps < \bar s\}<\delta, \quad \forall \eps\in(0,1].
%$$
%\end{proposition}

\begin{proposition}\label{l_loubou}
Let $\|v_0\|<R-r$ with $0<r<R$.  %for any $\delta>0$ there exists $\bar s=\bar s(\delta, R-\|v_0\|)>0$ such that
Then there exists a constant $C=C(R)$ such that, for all $T\in [0,1]$,
$$
\mathbf{P}\{\tau_R^\eps <T\}\leqslant \frac{C(R)(T+\sqrt{T})}{r^2}, \quad \forall \eps\in(0,1].
$$
\end{proposition}

\begin{proof}
As was explained in the beginning of the proof of Theorem \ref{th_ham_loc}, we may assume that the coefficients
 ${\mathbf{P}}_k(v)$ and $B_{kj}(v)$
are extended from the set $\bar\fB$
%$\{v\in\mathbb R^{2n}\,:\, I(v)\in \bar{B}\}$
 to the whole space $\mathbb R^{2n}$ in such a way that
the extensions  are bounded, uniformly  Lipschitz continuous, and the extended  matrix $B(v) B(v)^t$ is positive definite.
The function $W(I)$ also admits an extension to $\mathbb R^d$, the extended function is  bounded and uniformly Lipschitz continuous.
Then a solution of \eqref{v-equation1} is well defined for all $\tau>0$, and by the It\^o formula we have
\begin{equation}\label{new_p1}
  d\|v\|^2(\tau)=2\sum_{k=1}^n\bv_k^t{\mathbf{P}}_k(v)d\tau+\sum_{k=1}^n\sum_{j=1}^{n_1}\|B_{kj}(v)\|^2_{HS}d\tau+
2 \sum_{k=1}^n \sum_{j=1}^{n_1}\bv_k^tB_{kj}(v)d\bbeta_j(\tau), 
\end{equation} 
Therefore,  
\begin{eqnarray}
&\mathbf{E} \big( \sup\limits_{0\leqslant \tau\leqslant T} \big| \|v\|^2(\tau)-\|v_0\|^2\big|\big)
\nonumber \\ 
\label{new_p2}
&\leqslant \mathbf{E} \Big( \sup\limits_{0\leqslant \tau\leqslant T}\Big|\int_0^\tau 
\big[2\sum_{k=1}^n\bv_k^t(s){\mathbf{P}}_k(v(s))+\sum_{k=1}^n\sum_{j=1}^{n_1}\|B_{kj}(v)\|^2_{HS}(s)\big]ds\Big|
\Big)\\
\nonumber
&+\mathbf{E} \Big( \sup\limits_{0\leqslant \tau\leqslant T}\Big|\int_0^\tau
 2 \sum_{k=1}^n \sum_{j=1}^{n_1}\bv_k^t(s)B_{kj}(v(s))d\bbeta_j(s) \Big|
\Big)
\end{eqnarray}
Using estimate  \eqref{apriori-1} and evoking the BDG inequality to bound the second term in the r.h.s. of  \eqref{new_p2} we find that 
%\begin{equation}\label{new_p6}
$$
\mathbf{E} \big( \sup\limits_{0\leqslant \tau\leqslant T} \big| \|v\|^2(\tau)-\|v_0\|^2\big|\big) \le C(R)(T+\sqrt{T}), 
$$
for all $T\in [0,1]$. 
%\begin{equation} \label{new_p3}
%\mathbf{E}\sup_{0\leqslant \tau\leqslant 1} \|v\|^2(\tau)\leq C_0(R)\qquad \forall \eps\in(0,1].\end{equation}
%Using the boundedness of ${\mathbf{P}}_k$ and $B_{kj}$ we conclude that
%\begin{equation}\label{new_p4}
%\mathbf{E} \Big( \sup\limits_{0\leqslant \tau\leqslant T}\Big|\int_0^\tau
%\big[2\sum_{k=1}^n\bv_k^t(s){\mathbf{P}}_k(v(s))+\sum_{k=1}^n\sum_{j=1}^{n_1}\|B_{kj}(v)\|^2_{HS}(s)\big]ds\Big|
%\Big)\leqslant C_1(R)T\end{equation}
%for all $T\in [0,1]$.  The second term on the right-hand side of \eqref{new_p2} can be estimated with the help of the BDG inequality as follows:
%\begin{equation*}\label{new_p5}
%  \mathbf{E} \Big( \sup\limits_{0\leqslant \tau\leqslant T}\Big|\int_0^\tau
 %2 \sum_{k=1}^n \sum_{j=1}^{n_1}\bv_k^t(s)B_{kj}(v(s))d\bbeta_j(s) \Big|
%\Big)\leqslant  C\Big(\mathbf{E} \int_0^T\|v\|^2(s)ds\Big)^\frac12\leqslant C_2(R)\sqrt{T}  
%\end{equation*}
%Combining  \eqref{new_p2},  \eqref{new_p3}  and \eqref{new_p5} yields
%Combination of  the last inequality and \eqref{new_p2} and  \eqref{new_p3} yields that
%\begin{equation}\label{new_p6}
%\mathbf{E} \big( \sup\limits_{0\leqslant \tau\leqslant T} \big| \|v\|^2(\tau)-\|v_0\|^2\big|\big)\leqslant
%C(R)(T+\sqrt{T}), \end{equation}
%for all $T\in [0,1]$.  
Since $\{\omega\in\Omega\,:\,\tau_R\leqslant T\}\subset\{\omega\,:\, 
 \sup\limits_{0\leqslant \tau\leqslant T} \big| \|v\|^2(\tau)-\|v_0\|^2\big|\geqslant 2Rr-r^2\}$, then 
 from  this relation 
 %\eqref{new_p6}  
 by  Chebyshev's  inequality we obtain 
 $$
 \mathbf{P}\{ \tau_R\leqslant T\}\leqslant  C(R)\frac{T+\sqrt{T}}{2Rr-r^2}\leqslant C(R)\frac{T+\sqrt{T}}{r^2}.
 $$
\end{proof}

%\begin{proof}
%By Theorem~\ref{th_ham_main}
%the laws for the process $\{I(v^{\eps}(\tau;v_{0})),\tau\in[0,T]\}$, $0<\eps\leqslant 1$, are tight in $\mathcal{P}(C([0,T],\mathbb{R}_+^d))$.
%So by  Prokhorov's theorem for any $\delta>0$ there exists
%a compact set $\mathcal{M}_\delta$ in $C([0,T],\mathbb{R}_+^d)$  such that
%$\mathbf{P}\{I(v(\cdot))\in\mathcal{M}_\delta\}>1-\delta$.
%Now the required statement  is a consequence of the equicontinuity of functions
%from $\mathcal{M}_\delta$, granted by the Arzel\`a--Ascoli theorem.
%\end{proof}

\section{Mixing and  uniform convergence}\label{s_9}
In this section we establish  the uniform in time convergence  in distribution  of the actions of solutions for equation
 \eqref{v-equation1} to those for  solutions of  effective equation~\eqref{effective1},
 with respect to the dual-Lipschitz metric  (see Definition~\ref{d_dual-Lip})
in the space of  probability measures. The proof uses the approach developed in \cite{HGK22, GHKu}, where a similar result was obtained in the easier case when the frequency vector $W$ in equation \eqref{v-equation1} is
constant (cf. Remark~\ref{r_111}.1)).

\begin{proposition}\label{p_ran_ini}  Under the assumption of Amplification \ref{amp_ran_ini}
 let  the r.v. $v_0$ be  such that $\|v_0\| \le R$ a.s.,
for some $R>0$. Then the rate of convergence in \eqref{converge-i-v} with respect to the dual-Lipschitz distance depends only on $R$.
\end{proposition}

\begin{proof}
The proof of  Amplification \ref{amp_ran_ini}   shows that it suffices to verify that for a non-random initial vector $v_0 \in \bar B_R(\mathbb{R}^{2n})$ the rate of convergence in \eqref{converge-i-v} depends only on $R$.
 Assume the opposite. Then there exist a $\delta>0$, a sequence
$\eps_j\to0$ and vectors $v_j\in \bar B_R(\mathbb{R}^{2n})$ such that
\begin{equation}\label{big-delta}\|\cD\big(I(v^{\eps_j}(\cdot; v_j))\big)-\cD\big(I(v^0(\cdot; v_j))\big)\|_{L,C([0,T],\mathbb{R}^{n}_+)}^*\geqslant \delta.\end{equation}
By \eqref{apriori-v-2} and \eqref{ef-v-bound}, using  the same argument as in the proof of Lemma~\ref{l_compa}, we know that the two sets of probability measures
$\{\cD\big(I(v^{\eps_j}(\cdot; v_j))\big)\}$ and $\{\cD\big(v^0(\cdot; v_j)\big)\}$ are tight, respectively
 in  $\mathcal{P}(C([0,T],\mathbb{R}_+^n))$ and $\mathcal{P}(C([0,T],\mathbb{R}^{2n}))$. Therefore, there exists a sequence $k_j\to\infty$ such that  $\eps_{k_j}\to0$, $v_{k_j}\to v_0$,
$$
\cD(I(v^{\eps_{k_j}}(\cdot; v_{k_j})))\rightharpoonup  Q_0^I\;\text{ in }\;\cP(C([0,T],\mathbb{R}_+^n)),$$
and
$$ \cD(v^{0}(\cdot; v_{k_j}))\rightharpoonup Q^v_0\; \text{ in }\; \cP(C([0,T],\mathbb{R}^{2n})).
$$
Then due to \eqref{big-delta},
\begin{equation}\label{contr1}
 \;\| Q_0^I-I\circ Q^v_0\|_{L,C([0,T],\mathbb{R}_+^n)}^*\geqslant\delta.
\ee
Since in the well-posed eq. \eqref{effective1}  the drift and  dispersion   are  locally Lipschitz and its solutions satisfy
estimates  \eqref{ef-v-bound},  then  the law
$\cD(v^0(\cdot;v'))$ is continuous with respect to the law of the
initial condition $v'$. Therefore the limiting measure
$Q_0^v$ is the unique weak solution of the effective equation \eqref{effective1} with initial condition $v^0(0)=v_0$. By \eqref{conv_new} the measure
$ Q_0^I$ equals  $I\circ Q_0^v$.  This   contradicts  \eqref{contr1} and proves  the assertion.
\end{proof}

In this section, we make the following assumption concerning system
 \eqref{v-equation1} and the corresponding effective equation \eqref{effective1}.

   \begin{assumption}\label{a_8.1}
   The first three items (1)-(3)     of Assumption~\ref{assumption-v-3} hold, and
   \begin{enumerate}
   \item [(4)] For any $v_0\in\mathbb{R}^{2n}$ a unique strong solution $v^{\eps}(\tau;v_0)$ of \eqref{v-equation1} is such that for some
   $q_0>(q\vee4)$ we have
   \begin{equation}\label{v-upper-2}
   \mathbf{E}\sup_{T'\leqslant\tau\leqslant T'+1}\|v^{\eps}(\tau;v_0)\|^{2q_0}\leqslant C_{q_0}(\|v_0\|),
   \end{equation}
   for every $T'\geqslant0$ and ${\eps}\in(0,1]$, where $C_{q_0}$ is a continuous non-decreasing function.
\item[$(5)$]   Effective equation \eqref{effective1} is  mixing with a stationary  measure $\mu^0$ and a (strong) stationary solution
$v^{st}(\tau)$, $\tau\ge0$.

\item[$(6)$]
For any its solution
 $v(\tau)$, $\tau\ge0$, such that $\cD(v(0))=:\mu$ and  $ \lan  \|v\|^{2q_0}, \mu(dz) \ran = \EE \|v(0)\|^{2q_0}
  \le M'$ for some $M'>0$ (we recall  notation \eqref{recall})  we have
\begin{equation}
\label{mixing-quantify}
\|\cD(v(\tau))-\mu^0\|_{L, \R^{2n}}^*\leqslant g_{M'}(\tau, d)
%)\|\mu-\mu^0\|_{L}^*,
\quad \forall \tau\geqslant0, \;
 \text{ if } \| \mu- \mu^0\|_{L,\mathbb{R}^{2n}}^*\leqslant {d}\leqslant2.
\end{equation}
Here the function
$$
g: \R_+\times [0,2]\times  \R_+  \ni\!  (\tau, d, M) \mapsto g_M(\tau, d),
$$
is continuous,  vanishes with $d$, converges to zero when  $\tau\to\infty$
and is such that  for each fixed
$M\ge 0$ the function $ (\tau, d)\mapsto  g_M (\tau, d)$
 is uniformly continuous in $d$, for $(\tau, d) \in [0, \infty)\times[0,2]$. \footnote{So $g_M$
extends to a continuous function on $[0, \infty]\times [0,2]$ which vanishes when $\tau=\infty$ or $d=0$.}
 \end{enumerate}
  \end{assumption}

  We emphasize that   we assume the  mixing only for  the effective equation \eqref{effective1},  but not  for the
   original equation~\eqref{v-equation1}.  Since
  Assumption~\ref{a_8.1} implies   Assumption~\ref{assumption-v-3}, then the assertions of Section~\ref{lifting} with any $T>0$ hold
   for solutions of   equation  \eqref{v-equation1} which we analyse in this section.

Assumption (6) above    may seem rather restrictive. But it is not, as shows the next result:

  \begin{proposition}\label{p_suff_cond}
 If we keep all the conditions in Assumption \ref{a_8.1} except (6)
  and assume that  for each $M>0$ and any $v^1, v^2 \in \bar B_M(\mathbb{R}^{2n})$ we have
  \be\label{x1}
  \| \cD (v(\tau ; v^1)) - \cD (v(\tau ; v^2)) \|_{L,\mathbb{R}^{2n}}^*  \le {\mathfrak g}_M(\tau ),
  \ee
  where
   ${\mathfrak g}$ is a non-negative  continuous function of $(M,\tau)\in \R_+^2$  which
     goes to zero when $\tau \to\infty$ and is  non-decreasing  in $M$,
     then (6) holds with a suitable function $g$.
  \end{proposition}
    For a proof of the proposition we refer the reader to \cite[Section 7.1]{GHKu}.

  Note that \eqref{x1} holds (with ${\mathfrak g}$ replaced by $2{\mathfrak g}$) if
   $$
  \| \cD v(\tau ; v) - \mu^0 \|_{L,\mathbb{R}^{2n}}^*  \le {\mathfrak g}_M(\tau )\quad \forall\,  v \in \bar B_M(\mathbb{R}^{2n}) .
  $$
  Usually a proof of mixing for eq.  \eqref{effective1} in fact establishes the estimate above.
  So, given assumptions (1)-(5),
   condition (6)  is a rather mild restriction.

   \begin{theorem}\label{thm-uniform}
  Under Assumption \ref{a_8.1},
   for  any $v_0\in\mathbb{R}^{2n}$
     \[
  \lim_{\eps\to0}\,
  \sup_{\tau\geqslant0}\|I\circ\cD(v^{\eps}(\tau;v_0))-I\circ\cD(v^{0}(\tau;v_0))\|_{L,\mathbb{R}_+^n}^* =0,
  \]
  where $v^{\eps}(\tau;v_0)$ and $v^{0}(\tau;v_0)$ solve respectively \eqref{v-equation1} and \eqref{effective1} with the same
   initial condition $v_0$.
  \end{theorem}
%We first prove  the following lemma.

\begin{proof}
Below we abbreviate $\| \cdot\|_{L, \R^n_+}^*$ to $\| \cdot\|_{L}^*$.
Since $v_0$ is fixed, we  also  abbreviate $v^\eps(\tau; v_0)$ to  $v^\eps(\tau)$.
Due to \eqref{v-upper-2}
\be\label{M}
 {\EE}\| v^{\ep}(\tau)\|^{2q_0}\leqslant C_{q_0}(\|v_0\|) =:M^*    \quad \forall\, \tau\ge0.
\ee
By \eqref{M} and   \eqref{converge-i-v} also
\be\label{8.05}
 {\EE}\| v^{0}(\tau;v_0)\|^{2q_0} =
\lan \|v\|^{2q_0}, \cD ( v^0({\tau};v_0))\ran \le M^* \quad \forall\, \tau\ge0.
\ee
Since $\cD( v^0({\tau};0)) \strela \mu^0$ as ${\tau}\to\infty$, then  from the estimate above with $v_0=0$  we get that
\be\label{x0}
\lan \|v\|^{2q_0}, \mu^0 \ran \le C_{q_0} (0)\le M^*.
\ee
%{\color{blue}
%Probably we do not need the sentence below and it should be removed:
%For later usage we note that
%since to derive  estimates  \eqref{8.05} and
%\eqref{x0} we only used Assumptions~\ref{a_8.1}.(1)-(4) and the fact that eq.~\eqref{effective1} is mixing,
%then the two estimates hold under the assumptions of Proposition~\ref{p_suff_cond}. }
%%We recall that $v^{st}(\tau), \tau\ge0$, is a strong stationary solution of \eqref{effective1}, $\cD v^{st}(\tau) \equiv \mu^0$.

The  constants in estimates below depend on $M^*$, but usually this   dependence  is not indicated.
For any $T\ge0$ we denote by $v_T^0(\tau)$, $\tau\ge0$,  a weak solution of  effective equation
\eqref{effective1} such that
\be\label{v^0_T}
\cD (v^0_T(0) )= \cD (v^\eps(T)).
\ee
Thus $v^0_T(\tau)$ depends on $\eps$, and    $v^0_0(\tau) = v^0(\tau; v_0)$.

  Below in this proof by $\kappa^j(\cdot)$, $j=1,\dots,5$, we denote various monotonically increasing continuous functions $\mathbb{R}_+\to\mathbb{R}_+$, vanishing at zero and positive outside it.  These
  functions depend only on the constants in Assumption~\ref{a_8.1} and functions $\nu$ in
  \eqref{nu_Lip}, $C$ in \eqref{v-upper-2} and $g$ in \eqref{mixing-quantify}.
    A function $\kappa^j$, $j\ge2$,     may be constructed in terms of functions $\kappa^l$ with $l<j$.

  It is straightforward to see that the proof of Lemma~\ref{l_key_ham} implies that the process $I^{\eps}(\tau)=I(v^{\eps}(\tau;v_0))$
  satisfies the estimate below (recall that    $[I]=\min_{1\leqslant j\leqslant n} I_j$)
  \begin{equation}\label{i-eps-small}
  \int_{T'}^{T'+1}\mathbf{P}\big(\{[I^{\eps}(\tau)]<{\gamma}\}\big)\leqslant \kappa^1({\gamma}),\qquad \forall T'\geqslant0, \;\ \forall {\gamma}, \eps\in(0,1],
  \end{equation}
  for some fixed function $\kappa^1(\cdot)$, and that the action-vector of the stationary solution $v^{st}(\tau)$ of equation \eqref{effective1} also meet this estimate. Then
  \begin{equation}
  \label{i-st-small}
  \mathbf{P}\big(\{[I(v^{st}(\tau))]<{\gamma}\}\big)<\kappa^1({\gamma}),\; \forall \tau\geqslant0, \;\forall {\gamma}\in(0,1].
  \end{equation}
  We will say that a  moment of time  $\tau\geqslant 0$ is $({\gamma},{\eps})-typical$, where
   ${\eps},{\gamma}\in(0,1]$, if $\mathbf{P}\big(\{[I^{\eps}(\tau)]<{\gamma}\}\big)\leqslant \kappa^1({\gamma})$.
  In view of \eqref{i-eps-small}, we have

  \begin{lemma}\label{typical-time}For each $\gamma, {\eps} \in(0,1]$, every interval  $[T',T'+1]$, $T'\geqslant0$, contains
  an $(\gamma,{\eps})$-typical moment of time $\tau=\tau(T', \gamma, \eps)$.
  \end{lemma}
  We continue with  two technical lemmas, needed to prove  principal Lemma~\ref{random-initial-average}.  By \eqref{M} and \eqref{x0} and Chebyshev's
  inequality,   for   any $\tau\ge0$ and $R>0$,
  \be\label{kappa2}
  \lan \cD v^\eps(\tau), B_R\ran, \ \lan \mu^0, B_R\ran \ge 1- R^{-2q_0} M^* =: 1- \kappa^2(R^{-1}).
  \ee

  \begin{lemma}\label{i-v-measure}
  For any $\eps\in(0,1]$ and  $\tau\ge0$,
 % If $m_1,m_2\in\mathcal{P}(\mathbb{R}^{2n})$ are such that $m_j(B_R)\geqslant1-\kappa^2(R^{-1})$ for each $R>1$, $j=1,2$, then
  \begin{equation}\label{i-v-measure-q}
  \| \cD I(v^\eps(\tau)) - I\circ \mu^0\|_{L}^* \leq\kappa^3\big( \| \cD v^\eps(\tau) -  \mu^0\|_{L,\mathbb{R}^{2n}}^*\big)
  \end{equation}
  (see Definition \ref{d_dual-Lip} for the distance $\|\cdot\|_L^*$ and the norm $| \cdot |_L$, used below).
  \end{lemma}
  \begin{proof}
  Let us abbreviate $\cD v^\eps(\tau) =: m$.  For $R\geqslant2$, consider a function $G_R(t)$ on $\mathbb{R}_+$ as in Figure \ref{GR2}.
 \begin{figure}[ht]
  \centering

  \begin{tikzpicture}[x=0.75pt,y=0.75pt,yscale=-1,xscale=1]
\label{GR1}
%Shape: Axis 2D [id:dp6073793338608799]
\draw  (182,220) -- (440,220)(207.8,103) -- (207.8,233) (433,215) -- (440,220) -- (433,225) (202.8,110) -- (207.8,103) -- (212.8,110)  ;
%Straight Lines [id:da6616550575886496]
\draw    (208,150) -- (332,150) ;
%Straight Lines [id:da1683662328152976]
\draw    (332,150) -- (405,220) ;
%Straight Lines [id:da6485272141276033]
\draw  [dash pattern={on 0.84pt off 2.51pt}]  (332,150) -- (333,220) ;

% Text Node
\draw (195,141) node [anchor=north west][inner sep=0.75pt]    {$1$};
% Text Node
\draw (195,221) node [anchor=north west][inner sep=0.75pt]    {$0$};
% Text Node
\draw (318,223) node [anchor=north west][inner sep=0.75pt]    {$R-1$};
% Text Node
\draw (395,222) node [anchor=north west][inner sep=0.75pt]    {$R$};
% Text Node
\draw (424,221) node [anchor=north west][inner sep=0.75pt]    {$t$};
% Text Node
\draw (210,114) node [anchor=north west][inner sep=0.75pt]    {$G_{R}$};

\end{tikzpicture}
\caption{}\label{GR2}
\end{figure}

Then $|G_R|\leqslant1$ and $\text{Lip}\; G_R\leqslant1$.  For any $f\in C_b(\mathbb{R}^{n}_+)$, $|f|_L\leqslant1$ we have
  \[
  \begin{split}\langle f, I\circ  m \rangle-\langle f,I\circ \mu^0\rangle
  &=\langle (fG_R)\circ I, m \rangle-\langle(f G_R)\circ I,\mu^0\rangle\\
  &\quad+\langle (f (1- G_R))\circ I, m \rangle-\langle (f (1- G_R))\circ I,\mu^0\rangle\\
  &\leqslant 2R\| m -\mu^0\|_{L,\mathbb{R}^{2n}}^*+2\kappa^2(\tfrac{1}{R-1}),
  \end{split}
  \]
  since on the ball $B_R$ we have $|I|\leqslant R^2$ and $\text{Lip} \;I\leqslant 2R$. Minimizing the r.h.s in $R\geqslant2$, we get \eqref{i-v-measure-q}.
   \end{proof}

   In the lemma below $\Om$ is the original probability space and $\Om^1$ is the segment $[0,1]$ with the Borel sigma-algebra $\cF^1$  and the Lebesgue
measure (i.e., another probability space, independent from $\Om$).

   \begin{lemma}\label{measure-itov}
   Let $\eps \in(0,1]$. For any ${\bar T}\ge0$ consider solution $v^\eps(\tau)$ %and $v^{st}(\tau)$
    on the interval $J= [{\bar T}, {\bar T} +1]$ and denote
   $$
  \delta:= \sup_{\tau \in J}  \|\cD(I( v^\eps(\tau))- I\circ \mu^0 \|_{L}^*.
   $$
   Then there exists a function\,\footnote{It depends on functions $\kappa^1$ and $\kappa^4$. The latter appears below in the lemma's proof and depends
   only on the constants in Assumption~\ref{a_8.1}.}
    $\kappa^5(\delta)$ and for each $\delta>0$ exists $r=r(\bar T, \delta, \eps) \in J$ and a r.v.  $\theta=\theta(\bar T, \delta, \eps) \in \T^n$,
    %measurable with respect to $\cF_r$ (see \eqref{OmOm}),
      such that
   %\be\label{10.i}
   $$
   \hve(r) := \Phi_\theta \tilde v^\eps(r), \qquad \cD \tilde v^\eps(r) = \cD v^\eps(r)
   $$
   %where $\theta=\theta(\delta, \eps) \in \T^n$ is a r.v., measurable with respect to $\cF_r$ (see \eqref{OmOm}), and
   satisfies
   \be\label{10.ii}
   \| \cD \hve(r) - \mu^0 \|^*_{L, \R^{2n}} \le \kappa^5(\delta).
   \ee
   The random vector $\tilde v^\eps(r)$ and the random variable $\theta \in \T^n$ are defined on the extended probability space
   $\Om\times\Om_1$ and are
   measurable with respect to sigma-algebra  $\cF_r\times \cF^1$,
   where $\cF_r$  is the sigma-algebra of the past    (see \eqref{OmOm}) and $\cF^1$ is the sigma-algebra on the probability space
   $\Om^1$.
      \end{lemma}

   \begin{proof} For  some $\tau\in J$ let us denote
   $
   I(\vep(\tau) )= \Iet$,  $ I(\vst(\tau) )= \Ist$ and abbreviate
   $$v^\eps(\tau) =: v^\eps, \;\;v^{st}(\tau) = v^{st}, \;\; \mu^\eps = \cD(v^\eps).$$
   %(we recall that $\cD v^{st}(\tau) \equiv \mu^0$).
   %Then $\vep(\tau) =  V_{\fet}(\Iet)$ and  $\vst(\tau) =  V_{\fst}(\Ist)$ (see \eqref{V_theta}).
    By Lemma~\ref{typical-time} and \eqref{i-st-small} for
   any $\gamma \in(0,1]$ there exists
   $
   \tau =   \tau(\bar T,\gamma, \eps)\in J
   $
   such that
   \be\label{10.16}
   \PP ( [ \Iet]<\gamma),  \; \PP([ \Ist] <\gamma) \le \kappa^1( \gamma).
   \ee
   On $\cP(\R_+^n)$ consider the Kantorovich distance
   $$
   \| \mu-\nu\|_K = \sup_{ Lip(f) \le1} \big( \lan f, \mu\ran - \lan f, \nu\ran\big) \le\infty,
   \quad \mu, \nu\in \cP(\R_+^n).
   $$
   Clearly $\| \mu-\nu\|_K \ge  \| \mu-\nu\|^*_{L, \R^n_+}$, and  due to estimates
    \eqref{M} and \eqref{x0} also
   \be\label{Kant}
   \| \cD \Iet - \cD \Ist\|_K\le \kappa^4(  \| \cD \Iet - \cD \Ist\|_{L,\R^n_+}^*) \le \kappa^4(\delta),
   \ee
   for some function $\kappa^4$.
   See   \cite[Section 11.4]{BK} and    \cite[Chapter 7]{Vil}. By the Kantorovich--Rubinstein
   theorem (see \cite{Dud, Vil, BK}) there exist r.v.'s \, $\tilI, \tilst$ such that
   $
   \cD \tilI = \cD I^\eps_\tau$, $\cD \tilst =\cD I_\tau^{st}
   $
   and
   \be\label{KR}
   \EE | \tilI - \tilst| = \| \cD I^\eps_\tau  - \cD I^{st}_\tau\|_K \le \kappa^4(\delta).
   \ee

 Now consider the operator of  projecting   vectors to their actions:
 $$
 I: \R^{2n} \to \R^n_+, \quad v \mapsto I(v).
 $$
 Then
$
I\circ \mu^\eps = \cD( I^\eps_\tau) = : \nu^\eps.
$
On the  extended probability space
$
\Om\times \Om^1 =\{(\om, \om_1)\}
$
   exists a r.v. $\tilde v^\eps(\om, \om_1)\in \R^{2n}$ such that
   $$
   \cD \tilde v^\eps = \mu^\eps, \quad I\big( \tilde v^\eps(\om, \om_1)\big) = \tilI(\om), \;\; a.s.
   $$
   Indeed, by the regular conditional probability theorem (see \cite[Sec.~10.2]{Dud}, \cite[Sec.~5.3.C]{brownianbook})
   the measure $\mu^\eps$ may be disintegrated as
   $$
   \mu^\eps(dv) = \mu_I(dv) \nu^\eps(dI),
   $$
   where $\mu_I$ is a measure on $\R^{2n}$ which is a measurable function of $I\in \R_+^n$ and is such that
   $
   \mu_{I'} (I^{-1}(I')) =1
   $
   for every $I' \in\R_+^n$. Then there exists a measurable mapping
   $
   \eta: \R_+^n \times \Om_1 \to \R^{2n}, \ (I, \om_1) \mapsto \eta^{\om_1} (I),
   $
   such that $ \mu_I \equiv \cD \eta^\cdot(I)$, see \cite[Theorem~1.2.28]{KS12}. Now it remains to set
   $
   \tilde v^\eps(\om, \om_1) = \eta^{\om_1}( \tilI(\om))$. Indeed,  for any function $f\in C_b(\R^{2n})$ we have
   $$
   \EE^{\om}\EE^{\om_1} f( \eta^{\om_1}( \tilI(\om)) =\int \nu^\eps (dI) \Big( \int f(v) \mu_I (dv)\Big).
   $$

   Similar, on $
\Om\times \Om^1
$
exists a r.v. $\tilde v^{st}$ such that
 $$
   \cD \tilde v^{st} = \mu^0, \quad I\big( \tilde v^{st}(\om, \om_1)\big) = \tilst(\om), \;\; a.s.
   $$

Let us denote $\tilde \f^\eps_\tau = \f( \tilde v^\eps)$ and $\tilde \f^{st}_\tau = \f(\tilde v^{st})$. Then  $\tilde v^\eps = V_{\tilde\f^\eps_\tau}( \tilI)$ and
  $\tilde v^{st}= V_{\tilde\f^{st}_\tau}( \tilst)$. Let us set
  $$
   \hve  = V_{\tilde\f^{st}_\tau} (\tilde I^\eps) = \Phi_{\tilde\f^{st}_\tau - \tilde\f^\eps_\tau} (\tilde v^\eps).
   $$
   Then $\f(\hve) = \f(\tilde v^{st})$, so
    \be\label{10.18}
   \| \hve - \tilde v^{st} \|  \le \tfrac1{2\sqrt\gamma}  | \tilI - \tilst |  \quad\text{if}\quad [\tilI], \; [\tilst] \ge \gamma.
   \ee

   Now we estimate the distance between $\cD(\hve)$ and $\cD\tilde v^{st}= \mu^0$. To do that we take any function $f\in C_b(\R^{2n})$,
   $|f|_L \le1$, and consider
   $$
    X:= \EE \big( f(\hat v^\eps ) - f(\tilde v^{st}).
     $$
   Introducing the events
   $$
   Q^1_\gamma =\{ [I(\hat v^\eps )]<\gamma\} ,\;  Q^2_\gamma =  \{ [I(\hat v^{st} )]<\gamma\},\quad Q_\gamma = Q^1_\gamma\cup Q^2_\gamma
   %Q_\gamma^c = \Omega\setminus Q_\gamma,
   $$
   we write $X$ as
     $$
   X=X_1+ X_2, \quad X_1 =\EE  \big( f(\hat v^\eps ) - f(v^{st}) \big)\b1_{Q_\gamma}, \;\;
    X_2 =\EE  \big( f(\hat v^\eps ) - f(v^{st})\big) \b1_{Q^c_\gamma}.
   $$
   Since in view of \eqref{10.16}
    \be\label{77}
   \PP(Q^1_\gamma),  \;  \PP(Q^2_\gamma) \le 2 \kappa^1(\gamma)
   \ee
   and $|f|\le1$,    then $X_1\le 4\kappa^1(\gamma)$. As Lip$f \le1$, then using \eqref{10.18} and \eqref{KR} we get that
   $$
   X_2 \le  \EE  \big |\hve
    -\tilde v^{st} \big| \b1_{Q^{c}_\gamma}
     \le\tfrac1{2\sqrt\gamma}  \EE  | \tilI - \tilst| \le  \tfrac{\kappa^4(\delta)}{2 \sqrt\gamma}.
   $$

   Re-denoting $\hve$ to $\hve(\tau)$ and $\tilde v^{st}$ to $\tilde v^{st}(\tau)$ we get from the estimates on $X_1, X_2$ that
      \[
   \mathbf{E} \big(f({\hve(\tau)})-f({\vst(\tau)})\big)
  % \mathbf{E}(f((\hve(\tau)) -  f(\vst(\tau))
   \leqslant 8\kappa^1(\gamma)+\tfrac{1}{2\sqrt{\gamma}}\, \kappa^4(\delta).
   \]
   Minimizing in $\gamma\in (0,1]$ we achieve that the r.h.s. is less then $\kappa^5(\delta)$ for some $\gamma(\bar T,\delta, \eps)$,
   that is for some $\tau(\gamma(\bar T,\delta, \eps), \eps) =: r(\bar T, \delta, \eps)$. Since $f$ is any continuous function with $|f|_L \le1$, then
   $$
   \| \cD \hve( r(\bar T, \delta, \eps)) - \cD \vst( r(\bar T, \delta, \eps)) \|^*_{L, \R^{2n}}
    \le  \kappa^5(\delta).
   $$
   This relation and the formula for $\hve(\tau)$    prove the lemma.
      \end{proof}

      Now we state and prove a key lemma for the proof of the theorem. Below
      functions $\kappa^3$ and $\kappa^5$  are as in Lemmas~\ref{i-v-measure} and \ref{measure-itov}, and we
       recall notation \eqref{v^0_T}.
\begin{lemma}\label{random-initial-average}
\begin{enumerate}
\item For any $T>0$ and $\delta>0$  there exists $\eps_1 =\eps_1(\delta,T)>0$ such that if $\eps\le \eps_1$, then
\be\label{gr}
\sup_{\tau\in[0,T]}\|I\circ\cD(v^{\eps}({\bar T}+\tau)) - I\circ\cD(v_{{\bar T}}^{0}(\tau))\|_{L,\mathbb{R}_+^n}^*\leqslant\delta/2 \quad \forall\, {\bar T}\geqslant0.
\ee
%where $a^{0}(\tau)$ solves the effective equation \eqref{effective-equation} with initial condition $\cD(a^{0}(0))=\cD(a^{\eps}(T';v_0))$.
\item
For any $\delta>0$, choose a function $T^*=T^*(\delta)\ge0$ such that  $\kappa^3(g_{M^*}(T,2))\leqslant \delta/2$ for any $T\geqslant T^*(\delta)$.
Then there exists $\eps_2 =\eps_2(\delta) \in (0,1]$ with the following property: \\
assume that  a non-random   $T' =T'(\delta, \eps)\ge0$ is such that for every $0<\eps \le \eps_2$ it holds that
\begin{equation}\label{key-9-lemma1}
\|I\circ \cD(v^{\eps}(T'))-I\circ\mu^0\|_{L}^*\leqslant \delta
\end{equation}
and
\be\label{10.100}
\| \cD \hve  - \mu^0\|^*_{L, \R^{2n}} \le \kappa^5(\delta)
\ee
for  $\hve= \hve(T') = \Phi_\theta(\vep(T'))$, where  $\theta$ is some r.v.  (depending on $\delta$ and $\eps$), measurable with  respect to
$\cF_{T'}$.
Then for $0<\eps\le \eps_2$ we have
\begin{equation}\label{key-9-lemma2}
\sup_{\theta\in[0,1]}\|I\circ \cD(v^{\eps}(T'+T^*+\theta))-I\circ \mu^0\|_{L}^*\leqslant \delta,
\end{equation}
and
\begin{equation}\label{key-9-lemma3}
\sup_{\tau\in[T',T'+T^*+1]}\|I\circ \cD(v^{\eps}(\tau))-I\circ \mu^0\|_{L}^*\leqslant \frac{\delta}{2}+\kappa^3\big(\max_{\theta\in[0,T^*+1]}g_{M^*}(\theta,\kappa^5(\delta) )\big).
\end{equation}
\end{enumerate}
\end{lemma}

\begin{proof}
For a measure $\nu\in \cP(\mathbb{R}^{2n})$ we
denote by
$v^\eps(\tau;\nu)$ a weak solution of eq.~\eqref{v-equation1} such that $\cD (v^\eps(0)) =\nu$, and define $v^0(\tau;\nu)$
similarly. Since eq.~\eqref{v-equation1} defines a Markov process in $\mathbb{R}^{2n}$  (e.g. see \cite[Section~5.4.C]{brownianbook} and
\cite[Section~3.3]{ khasminskii}), then
$$
I\circ\cD (v^\eps(\tau;\nu)) = \int_{\mathbb{R}^{2n}} I\circ\cD (v^\eps(\tau;v )) \, \nu(dv),
$$
and a similar relation holds for $I\circ\cD (v^0(\tau;\nu))$.

(1) Denote $\nu^\eps = \cD (v^\eps({\bar T}))$. Then
\be\label{y1}
\cD (v^\eps({\bar T}+\tau)) = \cD (v^\eps(\tau; \nu^\eps)), \quad  \cD (v^0_{{\bar T}}(\tau) )= \cD (v^0(\tau; \nu^\eps)).
\ee
By  \eqref{kappa2}, for any $\delta>0$  there exists $K_\delta>0$ such that for each $\eps$,
 $\nu^\eps(\mathbb{R}^{2n} \setminus \bar B_{K_\delta} )\le \delta/8$, where
$\bar B_{K_\delta}:=\bar B_{K_\delta}(\R^{2n})$.
So
$$
\nu^\eps = A^\eps \nu^\eps_\delta +  \bar A^\eps \bar\nu^\eps_\delta , \quad  A^\eps = \nu^\eps( \bar B_{K_\delta}),\;
\bar A^\eps= \nu^\eps(\R^{2n}\setminus \bar B_{K_\delta}),
%A^\eps +\bar A^\eps =1,
$$
where  $ \nu^\eps_\delta $ and $ \bar\nu^\eps_\delta $ are the conditional probabilities $ \nu^\eps(\cdot \mid  \bar B_{K_\delta})$ and
$ \nu^\eps(\cdot \mid\R^{2n}\setminus   \bar B_{K_\delta})$.
Accordingly,
%\be\label{decompp}
$$
\cD (v^\kappa (\tau; \nu^\eps) )= A^\eps \cD (v^\kappa (\tau; \nu^\eps_\delta) )+ \bar A^\eps \cD (v^\kappa (\tau; \bar\nu^\eps_\delta)),
$$
where $\kappa=\eps$ or $\kappa=0$. Therefore,
\[\begin{split}&
\| I\circ\cD (v^\eps (\tau; \nu^\eps) )- I\circ\cD( v^0 (\tau; \nu^\eps)) \|_L^*\\
&\le
A^\eps \| I\circ \cD (v^\eps (\tau; \nu^\eps_\delta) )- I\circ \cD (v^0 (\tau; \nu^\eps_\delta) )\|_L^* +
\bar A^\eps \| I\circ\cD (v^\eps (\tau; \bar \nu^\eps_\delta)) - I\circ\cD (v^0 (\tau;\bar \nu^\eps_\delta) )\|_{L}^*.
\end{split}\]
The second term on the r.h.s  obviously is bounded by
$2\bar A^\eps\leqslant\frac{\delta}{4}$. While by  Proposition~\ref{p_ran_ini},
 there exists $\eps_1>0$, depending only on   $K_\delta$ and $T$,  such that for $0\le \tau\le T$  and  $\eps\in(0,\eps_1]$
 the first term in the r.h.s.   is $\leqslant\frac{\delta}{4}$.
 Due to \eqref{y1} this proves the first  assertion.
  \smallskip

(2) Let us choose ${\eps}_2(\delta):={\eps}_1(T^*(\delta)+1,\delta)$.
We have
\begin{equation}\label{i-difference1}
\begin{split}&\sup_{\tau\in[0,1]}\|I\circ \cD(v^{\eps}(T'+T^*+\tau))-I\circ \mu^0\|_{L}^*\\
&\leqslant\sup_{\tau\in[0,1]}\|I\circ \cD(v^{\eps}(T'+T^*+\tau)-I\circ \cD(v^0_{T'}(T^*+\tau))\|_{L}^*\\
&+\sup_{\tau\in[0,1]}\|I\circ \cD(v^0_{T'}(T^*+\tau))-I\circ \mu^0\|_{L}^*.
\end{split}
\end{equation}
By  \eqref{gr} and the choice of $\eps_2$,
 the first term in the r.h.s is less than $\frac{\delta}2$.
Let us examine the second one.
  By Proposition~\ref{p_inirotate},
% Remark \ref{remark-phase},
\[\cD\big(I ( v^0(T^*+\tau; {\hve}))\big) =\cD\big(I( v_{T'}^0(T^*+\tau))\big),\; \;\; \forall\tau\in[0,1].\]
Thus the second term in the r.h.s of \eqref{i-difference1} equals $\sup_{\tau\in[0,1]}\|I\circ \cD(v^0(T^*+\tau; {\hve}))-I\circ\mu^0\|_{L}^*$.  Since
$
\| \hve\| = \| v^\eps(T')\|,
$
then by \eqref{key-9-lemma1}  and \eqref{mixing-quantify}
\[\|\cD(v^0(T^*+\tau; {\hve}))-\mu^0\|_{L,\mathbb{R}^{2n}}^*\leqslant g_{M^*}(T^*+\tau,\kappa^5(\delta)). \]
So in view of Lemma~\ref{i-v-measure} the second term in the r.h.s of \eqref{i-difference1} is bounded by $\sup_{\tau\in[0,1]}\kappa^3(g_{M^*}(T^*+\theta,\kappa^5(\delta)))$, which is
$\leqslant \delta/2$ due to the definition of
$T^*(\delta)$.  This proves \eqref{key-9-lemma2}.

Similarly,
\[\begin{split}\sup_{\tau\in[T',T'+T^*+1]} &\|I\circ\cD(v^{\eps}(\tau)) -I\circ\mu^0\|_L^*\\
&\leqslant\sup_{\theta\in[0,T^*+1]}\|I\circ\cD(v^{\eps}(T'+\theta))-I\circ\cD(v^0(\theta;{\hve}))\|_L^*\\
&+\sup_{\theta\in[0,T^*+1]}\|I\circ\cD(v^{0}(\theta;{\hve}))-I\circ\mu^0\|_L^*
\end{split}\]
By \eqref{gr} and the definition of $\eps_2$,
the first term in the r.h.s is less than $\delta/2$, while by \eqref{mixing-quantify} and Lemma~\eqref{i-v-measure}, the second term is bounded by
$\kappa^3(\lambda)$, where $\lambda = \max_{\theta\in[0,T^*+1]}g_{M^*}(\theta,\kappa^5(\delta))$.
Thus we proved \eqref{key-9-lemma3}.
\end{proof}

Now we are ready to prove the theorem.
Let us fix arbitrary $\delta>0$ and take some $0<\delta_1\le \delta/4$. Below in the proof the functions $\eps_1$, $\eps_2$ and $T^*$ are as in
Lemma~\ref{random-initial-average}. We will abbreviate  $T^*(\delta_1) =: T^*$, $\eps_2(\delta_1) =: \eps_2$ and will always assume that
$$
0<\eps \le \eps_2.
$$

i) By the definition of $T^*$, \eqref{mixing-quantify}  and \eqref{M},
%\eqref{M} and \eqref{mixing-quantify} with $d=2$, there exists $T^*=T^*(\delta)$ such that
\be\label{84}
\begin{split}
\| \cD\big( v^0_{\bar T} (\tau)) - \mu^0\|_{L,\mathbb{R}^{2n}}^* &\le  g_{M^*}(\tau,2)\quad \forall \,\bar T\ge0, \\
\kappa^3(g_{M^*}(\tau,2))&\le\delta_1/2 \quad \forall\, \tau\ge T^*.
\end{split}
\ee

ii)  By \eqref{gr}
\be\label{85}
\sup_{0 \le \tau \le T^*+1}
\| I\circ\cD\big( v^\eps (\tau)) -   I\circ\cD\big( v^0 (\tau;v_0) \big) \|_L^* \le     \tfrac{\delta_1}2.
\ee
From  \eqref{84},  \eqref{85}  and Lemma~\ref{i-v-measure} we have
\be\label{855}
\sup_{\tau\in[0,1]}\| I\circ\cD\big( v^\eps (T^*+\tau)) -  I\circ\mu^0 \|_L^* <  \delta_1.
\ee

iii)  By  \eqref{85} and  Lemma~\ref{measure-itov} there exists $T'_1 = T'_1(T^*, \delta_1, \eps) \in [T^*, T^*+1]$ such that apart from
\eqref{855}, there exists a r.v. $\theta_1 \in \T^n$, measurable with respect to $\cF_{T_1'}$, and such that
$
\hve_1(T'_1) = \Phi_{\theta_1} v^\eps(T'_1)
$
satisfies
\be\label{990}
\| \cD \hve_1(T'_1) - \mu^0\|^*_{L, \R^{2n}} \le \kappa^5(\delta_1),
\ee
Considering $v^\eps(T'_1+\tau)$, $\tau\ge0$, we get that \eqref{855} holds with $T^*$ replaced by $T'_1 +T^*$, and an analogy of \eqref{990} holds for $T'_1$
replaced with  some $T_2'\in  [T'_1+T^*, T'_1+ T^*+1]$. Iterating this argument we construct a sequence
 $T_N'$, $N=1,\dots$, such that $T'_{N+1} \in [ T'_N +T_*, T'_N +T_*+1 ]$,
$$
\| I\circ\cD\big( v^\eps (T_N')) - I\circ \mu^0 \|_L^* \le  \delta_1 \quad \forall\, N,
$$
and
$$
\| \cD \hve_N(T'_N) - \mu^0\|^*_{L, \R^{2n}} \le \kappa^5(\delta_1) \quad \forall\, N,
$$
for $\hve_N(T'_N)  = \Phi_{\theta_N} v^\eps(T'_N)$ with a suitable $\cF_{T'_N}$-measurable $\theta_N$. So by \eqref{key-9-lemma3}
\be\label{857}
\sup_{\tau \in [T'_N, T'_{N+1}]}
\| I\circ\cD\big( v^\eps (\tau)) - I\circ\mu^0 \|_L^* \le \frac{\delta_1}{2}+\kappa^3\Big(\max_{\theta\in[0,T^*+1]}g_{M^*}\big(\theta,\kappa^5(\delta_1)\big)\Big),
\ee
for every $N$.

iv) Finally,  by \eqref{85} if $\tau\le T'_1$ and by \eqref{84} with $\bar T=0$ jointly with \eqref{857} if $\tau\ge T_1'$,  we have that
$$
\| I\circ\cD\big( v^\eps (\tau)\big) -   I\circ\cD\big( v^0 (\tau;v_0) \big) \|_L^* \le {\delta_1}+\kappa^3\Big(\max_{\theta\in[0,T^*+1]}g_{M^*}\big(\theta,\kappa^5(\delta_1)\big)\Big),\qquad \forall\, \tau\ge0.
$$
By the assumption, imposed in $(6)$ of Assumption \ref{a_8.1} on function $g_{M}$,   $g_{M} (t, d)$ is uniformly continuous in $d$ and
vanishes at $d=0$. Recall that $\kappa^3$ and $\kappa^5$ are both monotonically increasing continuous functions,
vanishing at $0$.  Therefore we have that  there exists $\delta^*>0$, which we may assume to be $\leqslant \delta/2$, such that if $\delta_1\le \delta^*$, then $\kappa^3(g_{M^*}(\theta,\kappa^5(\delta_1)))\leqslant \delta/2$ for all
$ \theta\geqslant0$ . Then by the estimate above,
$$
\| I\circ\cD\big( v^\eps (\tau)) -  I\circ \cD\big( v^0 (\tau;v_0) \big) \|_L^* \le
 \delta,\;\;\; \forall \tau\geqslant0  \quad\text{if} \quad \eps \le %\eps_*(\delta) :=
    \eps_2(\delta^*\big(\delta)\big)>0,
 %\min(\eps_1, \eps_2)(\delta^*, T^*),
$$
for every positive $\delta$. This proves the theorem's assertion.
\end{proof}

We end this section with a sufficient condition  for the validity of (4)-(6) in Assumption \ref{a_8.1}.
\begin{proposition}\label{sufficient1}
Assume that (2)-(3) of Assumption \ref{assumption-v-3} hold true and  the field $P(v)$ in \eqref{v-equation1}  is coercive. That is,
 there exist $\alpha_1>0$ and $\alpha_2\geqslant0$ such that
\begin{equation}\label{dissipative1}( P(v),v) \leqslant -\alpha_1\|v\|+\alpha_2,\; \forall v\in\mathbb{R}^{2n},
\end{equation}
where  $(v,w)=\sum_{j=1}^n \bv_j\cdot\mathbf{w}_j$ is the
 inner product on $\mathbb{R}^{2n}$. Then (4)-(6) in Assumption~\ref{a_8.1} hold true.
\end{proposition}
\begin{proof}
For the
drift term in \eqref{v-equation1}  $b(v)=({\eps}^{-1}W_k(I)\bv_k^\bot+{\mathbf{P}}_k(v), k=1,\dots,n)$, by \eqref{dissipative1}, we have
$$
%\begin{equation}\label{drift-bound}
(b(v),v)=(P(v),v)\leqslant -\alpha_1\|v\|+\alpha_2,\; \forall v\in\mathbb{R}^{2n}.
$$
For the drift term $R$  in  effective equation \eqref{effective1},
by \eqref{dissipative1} and the definition or $R= \lan P\ran$  we have
\[(R(v),v)=\sum_{k=1}^n\int_{\mathbb{T}^n}{\mathbf{P}}_k(\Phi_\theta v)\cdot \Phi_\theta^k\bv_kd\theta=\int_{\mathbb{T}^n}(P(\Phi_\theta v),\Phi_\theta v)d\theta\leqslant -\alpha_1\|v\|+\alpha_2,\;\forall v\in\mathbb{R}^{2n}.\]
By the definition of the diffusion matrix of equation
\eqref{effective1} we know that the uniform ellipticity condition as
 in (2) of Assumption \ref{assumption-v-3} holds for the effective equation \eqref{effective1}.
Then the assertion of proposition directly follows from \cite[Proposition 9.3]{GHKu}.
\end{proof}
\begin{remark}The assumption in Proposition \ref{sufficient1} also ensures  the mixing in equation
 \eqref{v-equation1} for each ${\eps}\in(0,1]$, see in \cite{GHKu}.
 In this case, for   corresponding stationary measures $\mu^{\eps}$,  the measures
 $I\circ \mu^{\eps}$ converge weakly to $I\circ \mu^0$ as ${\eps}\to0$.
\end{remark}

\subsection{On the proof of Theorem \ref{t_mixAP}} \label{ss_10.1}
Due to Proposition \ref{p_suff_cond} the laws of solutions $I^0(\tau)$ for equation \eqref{averaged-eq-1} obey estimates \eqref{mixing-quantify}
(where $v(\tau)$ is replaced by $I^0(\tau)$). Now, denoting $v^\eps(\tau) := (I^\eps(\tau), \f^\eps(\tau))$ we may repeat for
$
  I(v^\eps(\tau))= I^\eps(\tau)
$
the proof of Theorem~\ref{thm-uniform}, given above in Section~\ref{s_9}, with $I\circ v^0(\tau)$ replaced by $I^0(\tau)$. In fact, a proof
of  Theorem~\ref{t_mixAP} is simpler than that of Theorem~\ref{thm-uniform} since in the former theorem the mapping
$
v^\eps \mapsto (I^\eps, \f^\eps)
$
is a trivial isomorphism, while in the setting of the latter theorem this is the non-linear action-angle mapping \eqref{ac_an}, which is singular
when some $\bv_j^\eps$ vanishes. Accordingly, to prove an analogy of the key Lemma~\ref{random-initial-average} for solutions of
equation \eqref{ori_scaled} the technical argument, contained in Lemmas~\ref{typical-time}, \ref{i-v-measure} and
\ref{measure-itov} becomes redundant. We skip details of an exact realisation of this sketch.

\subsection{Damped/driven Hamiltonian systems.} \label{ss_10.2}
As an example let us consider system \eqref{v-equation1}, where the drift $P(v)$ is a damped Hamiltonian field and the diffusion
is a diagonal  additive random force:
\be\label{exampl}
d\bv_k=\ep^{-1}W_k(I)\bv_k^{\bot}d\tau -\nu_k \bv_k d\tau +J \nabla_{\bv_k} h(v) d\tau + \gamma_k d\bbeta_k(\tau), \quad
k=1,\dots,n,
\ee
$v(0) = v_0$. Here Hamiltonian $h(v)$ is $C^2$-smooth, $J=\left(\begin{array}{cc} 0&-1 \\
1&0\end{array}\right)$,
 $\nu_k$'s are nonnegative real numbers and $\gamma_k$'s are nonzero real numbers. Results in
 Sections~\ref{s_Birk_int}-\ref{s_9}  apply to
this system if Assumption~\ref{a_8.1} holds. There items (1)-(3) are easy to meet, while items (4)-(6) hold e.g. if  (1)-(3) are
satisfied, all  $\nu_k$'s are positive and Hamiltonian $h(v)$ commutes with $\| v\|^2$. Indeed, then the scalar product
$
(J\nabla h(v), v)
$
vanishes, so the drift in the equation is coercive, and items (4)-(6) hold in view of Proposition~\ref{sufficient1}.

To write down the effective equation we note that its  dispersion is the same as in eq.~\eqref{exampl}
 since now matrix $X(v)$ as in Section~\ref{construct-effective} is diag$\{\gamma_k^2\}$.
 The drift $R(v)$ is an averaging of the drift $-\{\nu_k \bv_k\}  +J \nabla h(v)$. It is easy to see that the averaging does not change the
first term. The averaging  of the second one is
$
J \nabla  \langle h \rangle (v),
$
where $ \langle h \rangle$ is the averaged   Hamiltonian,
$\
 \langle h \rangle (v) = \int_{\T^n} h(\Phi_{\theta} v) d\theta
$
(see Proposition~3.5 in \cite{GHKu}).  So the effective equation reads
\be\label{eff_eqq}
d\bv_k=-\nu_k \bv_k d\tau +J \nabla_{\bv_k} \langle h \rangle(v) d\tau + \gamma_k d\bbeta_k(\tau), \quad
k=1,\dots,n.
\ee
Obviously the function $ \langle h \rangle (v)$ depends only on the actions i.e. $ \langle h \rangle= h^0(I_1, \dots, I_n)$,
where  $h^0$ is a continuous function, and
$
h^0(, \dots, I_j, \dots) = h^0(, \dots, -I_j, \dots)
$
for each $j$. Since $h$ is $C^2$-smooth, then the averaged Hamiltonian
 $\langle h \rangle$ is $C^2$-smooth in $v$ as well. From here and Whitney's theorem (see in \cite{Whit} Theorem~1 with $s=1$
 and the remark, concluding that paper) we derive that the function $h^0(I)$ is $C^1$-smooth. So in effective equation
 \eqref{eff_eqq}
 \be\label{11.1}
 J \nabla_{\bv_k} \langle h \rangle(v) = (\p h^0(I)/ \p I_k) \bv_k^\perp, \quad h^0 \in C^1.
 \ee
 In particular,  applying Ito's formula to  an action $I_k=\tfrac12 | \bv_k|^2$, where $v(\tau)$ solves \eqref{eff_eqq}, we get that
  \be\label{11.2}
 \frac{d}{d\tau} \EE \tfrac12 | \bv_k|^2(\tau) = -2\nu_k  \EE \tfrac12  | \bv_k|^2(\tau) +{\gamma_k^2}.
 \ee
 So the stationary measure $\mu^0$  for \eqref{eff_eqq} is such that
 $
 \int \tfrac12 | \bv_k|^2 \mu^0(dv) = {\gamma_k^2} /(2\nu_k )
 $
 for each $k$, and  we get from Theorem~\ref{thm-uniform} that for any $v_0$ solution $v^\eps(\tau)$ of \eqref{exampl} satisfies
 $$
 \lim_{\tau\to\infty, \eps\to0}   \EE \tfrac12|\bv^\eps_k|^2(\tau) = \frac{{\gamma_k^2}}{2 \nu_k}, \quad \forall\, k.
 $$

 If all $\nu_k$'s vanish (but still $\gamma_k$'s are nonzero and $h(v)$ commutes with $\|v\|^2$), then Assumption~\ref{a_8.1}  does not hold since \eqref{v-upper-2} fails. Nonetheless,
  Assumption~\ref{assumption-v-3} remains valid (which can be readily verified by the same argument as in the proof of Proposition \ref{l_suff_assum_i4}), and the effective equation has the form \eqref{eff_eqq}, \eqref{11.1} with $\nu_k \equiv 0$. So we
  derive from \eqref{11.2} with $\nu_k=0$ that on any finite time-interval $[0,T]$ the averaged actions of solutions $v^\eps(\tau)$ for
  \eqref{exampl} approximately have linear growth with $\tau$:
  $$
  \EE \tfrac12|\bv^\eps_k|^2(\tau)  = \tfrac12 |\bv_{0k}|^2 + {\gamma_k^2}\, \tau + o_{\epsilon\to0}(1), \quad 0\le \tau\le T,\;\; \;\forall\, k,
  $$
 where the term $o_{\epsilon\to0}(1)\to0$ as $\epsilon\to0$.

    \appendix

  \section{Proof of % \eqref{estim_ogo1}
  Lemma~\ref{l_key_ham}} \label{app_proof_lem_not0}

We fix ${\eps}\in(0,1]$ and do not indicate the dependence on it. Relation \eqref{estim_ogo1} is already established. A proof of \eqref{estim_ogo2}
goes in 3 steps.

{\bf Step 1}: Constructing for a fixed $k$ and
any $\delta\in(0,1]$ an It\^o process $\bar\bv_k^\delta(\tau)$, $\tau\in[0,T]$, such that
$| \bar\bv_k^\delta| \equiv | \bar\bv_k| $, and
 if $\|\bar\bv_k^\delta\|\geqslant\delta$ then no $\frac{1}{{\eps}}$-term explicitly appear in the drift term.

  Denote by ${U}={U}(\zeta_1,\zeta_2):(\mathbb R^2\setminus\{0\})\times(\mathbb R^2\setminus\{0\})
  \mapsto SO(2)$ the unique  rotation of  $\mathbb R^2$
  that maps $\frac{\zeta_2}{|\zeta_2|}$ to  $\frac{\zeta_1}{|\zeta_1|}$. Then
   ${U}(\zeta_2,\zeta_1)= ({U}(\zeta_1,\zeta_2))^{-1}= ({U}(\zeta_1,\zeta_2))^t$.

Let $v(\tau)=(\bv_k(\tau),k=1,\dots,n)$ be a solution of  equation \eqref{v-equation1}.
We introduce the vector-functions
$$\bar {\mathbf{P}}_k(\bar \bv_k,v)=U(\bar \bv_k,\bv_k){\mathbf{P}}_k(v),\quad \bar B_{kj}(\bar\bv_k,v)=U(\bar \bv_k,\bv_k)B_{kj}(v),$$
where $k=1,\dots,n,\; j=1,\dots,n_1$. We fix some $k$ and
consider the following stochastic equation for $\bar \bv_k(\tau)\in\mathbb{R}^2$:
\begin{equation}\label{modi_rotated_pert1}
d\bar\bv_k=\bar {\mathbf{P}}_k(\bar\bv_k,v(\tau))d\tau+\sum_{j=1}^{n_1}\bar B_{kj}(\bar\bv_k,v(\tau))d\bbeta_j(\tau).
\end{equation}
Its coefficients are well defined for all non-zero $\bv_k$ and $\bar \bv_k$. The equation \eqref{modi_rotated_pert1}, given some
 initial data, has a unique solution as long as
$\|\bv_k\|,\|\bar\bv_k\|\geqslant\delta$ for any fixed $\delta>0$.

For an arbitrary $\delta\in(0,\frac12)$ we define the stopping times $\tau_j^\pm$ as follows:
$\tau_0^+=0$,
\[ \begin{split}
\tau_j^-&=\inf\left\{s\geq \tau_{j-1}^+\,:\, \min\limits_{1\leq k\leq n}\|\bv_k(s)\|\leq \delta\
\hbox{or }\|v(s)\|\geq\delta^{-1}\right\}, \quad j\geq 1,\\
\tau_j^+&=\inf\left\{s\geq \tau_{j}^-\,:\, \min\limits_{1\leq k\leq n}\|\bv_k(s)\|\geq 2\delta\
\hbox{ and } \|v(s)\|\leq(2\delta)^{-1}\right\}, \quad j\geq 1.
\end{split}
\]
Note that  $\tau_0^+\leq \tau_1^-$ and $\tau_j^-<\tau_j^+<\tau_{j+1}^-$ for $j\geq 1$.  See again Fig.~1, where now the line is the graph of the function
$ \|\mathbf{v}_k(\tau) \|$.

Since on each interval $\Lambda_j$ the norm of  solution $v(\tau)$ of \eqref{v-equation1} is bonded by $\delta^{-1}$, then
$\Lambda_j$ cannot be too short. So the sequence  $\tau_j^\pm$ stabilizes  at  $T$ after a finite random
number of steps.

Now we construct a continuous process $\bar{\bv}_k^\delta(\tau)$, $\tau\in[0,T]$. We set $\bar{\bv}_k^\delta(\tau_0^+)=\bv_k(\tau_0^+)$.
 For any $j\geqslant 0$ we define   $\bar{\bv}_k^\delta$ on   the segment
$\Lambda_{j}:=[\tau_{j}^+,\tau_{j+1}^-]$  as a solution of equation \eqref{modi_rotated_pert1},
 while on the complementary segments $\Delta_r=[\tau_r^-,\tau_r^+]$ we set
 \begin{equation}\label{bad_segm}
 \bar{\bv}^\delta_k(s)=
U\big(\bar\bv_k(\tau_r^-), \bv_k(\tau_r^-)\big)\bv_k(s).
 \end{equation}

 \begin{lemma}\label{lem_norm-eq}
   If  $\| \bar\bv_k^\delta(\tau_j^+)\|=\|\bv_k(\tau_j^+)\|$,   then $\|\bar\bv^\delta_k(s)\|=\|\bv_k(s)\|$
   for all $s\in\Lambda_{j}$, a.s.
 \end{lemma}
 \begin{proof}
 Denote $I_k^\delta=\frac12\| \bar\bv^\delta_k\|^2$.
   By  It\^o's formula, on the segment $\Lambda_{j-1}$,
   $$
   d I_k^\delta=\big(\bar\bv^\delta_k,  \bar  {\mathbf{P}}_k(\bar\bv^\delta_k,v(\tau))\big)d\tau +\sum_{i=1}^{n_1}\Big(\frac{1}{2}\|\bar B_{ki}(\bar\bv_k^\delta,v(\tau))\|^2_{HS}d\tau
   +\big(\bar\bv^\delta_k, \bar B_{ki}(\bar\bv_k^\delta,v)\big)d\bbeta_{i}(\tau)\Big),
   $$
and $I_k=\frac{1}{2}\|\bv_k\|^2$ satisfies
 $$
   d I_k  =\big(\bv_k,    {\mathbf{P}}_k(v)\big)d\tau +\sum_{i=1}^{n_1}\Big(\frac{1}{2}\| B_{ki}(v(\tau))\|^2_{HS}d\tau
   +\big(\bv_k,  B_{ki}(v)\big)d\bbeta_{i}(\tau)\Big).
   $$
By construction we have the following relations for the drift and diffusion term of these two equations:
\[\big(\bar\bv_k^\delta,\bar {\mathbf{P}}_k(\bar\bv_k^\delta,v)\big)+\frac{1}{2}\sum_{i=1}^{n_1}\|\bar B_{ki}(\bar \bv_k^\delta, v)\|_{HS}^2=\frac{\|\bar\bv_k^\delta\|}{\|\bv_k\|}\big(\bv_k,{\mathbf{P}}_k(v)\big)+\frac{1}{2}\sum_{i=1}^{n_1}\|B_{ki}(v)\|_{HS}^2,\]
\[\big(\bar\bv_k^\delta,\bar B_{ki}(\bar\bv_k^\delta,v)\big)=\frac{\|\bar\bv_k^\delta\|}{\|\bv_k\|}\big(\bv_k,B_{ki}(v)\big).\]
For the squared difference $(I_k-I_k^\delta)^2$ we have,
\begin{equation}\label{squa_ig}
  \begin{split}
d(I_k-I_k^\delta)^2&=\Big(2(I_k-I_k^\delta)\frac{|\bv_k|-|\bar\bv^\delta_k|}{|\bv_k|}
 \big(\bv_k, {\mathbf{P}}_k(\bv)\big)\\
 &+\frac{(|\bv_k|-|\bar\bv^\delta_k|)^2}{|\bv_k|^2}\sum\limits_{i=1}^{n_1}
 \|(\bv_k,B_{ki}(v))\|^2\Big)d\tau+d\mathcal{M}^\delta_\delta
 \end{split}
 \end{equation}
with $\mathcal{M}^\delta_\tau$ being a square integrable martingale. Since
$I_k-I_k^\delta=\frac12\big(\|\bv_k\|-\|\bar\bv^\delta_k\|\big)
\big(\|\bv_k\|+\|\bar\bv^\delta_k\|\big)$, letting
$\mathcal{J}_k^\delta(\tau)=(I_k-I^\delta_k)^2((\tau\vee \tau_j^+)\wedge\tau_{j+1}^-)$ and taking the expectation in
\eqref{squa_ig} we have
$$
\mathbf{E}\mathcal{J}_k^\delta(\tau)\leq \mathbf{E}J_k^\delta(0)+c(\delta) \int_0^\tau\mathbf{E}\mathcal{J}_k^\delta(s)\,ds.
$$
Since $J^\delta_k(\tau_{j}^+)=0$, then by Gronwall's lemma,   $\mathcal{J}_k^\delta(\tau)=0$ for
$\tau\in\Lambda_{j}$. The assertion of the lemma is proved.
 \end{proof}
By Lemma~\ref{lem_norm-eq} we have $\|\bv_k(s)\|=\|\bar\bv^\delta_k(s)\|$ for all $s\in\Lambda_{0}$.
By \eqref{bad_segm}  $\|\bv_k(s)\|=\|\bar\bv^\delta_k(s)\|$ for   $s\in \Delta_j$.
Iterating this procedure we conclude that  $\|\bv_k(s)\|=\|\bar\bv^\delta_k(s)\|$
on the whole interval $[0,T]$.

We define
$$
\widehat {\mathbf{P}}_k(\bar\bv^\delta_k, v,s)=\left\{
\begin{array}{ll}
  \bar  {\mathbf{P}}_k(\bar\bv^\delta_k,v), & \hbox{if } s\in\mathop{\cup}\limits_j \Lambda_{j},\\
  {U}_j\big[\frac1{\eps} W_k(I)\bv_k^\bot+{\mathbf{P}}_k(v)
  \big],
   & \hbox{if } s\in\mathop{\cup}\limits_j \Delta_j,
\end{array}
\right.
$$
where ${U}_j={U}(\bar\bv_k(\tau^-_j), (\bv_k(\tau_j^-)))$, and
$$
\widehat{B}_{ki}(\bar\bv^\delta_k, \bv,s)=\left\{
\begin{array}{ll}
  \bar {B}_{ki}(\bar\bv^\delta_k,v), & \hbox{if } s\in\mathop{\cup}\limits_j \Lambda_{j},\\
  {U}_jB_{ki}(v(s)),
   & \hbox{if } s\in\mathop{\cup}\limits_j \Delta_j.
\end{array}
\right.
$$
%and formally set $\frac{|\bar\bv_k^\delta|}{|\bv_k|}=1$ if $|\bv_k|=0$.
 Then $\bar\bv_k^\delta$
satisfies the equation
\begin{equation}\label{deltaeq}
\bar\bv_k^\delta(\tau)=\bv_k(0)+\int_0^\tau \widehat {\mathbf{P}}_k(\bar\bv^\delta_k(s), \bv(s),s)\,ds
+\int_0^\tau \sum_{i=1}^{n_1}\widehat{B}_{ki}(\bar\bv^\delta_k(s),\bv(s),s)\,d\bbeta_i(s).
\ee
Notice that under Assumption \ref{assumption-v-3} the diffusion coefficient in this Ito equation does not degenerate.

{\bf Step 2: } Truncation at a level $\|v\|=R$.

Define the stopping time $\tau_R=\inf\{\tau\in[0,T]:\|v\|\geqslant R\}$. We define the processes $\bv_k^R$ equal
 to $\bv_k$ for $\tau\in[0,\tau_R]$ and satisfying the trivial equation $d\bv_k^R(\tau)=d\bbeta(\tau)$ for $\tau \in[\tau_R, T]$.
 We also set $\bar\bv_k^{\delta,R}$ to be equal to $\bar\bv_k^\delta$ for $\tau\in[0,\tau_R]$ and for $\tau>\tau_R$
 equal to a solution of the equation $d\bar\bv_k^{\delta,R}=U(\bv_k(\tau_R),\bar\bv_k^{\delta,R}(\tau_R))d\bbeta(\tau)$. Clearly, $\|\bar\bv_k^{\delta,R}\|\equiv\|\bv_k^R\|$. By (4) of Assumption \ref{assumption-v-3} we have $\mathbf{P}\{\bv_k^{R}(\tau)\neq\bv_k(\tau) \text{ for some } \tau\in[0,T]\}\to0$ as $R\to\infty$. Therefore, it is sufficient  to prove the lemma for $\bv_k$ replaced by $\bv_k^R$ with arbitrary $R>0$.

{\bf Step 3:} Taking a limit  as $\delta\to0$.

Now fix $R>0$. We denote $\bar\bv_k^{\delta, R}(\cdot)$ still as $\bar\bv_k^\delta(\cdot)$.

Argue as in  Lemma~\ref{l_compa} we find  that  the family of processes $\bar\bv_k^{\delta,R}(\cdot)$, $\delta\in(0,1]$
 is tight in $C([0,T],\mathbb R^2)$.
Therefore, for a subsequence $\delta_j\to0$, $\cD(\bar\bv_k^{\delta_j}(\cdot))$ converges to some
 $\mathcal{Q}^0_k\in\mathcal{P}(C([0,T],\mathbb{R}^2))$. Consider the processes
$$
\mathcal{Y}_k^\delta(\tau):=\int_0^\tau \widehat {\mathbf{P}}_k(\bar\bv^\delta_k(s), \bv(s),s)\,ds, \quad
 \mathcal{M}_k^\delta(\tau):=\sum_{i=1}^{n_1}\int_0^\tau \widehat{B}_{ki}(\bar\bv^\delta_k(s),v(s),s)\,d\bbeta_i(s),
$$
with an obvious change for $\tau\geqslant \tau_R$.
 The sequence of  pairs
$(\mathcal{Y}_k^{\delta_j}(\tau),\mathcal{M}_k^{\delta_j}(\tau))$ is tight in $C(0,T;\mathbb R^4)$.
If $(\mathcal{Y}_k^0(\cdot),\mathcal{M}_k^0(\cdot))$ is a limiting in law
 process as $\delta_j\to0$, then $\cD \bv_k(\cdot) = \mathcal{Q}^0_k$, where
  $$
  \bar\bv_k(\tau)=\mathcal{Y}^0_k(\tau)+\mathcal{M}_k^0(\tau), \quad\tau\in[0,T].
  $$

Denote $C_{R,k}^P=\sup\{|{\mathbf{P}}_k(v)|,\|v\|\leqslant R\} $ and $C_{R,k}^W=\sup\{|W_k(I(v))|,\|v\|\leqslant R\} $.
Then for any $0\leq \tau'<\tau''\leq T$ we have
$$
|\mathcal{Y}^\delta_k(\tau'')-\mathcal{Y}^\delta_k(\tau')|\leq C_{R,k}^P|\tau''-\tau'|+
C_{R,k}^W\eps^{-1}\Big|\bigcup\limits_j [\tau_j^-,\tau_j^+]\cap [\tau',\tau'']\Big|.
$$
From the definition of $\tau_j^\pm$ and $\tau_j^+$ it follows that
$$
\mathbf{E}\Big|\bigcup\limits_j {\Delta_j}\cap [\tau',\tau'']\Big|\leq \mathbf{E}\Big|\bigcup\limits_j {\Delta_j}\cap [0,T]\Big|
\leq \mathbf{E}\int\limits_0^T \mathbf{1}_{\{|\bv_k(\tau)|\leq2\delta\}}\,d\tau.
$$
By Theorem~2.2.4 in \cite{Kry77} for each $\eps>0$ the term on the right-hand side of this inequality tends to zero as $\delta\to0$.
Therefore,
$\mathbf{E}\Big|\bigcup\limits_j {\Delta_j}\cap [\tau',\tau'']\Big|\to 0$. Using the fact that the set
$\{\phi\in C(0,T;\mathbb R^2)\,;\, |\phi(\tau')-\phi(\tau'')|\leq 2C_{R,k}^P|\tau'-\tau''|\}$ is closed we derive from the convergence
$
\cD(\mathcal{Y}_k^{\delta_j}(\cdot)) \strela \cD(\mathcal{Y}_k^{0}(\cdot))
$
 that
$$
\mathbf{P}\big\{|\mathcal{Y}_k^0(\tau')-\mathcal{Y}_k^0(\tau'')|\leq C_{P,k}|\tau'-\tau''|\big\}=1.
$$
So $\mathcal{Y}_k^0(\tau)=\int_0^\tau \rho_k^0(s)\,ds$ with $|\rho(s)|\leq C_{R,k}^P$.

The processes $\mathcal{M}_k^\delta$ are continuous square integrable martingales  with respect to the natural filtration.
Since their second moment are  bounded uniformly in $\delta$, then
 the limit process $\mathcal{M}_k^0$ is also a square integrable martingale. Denoting by
  $\left\langle\mathcal{M}_k^\delta\right\rangle(\tau)$   the quadratic characteristic of $\mathcal{M}_k^\delta(\tau)$,
from Corollary VI.6.7 in \cite{JSh} we deduce that $\left\langle\mathcal{M}_k^\delta\right\rangle(\tau)\to
\left\langle\mathcal{M}_k^0\right\rangle(\tau)$ as $\delta\to0$.

Under Assumption \ref{assumption-v-3} the quadratic characteristic of $\mathcal{M}_k^\delta$ satisfies the estimates
$$
C_m(\tau''-\tau') |\zeta|^2\leq \big([\left\langle\mathcal{M}_k^\delta\right\rangle(\tau'')-\left\langle\mathcal{M}_k^\delta\right\rangle(\tau')])\zeta,\zeta\big)
\leq C_m^{-1}(\tau''-\tau') |\zeta|^2, \quad\zeta\in\mathbb R^2,
$$
for some constant $C_m>0$. Then the quadratic characteristic of   $\mathcal{M}_k^0$ also meets these estimates. Thus
 there exists
a progressively measurable random matrix function $\sigma_k(s)$ with values in the space of symmetric $2\times2$ matrices such that
$C_m|\zeta|^2\leq (\sigma_k(t)\zeta,\zeta)\leq C_m^{-1}|\zeta|^2$ for all $\zeta\in\mathbb R^2$, and
$\mathcal{M}_k^0(\tau)=\int_0^\tau \sigma^{\frac12}(s)d B_s$, where $B_s$ is a standard Wiener process in $\mathbb R^2$.
Therefore, the process $\bar\bv_k$ admits the following representation:
\be\label{append}
\bar\bv_k(\tau)=\bar\bv_k(0)+\int\limits_0^\tau \rho_k^0(s)\,ds+\int\limits_0^\tau \sigma_k^{\frac12}(s)d B_s
\ee
The desired statement is now an immediate consequence of Theorem~2.2.4 in \cite{Kry77}.

\section{On Birkhoff integrability of Hamiltonian systems with one degree of freedom}\label{app_B}

Here we prove a global version of the well known fact that 1d Hamiltonian systems are integrable (e.g. see Example~3 in \cite[Section~5.3]{AKN}).
Consider the plane $\mathbb{R}^2=\{\mathbf{x}=(x, y)\}$, equipped with the standard area-form $dx\wedge dy$. By $B_r$
 we denote the open disc $\{\|\mathbf{x}\|<r\}$, $r>0$,  and set ${I}(\mathbf{x})=\frac{1}{2}\|\mathbf{x}\|^2=\frac{1}{2}(x^2+y^2)$.

\begin{theorem}\label{birkhoff}
Assume  that $H\in C^\infty(\mathbb{R}^2)$ satisfies the following:

i)  $dH(0)=0$ and $d^2H(0)$ is positively definite;

ii) for each $\mathbf{x}\neq0$, $dH(\mathbf{x})\neq0$;

iii) for each $a\in H(\mathbb{R}^2\setminus\{0\})$  the level set $M_a=\{\mathbf{x}\in\mathbb{R}^2: H(\mathbf{x})=a\}$ is a connected   loop.

Then there exist a smooth canonical  change of coordinates (SCCC) $\Psi: \mathbb{R}^2\to\mathbb{R}^2$, $\Psi(0)=0$ and a smooth function $h$,  $h'(0)\neq0$, such that $H(\mathbf{x})=h\big({I}(\Psi(\mathbf{x}))\big)$.
\end{theorem}
\begin{proof}
{\bf Step 1:} By Vey's Theorem (see \cite{Eliasson} and see a 1d version of the theorem
 in \cite[Appendix~D]{KKPS}), there exist $\delta \in [0, 1/2]$ and a SCCC $Q_\delta: B_\delta\to\mathbb{R}^2$, $(p,q)\mapsto (x,y)$,
  such that $Q_\delta(0)=0$ and
$H\circ Q_\delta(p,q)=f(\frac{p^2+q^2}{2})$, where $f$ is a smooth function, satisfying $f'(0)\neq0$.

{\bf Step 2: } Now we   construct a SCCC $Q$ defined on the whole plane, such that for some $0<\delta'<\delta$ we have
 $Q|_{B_{\delta'}}=Q_\delta|_{B_{\delta'}}$ and
$Q|_{\mathbb{R}^2\setminus B_1}=L$, where $L$ is the linear symplectic transformation $L=dQ(0)=dQ_\delta(0)$.

Let $\bar Q_\delta=L^{-1}\circ Q_\delta$, $(p,q) \mapsto (x,y)$. Then
\begin{equation}\label{qxnot0}d\bar Q_\delta(0)=Id\;\text{ and }\;\partial_px(p,q)>0,
\end{equation}
if $\|(p,q)\|$ is small. Then,  decreasing $\delta$ if needed, we achieve that the transformation $\bar Q_\delta$ admits a smooth generating function $S(x,q)$,
so $\bar Q_\delta(p,q) = (x,y)$ if and only if $p= \p_qS, y=\p_xS$ (e.g. see \cite[Section~1.3]{AKN}). Since $d\bar Q_\delta(0) = \,$id, then
$
S(x,q)= xq + o(\| (x,q)\|^2).
$
Now we extend $S(x,q)$ from a small neighbourhood of the origin to the whole $(x,q)$-plane in such a way that $S(x,q)=xq$  for $\|(x,q)\|\geqslant1$, keeping
the  condition $\frac{\partial^2 S}{\partial x\partial q}>0$ for $\|(x,q)\|\leqslant1$. The extended $S$ is as a generating function of
a SCCC $Q': \mathbb{R}^2\to\mathbb{R}^2$, $(p,q)\mapsto(x,y)$ . Then  $Q'|_{\mathbb{R}^2\setminus B_{1}}=Id$  and
 $Q'|_{B_{\delta'}}=\bar Q_\delta|_{B_{\delta'}}$ for a small enough $\delta'<\delta$.
The required SCCC is obtained as $Q=L\circ Q'$.

{\bf Step 3: } Denote  $H_1:=H\circ Q(p,q)$. Clearly, conditions i)-iii) stay true for $H_1$ and in $B_{\delta'}$,
\begin{equation}\label{h1-f}
H_1(p,q)=f(\tfrac{p^2+q^2}2 ).
\end{equation}

Let $a_0\in H_1(\mathbb{R}^2\setminus\{0\})$. By  the Loiuville--Arnold theorem in a small neighbourhood   of the curve $M_{a_0} =\{ H_1 =a_0\}$
   exists a SCCC $Q_1$: $(p,q)\mapsto (I,\varphi)\in \mathbb{R}\times \mathbb{T}$, $\T = \R/2\pi$,
    such that $dp\wedge dq =dI\wedge d\varphi$ and $H_1(p,q)=h(I(p,q))$ for a smooth function $h$.
Moreover, $M_a\cong\{I=I(a)\}\times\mathbb{T}$ and $I(a)=\frac{1}{2\pi}\oint_{M_a}pdq$. By Green's formula, \begin{equation}\label{action-area}I(a)=\frac{1}{2\pi}\int_{M_a}pdq=\frac{1}{2\pi}\iint_{D(M_a)}dp\wedge dp,
\end{equation}
 where $D(M_a)$ is the domain enclosed by $M_a$.  Hence for each $a\in H_1(\mathbb{R}^2\setminus\{0\})$, $2\pi I(a)$ is the area enclosed by $M_a$.
 So the  action variable  $I(p,q)=I(a(p,q))$ is well defined  globally on $\mathbb{R}^2\setminus\{0\}$.

 The angle variable $\f\in\T$ is defined modulo a shift $\f \mapsto \f+ \zeta(I)$,  where $\zeta$ is any smooth function. To specify its choice we find a smooth
 curve $l_0\subset\R^2$  from the origin to infinity, such that $l_0\cap B_{\delta'}\subset\{p>0,q=0\}$ and $l_0$ insects
  each level set $M_a$ in exactly  one point.  Setting $\varphi\mid_{l_0}=0$ we get a uniquely defined  smooth angle variable   $\varphi$ on $\mathbb{R}^2\setminus \{0\}$.  The constructed  variables $(I,\f)$  define an action-angle transformation    $\R^2\setminus \{0\} \to \R_{>0}\times \T$,
  $(p,q) \mapsto (I,\f)$.

  Now let $G(p,q)=\sqrt{2I}(\cos\varphi,\sin\varphi)$. We then have a SCCC $G$:  $\mathbb{R}^2\setminus \{0\}\rcirclearrowleft$ such that
   $H_1(p,q)=h({I}(G(p,q)))$.
From \eqref{h1-f}, \eqref{action-area} and the normalisation $\f(p,0) = 0$ if $0<p<\delta'$
we have that  $G(p,q)=(p,q)$ on $B_{\delta'}\setminus \{0\}$.  Therefore $G$ extends  $Q\mid_{B_{\delta'}}$  to a SCCC of  $\mathbb{R}^2$.

We define the wanted SCCC as $\Psi=G\circ Q^{-1}$
and the assertion of the theorem follows.
\end{proof}

\noindent
{\bf \large Acknowledgement} We thank Anatoli Neishtadt for discussions of the deterministic averaging.
This research was  supported  by the Ministry of Science and Higher Education of the Russian Federation (megagrant No. 075-15-2022-1115). GH was also supported by National Natural Science Foundation of China (No. 20221300605).
%supported by BFS/TFS project "Pure Mathematics in Norway".

\bibliography{reference}{}
\bibliographystyle{plain}

\end{document}